\documentclass[11pt]{article}

\usepackage{calc}
\usepackage{color}
\usepackage{graphicx}
\usepackage{float}
\usepackage{placeins}
\usepackage{url}

\usepackage{booktabs}
\usepackage{tabularx}
\usepackage{multirow}
\usepackage{makecell}
\usepackage{threeparttable}
\usepackage{wrapfig}
\usepackage[export]{adjustbox}

\usepackage{mathtools}
\usepackage{bm}
\usepackage[bbgreekl]{mathbbol}
\usepackage{siunitx}

\usepackage{algorithmic}
\usepackage{algorithm}

\usepackage[superscript,biblabel]{cite}
\usepackage{titleref}

\usepackage[left=0.5in,right=0.5in,top=0.5in,bottom=0.5in,includefoot,heightrounded]{geometry}
\usepackage{lscape}
\usepackage[margin=15pt,font=small,labelfont={bf,sf},justification=justified]{caption}
\usepackage{subcaption} 
\usepackage{titling}
\usepackage{authblk}
\usepackage[pagewise]{lineno}
\usepackage{soul}
\usepackage{blindtext}
\usepackage{lipsum}
\usepackage{easylist}
\usepackage{tikz}

\usepackage{hyperref}
\usepackage{cleveref}

\newif\ifbembo
\newif\ifcharter
\newif\iferewhon
\newif\iflibertine
\newif\iflibertinealt
\newif\ifpalantino
\newif\iftimesnewroman

\bembofalse
\chartertrue
\erewhonfalse
\libertinefalse
\libertinealtfalse
\palantinofalse
\timesnewromanfalse

\ifbembo
  \usepackage[p,osf]{ETbb} 
  \usepackage[scaled=.95,type1]{cabin} 
  \usepackage[varqu,varl]{zi4}
  \usepackage[T1]{fontenc}
  \usepackage[libertine,vvarbb]{newtxmath}
  \usepackage[cal=boondoxo,bb=boondox]{mathalfa}
\fi

\ifcharter
  \usepackage[scaled=.98,p]{XCharter}
  \usepackage[scaled=1.04,varqu,varl]{inconsolata}
  \usepackage[type1]{cabin}
  \usepackage[uprightscript,xcharter,vvarbb,scaled=1.05]{newtxmath}
  \usepackage[cal=euler,scr=boondoxo,bb=boondox]{mathalpha}  
\fi

\iferewhon
  \usepackage[p,osf,scaled=.98,space]{erewhon} 
  \usepackage[varqu,varl]{inconsolata} 
  \usepackage[type1,scaled=.95]{cabin} 
  \usepackage[utopia,vvarbb]{newtxmath}
\fi

\iflibertine
  \usepackage[sb]{libertinus}
  \usepackage[T1]{fontenc}
  \usepackage{textcomp}
  \usepackage[varqu,varl]{zi4}
  \usepackage{libertinust1math} 
  \usepackage[scr=boondoxo,bb=boondox]{mathalpha} 
\fi

\iflibertinealt
  \usepackage[sb]{libertinus}
  \usepackage[T1]{fontenc}
  \usepackage{textcomp}
  \usepackage{libertinust1math} 
  \usepackage[scr=boondoxo,bb=boondox]{mathalpha} 
\fi

\ifpalantino
  \usepackage[largesc]{newpxtext}
  \usepackage[vvarbb]{newpxmath}
\fi

\iftimesnewroman
  \usepackage{newtxtext} 
  \usepackage[scaled=0.95]{inconsolata} 
  \usepackage[vvarbb]{newtxmath} 
\fi

\usepackage{microtype}
\pretitle{\begin{center}\bfseries\Large}
\posttitle{\end{center}}

\predate{\begin{center}\normalsize}
\postdate{\end{center}}

\pretitle{\begin{center}\bfseries\sffamily\LARGE}
\posttitle{\end{center}}
\usepackage{abstract}

\allowdisplaybreaks
\usepackage{sectsty} 
\allsectionsfont{\sffamily}

\makeatletter
\patchcmd{\LS@rot}{90}{-90}{}{}
\patchcmd{\endlandscape}{90}{-90}{}{}
\makeatother

\def \grad{\nabla}

\def \p{\partial}

\newcommand{\D}[2]{\frac{\p #1}{\p #2}}

\renewcommand{\vec}[1]{\bm{\mathrm{#1}}}

\def \CC{\mathbb{C}}

\def \FF{\mathbb{F}}

\def \PP{\mathbb{P}}
\def \PPs{\PP^{\mathrm s}}

\def \sigmaf{\bbsigma^{\mathrm f}}

\def \sigmae{\bbsigma^{\mathrm e}}

\def \F{\vec{F}}

\def \X{\vec{X}}
\def \e{\vec{e}}

\def \u{\vec{u}}

\def \x{\vec{x}}

\def \vchi{\vec{\chi}}

\def \Dx{{\mathrm d} \x}

\def \DX{{\mathrm d} \X}

\newcommand{\euleriandx}{{\ensuremath{\Delta x}}}

\def \mfac{\text{MFAC}}
\def \MFAC{\text{MFAC}}

\def \L2{L^2}

\newcommand{\trans}{^{\mkern-1.5mu\mathsf{T}}} 
\newcommand{\negtrans}{^{\mkern-1.5mu\mathsf{-T}}} 

\definecolor{vargreen}{rgb}{0.0, 0.5, 0.0}
\definecolor{varpurp}{rgb}{0.5, 0.0, 0.5}

\newcommand{\ub}{\vec{u}}

\newcommand{\xb}{\vec{x}}
\newcommand{\Xb}{\vec{X}}
\newcommand{\fb}{\vec{f}}
\newcommand{\Fb}{\vec{F}}

\newcommand{\Ub}{\vec{U}}

\newcommand{\Nb}{\vec{N}}
\newcommand{\nb}{\vec{n}}

\newcommand{\cauchy}{\bbsigma}

\newcommand{\cauchyv}{\bbsigma^{\text{v}}}
\newcommand{\cauchys}{\bbsigma^{\text{e}}}

\newcommand{\ztensor}{\mathbb{0}}
\newcommand{\pstab}{\pi_{\text{stab}}}

\newcommand{\soliddom}{{\Omega^{\text{s}}_t}}

\newcommand{\fluiddom}{\Omega^{\text{f}}_t}

\newcommand{\dev}{\text{dev}}

\newcommand{\tr}{\text{tr}}
\newcommand{\nus}{\nu_\text{stab}}
\newcommand{\kappas}{\kappa_\text{stab}}

\newcommand{\Qone}{\mathcal{Q}^1}

\title{Local Divergence-Free Immersed Finite Element-Difference Method Using Composite B-Splines}
\author[1]{Lianxia ~Li} 
\author[1]{Cole ~Gruninger} 
\author[1]{Jae H.~Lee}
\author[1-5$^{*}$]{Boyce E.~Griffith} 

\affil[1]{Department of Mathematics, University of North Carolina, Chapel Hill, NC, USA}
\affil[2]{Department of Biomedical Engineering, University of North Carolina, Chapel Hill, NC, USA}
\affil[3]{Carolina Center for Interdisciplinary Applied Mathematics, University of North Carolina, Chapel Hill, NC, USA}
\affil[4]{Computational Medicine Program, University of North Carolina School of Medicine, Chapel Hill, NC, USA}
\affil[5]{McAllister Heart Institute, University of North Carolina School of Medicine, Chapel Hill, NC, USA}
\affil[*]{\texttt{Corresponding author, boyceg@email.unc.edu}}

\date{} 

\begin{document}

\maketitle

\begin{abstract}

In the class of immersed boundary (IB) methods, the choice of the regularized delta function plays a crucial role in transferring information between fluid and solid domains through interpolation and spreading operators. Most prior work using the IB method has used isotropic kernels that do not preserve the divergence-free condition of the velocity field, leading to loss of incompressibility of the solid when interpolating the Eulerian velocity to Lagrangian markers. To address this issue, in simulations involving large deformations of incompressible hyperelastic structures immersed in fluid, researchers often use a volumetric stabilization approach such as adding a volumetric energy term and using modified invariants in the constitutive model of the immersed structure. Composite B-spline (CBS) kernels offer an alternative approach by inherently maintaining the discrete divergence-free property.  This work evaluates the performance of CBS kernels in terms of their volume conservation and accuracy, comparing them with several traditional isotropic kernel functions using a construction introduced by Peskin (referred to as IB kernels) and B-spline (BS) kernels. Benchmark tests include pressure-loaded and shear-dominated flows, such as an elastic band under differential pressure loads, a pressurized membrane, a compressed block, Cook's membrane, a slanted channel flow, and a modified Turek-Hron problem. Additionally, we validate our methodology using a complex fluid-structure interaction model of bioprosthetic heart valve dynamics in a pulse duplicator. Results demonstrate that CBS kernels achieve superior volume conservation  compared to conventional isotropic kernels, eliminating the need for additional volumetric stabilization techniques typically required to address instabilities arising from volume conservation errors.
Further, CBS kernels achieve convergence on coarser fluid grids, while IB and BS kernels need finer grids for comparable accuracy. Unlike IB and BS kernels, which perform better with larger mesh ratio factor between solid and fluid grids, CBS kernels show improved results with smaller mesh ratio factors. Additionally, the study reveals that although wider kernels provide more accurate results across all methods, CBS kernels are less sensitive to variations in relative grid spacings than isotropic kernels. This research highlights the advantages of CBS kernels in achieving stable, accurate, and efficient FSI simulations without requiring specialized volumetric stabilization treatments when simulating large deformations of elastic solids immersed in fluid.

\end{abstract}

\noindent \textbf{Keywords:} immersed boundary, fluid structure interaction, composite B-spline kernels (CBS), immersed finite element/finite difference method, isotropic kernel,  volume conservation, volumetric stabilization

\newpage

\section{Introduction}
\label{sec:introduction}
The immersed boundary (IB) method provides both a mathematical formulation and numerical framework for modeling fluid-structure interaction (FSI)\cite{peskin2002}. Originally introduced by Peskin to study blood flow through heart valves\cite{peskin1972}, the IB method provides a hybrid treatment of the fluid-structure system: the fluid is described in an Eulerian frame of reference while the immersed structure is tracked in a Lagrangian frame. In its continuous formulation, the fluid-structure coupling is mediated by integral operators with singular delta function kernels supported at the immersed structure's location. The numerical implementation of the IB method approximates these singular delta functions with regularized versions and discretizes the Eulerian and Lagrangian degrees of freedom separately. The Eulerian degrees of freedom are solved for on a fixed, typically structured grid, while the Lagrangian degrees of freedom are discretized using either Lagrangian marker points (as in Peskin's classical IB method\cite{peskin1972,peskin2002}) or finite elements\cite{Griffith2017}.The fluid-structure coupling of the discrete formulation is achieved by using the regularized delta function to spread Lagrangian force densities to nearby Eulerian grid points and to interpolate fluid velocities back to Lagrangian markers so that the marker moves with the local fluid velocity. 

Since its introduction, the IB method has emerged as a powerful tool for modeling biological FSI problems, with applications spanning from animal locomotion \cite{jones2015,santhanakrishnan2018,alben2013,bhalla2014,tytell2014,bale2015,hoover2017,nangia2017} to cardiac mechanics \cite{griffith_luo_2009,hasan2017,griffith_heart2012,chen2016,crowl2011,kaiser2019,davey2023}. Because the IB method discretizes the structural and Eulerian degrees of freedom separately, it avoids the difficulties of grid generation associated with body-fitted approaches. This makes the method particularly competitive for problems involving large deformations, where body-fitted approaches require computationally intensive remeshing procedures as the structure evolves from its initial configuration.

While Peskin's classical IB method is effective for modeling fiber-based elasticity in which the material behaves as a collection of one-dimensional elastic arrays immersed in fluid, it has limitations when applied to thick bodies governed by general hyperelastic constitutive laws. To address this limitation, several researchers have replaced the Lagrangian marker description with finite element (FE) treatments of the structural equations, which provide a natural framework for modeling thick immersed bodies with general structural constitutive laws. Among the earliest works in this direction were those of Zhang et al.\cite{zhang2004} and Wang and Liu\cite{Wang2004}, who coupled fluid and solid degrees of freedom using regularized delta functions adapted from reproducing kernel particle methods\cite{Liu1995}. Taking a different approach, Boffi et al.\cite{Boffi2008} developed a fully variational FE-based formulation that treated the momentum equations in weak form, eliminating the need for regularized delta functions altogether. Griffith and Luo \cite{griffith_luo_2009} introduced a hybrid approach that combines FE discretization of the Lagrangian structure with a staggered Cartesian grid finite difference scheme for Eulerian degrees of freedom. Their method treats the Lagrangian momentum equations in weak form while maintaining the strong form of the Eulerian equations, thus requiring regularized delta functions for force spreading and velocity interpolation. In our work, we employ this hybrid approach of Griffith and Luo, which we refer to as the immersed finite element/finite difference method (IFED) due to its combination of FE and finite difference discretizations\cite{Griffith2017,wells2023nodal,vadalaroth2020,davey2023,lee2022lagrangian}.

A critical aspect of the IFED method, and indeed many IB formulations, is the choice of kernel function used to construct the regularized delta function. This choice significantly impacts the method's accuracy, particularly for relatively coarse discretizations\cite{lee2022lagrangian}. In IB formulations of fluid-structure interaction, these kernels are typically constructed as tensor products of one-dimensional functions\cite{peskin2002}. Most commonly, these kernels are isotropic, meaning that identical one-dimensional kernels are used in each coordinate direction to form the complete kernel function.

Despite their widespread success in both IB and IFED methods, isotropic regularized delta functions can lead to gradual volume loss in structures under pressure loads\cite{peskin1993,gruninger2024local,cortez2000,peskin2002,bao2017,griffith2012,long2003}. This behavior may seem counterintuitive since the continuous equations of motion guarantee, via the Reynolds transport theorem, that the immersed body should remain incompressible because the Eulerian velocity field that determines the body's motion is constructed to be divergence free. Two major factors contribute to this lack of volume conservation. First, although the velocity field may be discretely divergence-free (meaning its finite difference approximation of the divergence operator is zero), the interpolated velocity used to move the structure is not guaranteed to be continuously divergence-free. This limitation, combined with the fact that the interpolated velocity field is only first-order accurate near interfaces, can lead to substantial volume changes over the course of a simulation. Second, even when the Lagrangian force density should yield an Eulerian force density that is the gradient of a scalar field, the spreading operation does not typically preserve this gradient structure at the discrete level. As a result, the pressure, which can only balance discrete gradient forces, cannot fully compensate for all compressive forces spread onto the background Cartesian grid. These unbalanced forces then inappropriately influence the discrete fluid velocity field, leading to spurious flows normal to the surface of the immersed body that result in volume changes. While time stepping errors can also contribute to spurious volume changes, these effects are typically an order of magnitude smaller than the errors arising from the force spreading and interpolation operations\cite{gruninger2024local}.

In the context of the IFED method, Vadala-Roth et al.\cite{vadalaroth2020} addressed spurious volume changes by adopting approaches from computational mechanics used to model nearly incompressible hyperelastic structures. Their method separates the elastic strain energy functional into deviatoric and volumetric components and expresses it using modified invariants of the right Cauchy-Green tensor\cite{bonet2008}. While effective at preserving volume, this approach inatroduces significant computational challenges. The modified invariants increase the nonlinearity of the stress response, which can complicate implicit solvers. Additionally, the volumetric penalty terms are known to impose more severe time step restrictions for explicit timestepping schemes\cite{Devendran2012}. Moreover, these modifications not only affect the stress prefactor but also introduce additional isotropic stresses that alter the pressure response\cite{Devendran2012, vadalaroth2020}.

Recently, Gruninger and Griffith\cite{gruninger2024local} introduced regularized delta functions based on composite B-splines to mitigate volume conservation errors in Peskin's original IB method formulation. The key property of composite B-splines is their ability to transform finite differences into continuous derivatives, a property that follows from the convolution identity that generates B-spline family\cite{unser1992,schroeder2022,handscomb1984}. On staggered MAC grids, this property ensures that composite B-splines interpolate discretely divergence-free vector fields to continuously divergence-free vector fields\cite{schroeder2022,handscomb1984,gruninger2024local}. Gruninger and Griffith demonstrated that in the context of Peskin's classical IB method sufficiently regular composite B-splines can control volume conservation errors to the magnitude of those introduced due to discrete timestepping.

This study extends the work of Gruninger and Griffith \cite{gruninger2024local} by examining CBS regularized delta functions in the context of the IFED method. We demonstrate that these kernels' superior volume conservation properties often eliminate the need for both modified invariants and artificial (numerical) bulk modulus traditionally required for modeling incompressible materials. Through a suite fluid-structure interaction benchmarks, we show that while these kernels significantly improve volume conservation in pressure-loaded cases, they also yield results consistent with established studies for shear-dominated flows, including flow through an inclined channel and the Turek-Hron benchmark. Additionally, we evaluate these composite kernels in a complex fluid-structure interaction model of bioprosthetic heart valve dynamics in a pulse duplicator. 



\section{Continuous Equations of Motion}
\label{sec:cts_eqs_of_motion}
The continuous formulation of fluid-structure interaction within the immersed boundary framework considers the coupled dynamics of a fluid and an immersed structure in a fixed domain $\Omega$. At each time $t$, the domain is partitioned into fluid and solid subdomains, $\Omega^\text{f}_t$ and $\Omega^\text{s}_t$ respectively, such that $\Omega = \Omega^\text{f}_t \cup \Omega^\text{s}_t$. 

The IB framework employs both Eulerian and Lagrangian descriptions. Physical points in the domain are described by coordinates $\xb = (x_1, x_2, x_3) \in \Omega$, and material points of the structure are labeled by reference coordinates $\vec{X} = (X_1, X_2, X_3) \in \Omega^\text{s}_0$. The mapping $\vchi(\vec{X},t) \in \Omega^\text{s}_t$ tracks the physical position of material point $\vec{X}$ at time $t$, with $\Nb(\vec{X})$ denoting the outward unit normal to $\partial\Omega^\text{s}_0$ at $\vec{X}$. The Eulerian velocity and pressure fields are denoted using $\vec{u}(\xb,t)$ and $p(\xb,t)$, respectively. The Eulerian force density arising from the presence of the immersed structure is denoted $\fb(\xb,t)$. The structure is characterized by two Lagrangian quantities: the force density $\vec{F}(\vec{X},t)$ generated by hyperelastic stresses and the velocity field $\Ub(\vec{X},t)$. The material derivative is denoted by $\frac{\text{D}}{\text{D}t} = \frac{\p}{\p t} + \vec{u} \cdot \nabla$. Throughout, we assume both the fluid and immersed viscoelastic structure are of the same uniform density $\rho$ and viscosity $\mu$. 

The fluid-structure interaction is encoded mathematically using the following coupled system of equations
\begin{align}
    \label{eq:momentum}
    \rho\frac{{\mathrm D}\u}{{\mathrm D}t}(\xb,t) &= - \grad p(\xb,t) + \mu \grad^2 \u(\xb,t) + \fb(\xb,t), \\
    \label{eq:continuity}
    \grad \cdot \u(\xb,t) &= 0, \\
    \label{eq:fsiconstraint}
    \fb(\xb,t) &= \int_{\Omega^s_0}\vec{F}(\vec{X},t)\,\delta(\xb - \vchi(\vec{X},t)) \, \DX,  \\
    \label{eq:noslip}
    \D{\vchi}{t}(\vec{X},t) &= \Ub(\vec{X},t) = \int_\Omega \u(\xb,t) \, \delta(\xb - \vchi(\vec{X},t)) \, \Dx = \u(\vchi(\vec{X},t),t).
\end{align} 
The Dirac delta function $\delta(\xb)$ facilitates the interaction between Eulerian and Lagrangian descriptions of the fluid and structure, respectively. Through equation~\eqref{eq:fsiconstraint}, the delta function spreads Lagrangian force densities from the immersed structure to the Eulerian fluid, while in Equation~\eqref{eq:noslip}, the delta function evaluates the fluid velocity at the Lagrangian structure to enforce the no-slip condition ensuring that the fluid and structure move in concert. Since the fluid is incompressible, the Reynolds transport theorem guarantees that the immersed structure is also incompressible.

\subsection{Immersed Elastic Bodies}
For immersed bodies described as hyperelastic structures the total Cauchy stress tensor $\bbsigma$ of the fluid-structure system is spatially decomposed as
\begin{equation}
\label{eq:stress_decomp}
\bbsigma\left(\vec{x},t\right) = \begin{cases} 
    \sigmaf\left(\vec{x},t\right) & \text{in }\Omega\setminus\Omega^{\text{s}}_t ,\\
    \sigmaf\left(\vec{x},t\right) + \sigmae\left(\vec{x},t\right) & \text{in }\Omega^{\text{s}}_t,
\end{cases}
\end{equation}
in which $\sigmaf = -p\mathbf{I} + \mu(\nabla\vec{u} + \nabla\vec{u}\trans)$ is the viscous stress tensor of the Newtonian fluid and $\sigmae$ is the elastic Cauchy stress tensor of the immersed body.

The continuity of total stress $\bbsigma$ throughout $\Omega$, which follows from Newton's third law of motion, induces a jump condition in the fluid stress tensor across the fluid-solid interface $\partial\Omega^{\text{s}}_t$
\begin{equation}
\label{eq:jump_condition}
 [\![\sigmaf]\!](\xb,t)\nb(\xb,t) = \sigmae(\xb,t)\nb(\xb,t), \;\; \xb \in \partial\Omega^{\text{s}}_t
\end{equation}
with $[\![\sigmaf]\!](\xb,t) = \lim_{\epsilon\to 0^+} \sigmaf(\xb+\epsilon\nb,t)- \sigmaf(\xb-\epsilon\nb,t)$ and $\nb$ the unit outward normal to $\partial\Omega^{\text{s}}_t$.

When implementing the IB framework, it is advantageous to express the elastic stress using Lagrangian reference variables. This is achieved through the first Piola--Kirchhoff stress tensor $\PP^{\text{e}}$, which is related to the Cauchy stress through the kinematic transformation
\begin{equation}
\label{eq:kinematic_transform}
\PP^{\text{e}}(\vec{X},t) = J\left(\vec{X},t\right)\sigmae\left(\vchi(\vec{X},t),t\right)\FF\left(\vec{X},t\right)\negtrans,
\end{equation}
in which $\FF = \p\vchi/\p\X$ is the deformation gradient tensor and $J(\vec{X},t) = \text{det}\left(\FF(\vec{X},t)\right)$ is its determinant. For hyperelastic materials, $\PP^{\text{e}}$ can be derived from a strain energy density functional $\Psi$ via $\PP^{\text{e}} = \D{\Psi}{\FF}$. 

Following the approach of Boffi et al.~\cite{Boffi2008}, the coupled fluid-structure system can be expressed in distributional form, consistent with Peskin's original immersed boundary formulation
\begin{align}
\rho\frac{\mathrm{D}\u}{\mathrm{D}t}(\xb,t) &= -\nabla p(\xb,t) + \mu\nabla^2\vec{u}(\xb,t) + \fb(\xb,t), \\
\nabla\cdot\u(\xb,t) &= 0, \\
\label{eq:elastic_eul_force}
\fb^\text{e}(\xb,t) &= \int_{\Omega^{\text{s}}_0}\left(\nabla_{\vec{X}}\cdot\PP^{\text{e}}\right)\,\delta\left(\xb - \vchi\left(\vec{X},t\right)\right)\,\text{d}\vec{X} \nonumber \\ 
           &\quad - \int_{\p\Omega^{\text{s}}_0}\PP^{\text{e}}\Nb\,\delta\left(\xb - \vchi\left(\vec{X},t\right)\right)\,\text{d}A,\\
\label{eq:coupled_motion}
\frac{\p\vchi}{\p t}\left(\vec{X},t\right) &= \int_{\Omega} \u(\xb,t)\,\delta\left(\xb-\vchi\left(\vec{X},t\right)\right)\,\text{d}\xb.
\end{align}
The Eulerian elastic force density $\fb^e(\xb,t)$ comprises two distinct Lagrangian contributions: an internal volumetric force density $\nabla_{\vec{X}}\cdot\PP^{\text{e}}$ distributed throughout the elastic body and a surface transmission force density $-\PP^{\text{e}}\Nb$ concentrated on $\partial\Omega^{\text{s}}_0$. This surface force emerges from the stress jump condition \eqref{eq:jump_condition}, which ensures continuity of the total fluid-structure traction vector.

We re-emphasize that since the elastic structure moves according to the local Eulerian velocity through \eqref{eq:coupled_motion} and the Eulerian velocity field satisfies global incompressibility ($\nabla\cdot\u = 0$), all deformations of the elastic structure necessarily preserve volume, maintaining $J(\vec{X},t) = 1$. Consequently, the Eulerian pressure $p$ defined within the elastic body serves as a Lagrange multiplier enforcing incompressibility of the structural deformation. This renders any volumetric terms in the elastic strain energy density (and thus in the first Piola-Kirchhoff stress tensor $\PP^{\text{e}}$) technically redundant, as the incompressibility constraint is already enforced through the fluid pressure. 

To facilitate the finite element discretization of the Lagrangian variables as employed in the IFED method, we derive a variational formulation for computing the Lagrangian force density $\vec{F}(\vec{X},t)$. This force density determines the Eulerian force density through equation \eqref{eq:fsiconstraint} and arises from the volumetric and surface force contributions in equation \eqref{eq:elastic_eul_force}. The variational formulation is obtained by equating the elastic Eulerian force density with the Eulerian force density appearing in equation \eqref{eq:fsiconstraint}, multiplying them by a suitable test function $g(\xb)$ and integrating over $\Omega$. Upon performing these operations, we obtain the following variational characterization of the Lagrangian force density
\begin{align}
\label{eq:partitioned_form}
\int_{\Omega^{\text{s}}_0} \vec{F}(\vec{X},t)G(\vec{X})\,\text{d}\vec{X} &= \int_{\Omega^{\text{s}}_0} (\nabla_{\vec{X}}\cdot\PP^{\text{e}})G(\vec{X})\,\text{d}\vec{X} - \int_{\partial\Omega^{\text{s}}_0} (\PP^{\text{e}}\Nb)G(\vec{X})\,\text{d}A,
\end{align}
in which $G(\vec{X}) = \int_{\Omega}g(\xb)\,\delta\left(\xb-\vchi\left(\vec{X},t\right)\right)\,\Dx = g\left(\vec{\chi}(\vec{X},t)\right)$. By integrating equation \eqref{eq:partitioned_form} by parts, an equivalent variational formulation of the Lagrangian force density may be attained
\begin{align}
\label{eq:unified_formulation}
\int_{\Omega^{\text{s}}_0} \vec{F}(\vec{X},t)G(\vec{X}),\text{d}\vec{X} &= -\int_{\Omega^{\text{s}}_0} \PP^{\text{e}}(\vec{X},t)\nabla_{\vec{X}}G(\vec{X}),\text{d}\vec{X}.
\end{align}
The variational formulation of the Lagrangian force density appearing in equation \eqref{eq:unified_formulation} was characterized by Griffith and Luo as the \textit{unified} formulation. The unified formulation offers implementation advantages through its simpler structure and reduced regularity requirements on the deformation gradient and first Piola-Kirchhoff stress, making it particularly suitable for nodal finite element approximations. Most finite element-based IB methods have adopted this formulation for approximating the Lagrangian force density. In this work, we employ the unified formulation to evaluate the effectiveness of composite B-spline regularized delta functions. 

\subsection{Immersed Rigid Structures}
For stationary structures with prescribed velocity $\vec{V}(\vec{X},t)$ and position $\vec{\psi}(\vec{X},t)$, the Lagrangian force density acts as a Lagrange multiplier enforcing the kinematic constraints $\Ub(\vec{X},t) = \vec{V}(\vec{X},t)$ and $\vchi(\vec{X},t) = \vec{\psi}(\vec{X},t)$. A direct implementation of these constraints requires solving an extended saddle-point system at each timestep. Such systems, involving many degrees of freedom, present significant computational challenges due to the relative scarcity of effective preconditioners for iterative based solver approaches. 
To circumvent the difficulty associated with solving the extended saddle-point system, we employ a penalty formulation that approximates the Lagrange multiplier through via
\begin{equation}
\vec{F}(\X,t) = \kappa\left(\vec{\psi}(\vec{X},t) - \vchi(\vec{X},t)\right) + \eta\left(\vec{V}(\vec{X},t) - \Ub(\vec{X},t)\right).
\label{eq:tether_force}
\end{equation}
This penalty formulation can be interpreted physically as combining elastic restoring force with stiffness $\kappa\geq 0$ that penalizes positional deviations, and a viscous damping force with coefficient $\eta\geq 0$ that enforces the velocity constraint. The penalty formulation converges to the fully constrained Lagrange multiplier formulation as either $\kappa \to \infty$ or $\eta\to\infty$. Although either penalty term alone is sufficient to approximate the fully constrained formulation as its corresponding parameter approaches infinity, the combination of both terms often yields superior results for finite $\kappa$ and $\eta$. Empirically, it is observed that the viscous dampening force benefits the numerical stability and reduces the abundance of spurious force that arise when local Lagrangian degrees of freedom effectively over constrain their Eulerian counterparts. The elastic restoring force on the other hand provides a more robust enforcement of structural rigidity, particularly when the body or surface is subjected to substantial pressure loads. 

\section{Numerical Discretization}
The discretization of the equations of motion follows that of Grifftih and Luo~\cite{Griffith2017}. For simplicity, we describe the numerical scheme in two spatial dimensions. The extension of the description to three spatial dimensions is straightforward. 
\subsection{Spatial Discretization}
The Eulerian variables are discretized using a staggered-grid marker and cell approach~\cite{harlow1965}. The domain is partitioned into uniform grid cells of width $h>0$, with each cell indexed by a pair of integers $(i,j)$. The pressure $p_{i,j}$ is approximated at cell centers, located at $\xb_{i,j} = \left((i-1/2)h,(j-1/2)h\right)$. The horizontal velocity and force components $u_{i-1/2,j}$ and $\left(f_1\right)_{i-1/2,j}$ are approximated at the centers of vertical cell faces, located at $\xb_{i-1/2,j} = \left(ih,(j-1/2)h\right)$, and the vertical components $v_{i,j-1/2}$ and $\left(f_2\right)_{i,j-1/2}$ are approximated at the centers of horizontal cell faces, located at $\xb_{i,j-1/2} = \left((i-1/2)h,jh\right)$. The Laplacian of the velocity components and gradient of the pressure are approximated using standard second-order accurate central differences. The nonlinear convective term is approximated using the staggered grid version\cite{griffith2009accurate} of the xsPPM7\cite{rider2007} variant of the piecewise parabolic method\cite{colella1984}. In some of the numerical examples we showcase below, we a Cartesian grid adaptive mesh refinement (AMR) approach described by Griffith\cite{griffith_heart2012}.

The Lagrangian structure $\Omega_0^s$ is discretized using a triangulation $\mathcal{T}^h$ with $M$ nodes, and the Lagrangian variables are approximated using a nodal finite element approach. The nodes of the finite element mesh are denoted by $\{\Xb\}_{l=1}^M$ with corresponding shape functions $\{\phi_l \left(\Xb\right)\}_{l=1}^M$. The time-dependent physical positions of the mesh nodes are denoted by $\{\vchi_l\left(t\right)\}_{l=1}^M$. The discrete deformation $\vchi\left(\Xb,t\right)$ is described using the finite element shape functions:
\begin{equation}
\vchi_h \left(\Xb,t\right) = \sum_{l=1}^M\phi_l(\Xb)\vchi_l(t).
\end{equation}
The discrete approximation of the deformation gradient tensor is obtained by computing the total derivative of the discrete deformation
\begin{equation}
\FF_h\left(\Xb,t\right) = \frac{\partial \vchi_h}{\partial \Xb} = \sum_{l=1}^M \vchi_l(t)\frac{\partial \phi_l(\Xb)}{\partial \Xb}.
\end{equation}

Because $\vchi_h(\Xb,t)$ is approximated using nodal shape functions, which are $C^0$ but not $C^1$, $\FF_h$ is generally discontinuous across element boundaries. The approximation to the first Piola-Kirchhoff stress $\PP^{\text{e}}_h$ is obtained using $\FF_h$. The Lagrangian force density $\Fb(\Xb,t)$ is likewise approximated using the nodal finite element shape functions
\begin{equation}
\Fb_h(\Xb,t) = \sum_{l=1}^M \Fb_l(t)\phi_l\left(\Xb\right).
\end{equation}

For an immersed elastic body, the time-dependent nodal coefficients of the Lagrangian force densities are determined using a standard Galerkin projection:
\begin{equation}
\label{eq:disc_uni_formulation_elas}
\sum_{l=1}^M\left(\int_{\Omega_0^s}\phi_l(\Xb)\phi_m(\Xb)\,\text{d}\Xb\right)\Fb_l(t) = -\int_{\Omega_0^s}\PP^{\text{e}}\left(\Xb,t\right) \nabla_{\Xb}\phi_m\left(\Xb\right)\,\text{d}\Xb,
\end{equation}
for $m = 1,\ldots,M$. In practice, the integrals in \eqref{eq:disc_uni_formulation_elas} are evaluated by applying quadrature rules over individual elements. We employ the efficient nodal coupling scheme of Wells et al.\cite{wells2023nodal}, which yields the diagonal Galerkin projection
\begin{equation}
\Fb_m(t) = -\frac{1}{w_m}\sum_{q=1}^M\PP^{\text{e}}\left(\Xb_q,t\right) \nabla_{\Xb}\varphi_m\left(\Xb_q\right)w_q,
\end{equation}
in which $w_q$ represents nodal weights that can be defined using standard approaches for constructing diagonal or lumped mass matrices. For higher-order elements such as P2 elements, these weights are typically constructed by taking the diagonal of the elemental mass matrix and rescaling it to preserve the total mass of each element. Importantly, Wells et al. demonstrated that the specific numerical values of these weights cancel out when the same nodal quadrature scheme is used to approximate the force spreading operator associated with the IFED method\cite{wells2023nodal}.

For an immersed rigid structure, the discretized Lagrangian penalty force $\Fb_h(\Xb,t)$ is computed directly from the discrete deformation approximation $\vchi_h\left(\Xb,t\right)$:
\begin{equation}
\Fb_h\left(\Xb,t\right) = \kappa\left(\vec{\psi}\left(\Xb,t\right) - \vchi_h\left(\Xb,t\right)\right) + \eta\left(\vec{V}(\Xb,t) - \frac{\partial\vchi_h}{\partial t}(\Xb,t)\right).
\end{equation}
Because $\Fb_h$ and $\vchi_h$ share the same finite element basis, the nodal coefficients of $\Fb_h$ are determined directly by those of $\vchi_h$:
\begin{equation}
\Fb_l(t) = \kappa\left(\vec{\psi}_l\left(t\right) - \vchi_l\left(t\right)\right) + \eta\left(\vec{V}_l(t)  - \frac{\partial\vchi_l}{\partial t}(t)\right).
\end{equation}
Throughout the rest of the paper, we drop the subscript $h$ and assume each Lagrangian variable is represented by its discrete approximation to prevent the notation from becoming overly cumbersome.

\subsubsection{Lagrangian-Eulerian Coupling}
\label{subsec:lag_eul_coupling}

Following standard immersed boundary methodology, we replace the singular delta function coupling the Eulerian and Lagrangian degrees of freedom in equations~\eqref{eq:fsiconstraint} and~\eqref{eq:noslip} by a regularized version $\delta_h(\xb)$. The regularized delta functions are constructed as tensor products of one-dimensional regularized delta functions, which we detail in section \ref{subsec:choices_of_delta}. 

Using this regularized delta function, we approximate the Eulerian force density $\fb = \left(f_1,f_2\right)$ by discretizing the integral in equation \eqref{eq:fsiconstraint}. This discretization involves decomposing the integral into a sum over elements in the triangulation $\mathcal{T}^h$ and applying a quadrature scheme to each element. We use a nodal quadrature scheme as done by Wells et al.\cite{wells2023nodal}, which yields the following approximations for the components of the Eulerian force density:
\begin{align}
    \label{eq:force_spread}
\left(f_1\right)_{i-1/2,j} &= \sum_{l=1}^M (F_1)_l(t)\delta_h\left(\xb_{i-1/2,j}-\vchi_l(t)\right)w_l, \\
\left(f_2\right)_{i,j-1/2} &= \sum_{l=1}^M (F_2)_l(t)\delta_h\left(\xb_{i,j-1/2}-\vchi_l(t)\right)w_l,
\end{align}
in which $\Fb = (F_1,F_2)$ denotes the components of the discrete Lagrangian force density. For conciseness, we express the force spreading operation \eqref{eq:force_spread} using the compact notation
\begin{equation}
\fb(\x,t) = \boldsymbol{\mathcal{S}}[\vchi(\cdot,t)]\,\F(\X,t),
\end{equation}
and we refer to $\boldsymbol{\mathcal{S}}$ as the force spreading operator.

For the resulting semi-discrete formulation of the equations of motion to satisfy energy conservation\cite{peskin2002}, the operator mapping discrete Eulerian velocities to discrete Lagrangian operators must be discretely adjoint to the force spreading operator $\boldsymbol{\mathcal{S}}$. This relationship is expressed as
\begin{equation}
\Ub(\X, t) = \D{\vchi}{t}(\Xb,t) = \boldsymbol{\mathcal{S}}^{\ast}[\vchi(\cdot,t)]\,\u(\x,t).
\end{equation}
We define this velocity interpolation operator as $\boldsymbol{\mathcal{J}} = \boldsymbol{\mathcal{S}}^{\ast}$. Using the nodal quadrature scheme, each component $\Ub = \left(U,V\right)$ of the Lagrangian velocity can be written as
\begin{align}
U_l &= h^2\sum_{i,j}u_{i-1/2,j}\,\delta_h\left(\xb_{i-1/2,j} - \vchi_l(t)\right), \\
V_l &= h^2\sum_{i,j}v_{i,j-1/2}\,\delta_h\left(\xb_{i,j-1/2} - \vchi_l(t)\right).
\end{align}
As shown by Wells et al.\cite{wells2023nodal}, when using consistent nodal quadrature rules for both force computation and velocity interpolation, the interpolated velocity equals the projected interpolated velocity, eliminating the need for an additional projection step onto the finite element space.
\subsection{Choices of Regularized Delta Functions}
\label{subsec:choices_of_delta}
In this section, we detail the regularized delta functions used in our numerical experiments. These functions are constructed using tensor products of one-dimensional regularized delta functions, which themselves are built from dilated and scaled versions of a kernel function $\varphi(x)$. Following standard immersed boundary formulations, we choose the dilation and scale factors to match the background Cartesian mesh-width $h>0$. In one dimension, the regularized delta function takes the form
\begin{equation}
\delta_h(x) = \frac{1}{h}\varphi\left(\frac{x}{h}\right),
\end{equation}
and in two dimensions
\begin{equation}
\delta_h(\xb) = \frac{1}{h^2}\varphi\left(\frac{x}{h}\right)\zeta\left(\frac{y}{h}\right),
\end{equation}
in which $\zeta$ is another choice of one-dimensional kernel function. Here, we assume a uniform grid increment $h>0$ in both $x$ and $y$ directions. This formulation can be readily extended to higher dimensions.
\label{subsec:choices_reg_delta}
\label{subsec:choices_reg_delta}
\subsubsection{Isotropic Regularized Delta Functions}
Isotropic regularized delta functions have been most commonly used with immersed boundary methods. An isotropic regularized delta function uses the same one-dimensional kernel function $\varphi$ in each dimension of its tensor product formulation. For example, in two spatial dimensions:
\begin{equation}
\delta_h(\mathbf{\xb})=\frac1{h^2}\varphi\left(\frac xh\right)\varphi\left(\frac yh\right).
\end{equation}
The isotropic regularized delta functions we employ here are constructed from two families of kernel functions. The first family, which we call the ``IB family'', was originally postulated by Peskin based on constraints ensuring accuracy, physicality, and stability. The most widely used member of the IB family is the 4-point kernel function $\varphi_{\text{IB4}}$\cite{peskin2002}, given by:
\begin{equation}
    \varphi(r) = \begin{cases}
    0, & r \leq -2 \\
    \frac{1}{8}\left(5 + 2r - \sqrt{-7 - 12r - 4r^2}\right), & -2 \leq r \leq -1 \\
    \frac{1}{8}\left(3 + 2r + \sqrt{1 - 4r - 4r^2}\right), & -1 \leq r \leq 0 \\
    \frac{1}{8}\left(3 - 2r + \sqrt{1 + 4r - 4r^2}\right), & 0 \leq r \leq 1 \\
    \frac{1}{8}\left(5 - 2r - \sqrt{-7 + 12r - 4r^2}\right), & 1 \leq r \leq 2 \\
    0. & 2 \leq r
    \end{cases}
    \end{equation}
Setting $r = |x|$ ensures the kernel function is even. Although Peskin has derived many other kernels based on his original postulates\cite{bao2016,Roma1999}, we direct readers to the study by Lee and Griffith\cite{lee2022lagrangian} rather than providing a full review here. Going forward, we use the notation $\text{IB}_{\#}$ to indicate an IB kernel supported at $\#$ number of grid points. We specifically note that $\text{IB}_6$ refers to the three-times continuously differentiable kernel derived by Bao et al.\cite{bao2016}

The second family of isotropic regularized delta functions we consider is derived from the B-spline family of kernel functions. These isotropic B-spline regularized delta functions are constructed using dilated and scaled B-splines centered at zero. Each B-spline kernel corresponds to an order $n \ge 0$ and is a degree $n$ piecewise polynomial with $n-2$ continuous derivatives, compactly supported on the closed interval $\left[-\frac{1}{2}-\frac{n}{2},\frac{1}{2}+\frac{n}{2}\right]$.

B-spline piecewise polynomials have a rich history\cite{schoenberg1973}, though they are primarily attributed to Schoenberg, who developed them during his ballistics research in World War II\cite{schoenberg1946}. Later, Carl de Boor\cite{deBoor1978} popularized B-splines for computer-aided design (CAD) and formulated a recursive method for constructing order $n$ B-splines from the divided differences of order $n-1$ B-splines\cite{deBoor1978}. This formula has become widely known as the Cox-de Boor recursion formula.
Using this recursion formula, one can show that any order $n$ B-spline centered at zero can be generated by convolving the order $n-1$ B-spline with the zeroth order B-spline (which is simply the characteristic function of the interval $\left(-\frac{1}{2},\frac{1}{2}\right)$):
\begin{equation}
    \label{eq:conv_identity}
    \varphi^{\text{BS}}_n(x) = \varphi^{\text{BS}}_{n-1}\ast\varphi^{\text{BS}}_1 = \int_{-\frac{1}{2}}^{\frac{1}{2}} \varphi^{\text{BS}}_{n-1}(x-y)\,\text{d}y.
\end{equation}
Differentiating this formula, the fundamental theorem of calculus implies that the derivative of an order $n$ B-spline is the central difference of order $n-1$ B-splines
\begin{equation}
\frac{\text{d}\varphi^{\text{BS}}_n}{\text{d}x}(x) = \varphi^{\text{BS}}_{n-1}\left(x + \frac{1}{2}\right) - \varphi^{\text{BS}}_{n-1}\left(x - \frac{1}{2}\right).
\end{equation}
Moreover, equation \eqref{eq:conv_identity} shows that each B-spline is an even function. Therefore, used as a regularized delta function to interpolate a $C^2$ function $g(x)$ at a point $y$, it provides a second-order accurate interpolant:
\begin{equation}
\left|g(0)-\int_{-\infty}^{\infty}\frac{1}{h}\varphi^{\text{BS}}_n\left(\frac{x-y}{h}\right)g(y)\,\text{d}y\right| = \mathcal{O}(h^2).
\end{equation}
This second-order accuracy is sufficient for the IB method, as higher-order interpolation schemes provide little practical benefit because the velocity field in the IB method is not generally smooth, containing jumps in its normal derivatives across boundaries and generally limits the interpolation accuracy of regularized delta functions to first order.

Lee and Griffith were the first to extensively test isotropic regularized delta functions based on B-spline kernel functions\cite{lee2022lagrangian}. For cases involving stationary structures in shear-dominated flows, they found that the piecewise linear kernel yielded the most accurate results, particularly when using a Lagrangian mesh that was coarse relative to the Eulerian grid. However, for more general cases involving pressure-loaded Lagrangian structures, the three-point B-spline kernel yielded the most accurate and consistent results across various Eulerian and Lagrangian grid ratios --- typically outperforming even the IB family of kernel functions typically employed with the IB method. Gruninger et al.\cite{gruninger2024benchmarking} reached similar conclusions when applying the immersed boundary method to viscoelastic flow problems.
Throughout this study, we denote B-spline kernels as ${BS}_{\#}$, with $\#$ indicating the number of grid points in each dimension in the tensor product construction of the isotropic regularized delta function.

\subsubsection{Composite B-spline Regularized Delta Functions}
As observed in prior studies, isotropic regularized delta functions generally do not provide continuously divergence-free interpolants, even when interpolating discretely divergence-free velocity fields\cite{bao2017,gruninger2024benchmarking,peskin1993,cortez2000,griffith2012}. When combined with time-stepping and discretization errors, these interpolation errors can lead to significant volume conservation errors in IB method simulations, particularly for closed, pressurized membranes. Recently, Gruninger and Griffith constructed regularized delta functions for the IB method using composite B-spline kernels that provide continuously divergence-free interpolants of discretely divergence-free velocity fields on the MAC grid\cite{gruninger2024benchmarking,schroeder2022}.When used in conjunction with force-spreading operators, these delta functions transform continuous gradients to discrete ones with accuracy determined by the quadrature rule used in the force-spreading operation. This approach prevents the spurious, jet-like flows that typically arise when using isotropic regularized delta functions in the IB method\cite{kaiser2019,gruninger2024benchmarking}. The improved volume conservation derives from a fundamental property of B-splines: the central difference of an order $n$ B-spline equals the exact derivative of an order $n+1$ B-spline, as shown in equation \eqref{eq:conv_identity}.

The use of composite B-splines for continuously divergence-free interpolation of discretely divergence-free data appears to first date back to Hanscomb\cite{handscomb1984}, who originally applied them for this purpose. More recently, Schroeder et al.\cite{schroeder2022} later expanded this work by developing a family of 'spline chain' interpolation schemes. These schemes produce continuously divergence-free interpolants that not only achieve higher-order accuracy but also exactly reproduce the discrete vector fields at their interpolation points. Schroeder and coworkers extended these techniques to create continuously curl-free interpolants and divergence-free interpolants that are compatible with various finite difference stencils of the divergence and curl operators on the MAC grid\cite{chowdhury2024}.

Composite B-spline regularized delta functions achieve continuously divergence-free interpolation by using different one-dimensional B-spline kernels for different velocity components during interpolation and force-spreading operations\cite{gruninger2024local}. This approach uses a one-order higher B-spline kernel for the primary component direction, while applying normal-order kernels to other components. In the context of the IFED method, this yields the velocity interpolation operation
\begin{align}
    U_l&=\sum_{i,j}u_{i-1/2,j}\varphi^{\text{BS}}_{n+1}\left(\frac{x_{i-\frac{1}{2}}-\chi_l^1(t)}{h}\right)\varphi^{\text{BS}}_{n}\left(\frac{y_{j}-\chi_l^2(t)}{h}\right),\\    
    V_l&=\sum_{i,j}v_{i,j-1/2}\varphi^{\text{BS}}_{n}\left(\frac{x_i-\chi_l^1(t)}h\right)\varphi^{\text{BS}}_{n+1}\left(\frac{y_{j-\frac{1}{2}}-\chi_l^2(t)}h\right),   
\end{align}
in which $\vchi_l = \left(\chi_l^1,\chi_l^2\right)$. To ensure energy conservation\cite{peskin2002}, the force spreading operation is taken to be the discrete adjoint of the velocity interpolation operation. Because we use a consistent nodal quadrature rule with the IFED method, we avoid the need to solve linear systems of equations to project the interpolated velocity onto the finite element basis functions\cite{wells2023nodal}, which could compromise the divergence-free interpolation. We denote these composite B-spline kernels as CBS kernels, using the notation $\text{CBS}_{\#\#}$ to indicate the orders $\#$ of the B-spline kernels used in the composite construction of the regularized delta function.

\subsection{Volumetric Stabilization Approaches}
To develop a robust immersed method for incompressible hyperelastic structures undergoing large deformations, Vadala-Roth et al. \cite{vadalaroth2020} introduced a volumetric energy term and adopted the modified invariants when modeling the solid material. This stabilization approach, which adds consistency terms that vanish under grid refinement, reduces spurious volume changes in the numerical solution while maintaining the convergence properties of the underlying formulation.

\subsubsection{Volumetric Penalization}
\label{sec:volumetric_stabilization}

For flexible structures, traditional IB formulations  decompose the Cauchy stress as \cite{Griffith2017, Boffi2008}
\begin{equation}
    \cauchy = \cauchyv -p\mathbb{I} + \begin{cases} 
        \ztensor & \xb \in \fluiddom, \\ 
        \cauchys & \xb \in \soliddom, 
    \end{cases} \label{cauchy-nodev}
\end{equation}
where $\cauchyv = \mu\left(\nabla \ub + \nabla \ub \trans \right)$ is the deviatoric viscous stress and $\cauchys$ is the elastic stress. 

This formulation can lead to poor numerical results in the discretized equations, including unphysical and sometimes extreme contractions of the immersed structure \cite{vadalaroth2020}. To address this issue, the elastic stress tensor is decomposed into deviatoric and volumetric components:
\begin{equation}
    \cauchys = \dev[\cauchys] - \pstab\mathbb{I}  \label{dev-stab}
\end{equation}
where $\pstab\mathbb{I}$ acts as a stabilization term, providing an additional pressure in the solid domain that counteracts spurious compressible motions.

Following approaches used in nearly incompressible elasticity, the volumetric penalization term  $\pstab$, can be derived from volumetric energy $U(J)$ that depends only on volumetric changes in the structure:
\begin{equation}
    \pstab = -\frac{\partial U(J)}{\partial J}
\end{equation}
Note that the negative sign is due to the pressure convention in the fluid mechanics.

To control the stabilization strength, a numerical Poisson ratio $\nus$ is introduced to modulate the numerical bulk modulus $\kappas$ through the relationship:
\begin{equation}
    \kappas = \frac{2G(1 + \nus)}{3(1 - 2\nus)} \label{kappa-nu}
\end{equation}
in which $G$ is the shear modulus. This relationship mirrors the connection between the physical Poisson ratio $\nu$ and bulk modulus $\kappa$ in compressible materials. Setting $\nus = -1$ yields $\pstab = 0$, which recovers the unstabilized formulation. It is important to note that both $\kappas$ and $\nus$ are numerical parameters rather than physical ones, as the immersed structure remains incompressible in all cases.

\subsubsection{Modified Invariants}

The elastic Cauchy stress is related to the first Piola-Kirchhoff stress through
\begin{equation}
    \cauchys = \frac{1}{J}\PPs \FF\trans.
\end{equation}
For hyperelastic materials, the first Piola-Kirchhoff stress $\PPs$ can be derived from a strain energy functional $\Psi(\FF)$:
\begin{equation}
    \PPs = \frac{\partial \Psi(\FF)}{\partial \FF}.
\end{equation}

To achieve the desired decomposition of the Cauchy stress, the strain energy functional into volume-preserving and volume-changing components:
\begin{equation}
    \Psi(\FF) = W(\FF) + U(J). \label{decoupled_energy}
\end{equation}

To ensure material frame invariance \cite{bonet2008}, $\Psi$ for isotropic materials is typically expressed in terms of the first two invariants of the right Cauchy-Green tensor $\CC = \FF\trans \FF$ :
\begin{equation}
    \Psi(I_1,I_2) = W(I_1,I_2) + U(J)
\end{equation}
with $I_1 = \tr(\CC)$ and $I_2 = \frac{1}{2}(I_1^2 - \tr(\CC^2))$. This is referred to as the unmodified invariants-based model.

For nearly incompressible elasticity, it is common to use invariants that only encode shearing deformations but not volume change. To eliminate volume change information, the modified invariants are introduced\cite{flory1961thermodynamic}:
\begin{align}
    \bar{I}_1 &= J^{-2/3}I_1 \\
    \bar{I}_2 &= J^{-4/3}I_2
\end{align}
These are the invariants of the modified tensor $\bar{\CC} = \bar{\FF}\trans\bar{\FF}$. The resulting model takes the form:
\begin{equation}
    \Psi(\bar{I}_1,\bar{I}_2) = W(\bar{I}_1,\bar{I}_2) + U(J)
\end{equation}
This is referred to as the modified invariants-based model.

\subsection{IBAMR}
Fluid-structure interaction simulations described in this work use the IBAMR software infrastructure, which is a distributed-memory parallel implementation of the immersed boundary (IB) method with adaptive mesh refinement (AMR) \cite{griffith2007, griffith2009, ibamr-web-page}. IBAMR uses SAMRAI \cite{samrai-paper} for Cartesian grid discretization management and PETSc \cite{petsc-web-page, petsc-user-ref, petsc-efficient} for linear solver infrastructure.

\section{Benchmarks}
\label{sec:benchamarks}

We evaluate the performance of CBS kernels by examining their volume conservation properties and accuracy compared to those of isotropic kernels across various benchmark tests. These tests include both pressure-loaded and shear-dominated flow cases, encompassing rigid and elastic bodies. The benchmarks are designed to assess CBS kernel behavior within traditional IB and IBFE frameworks.

Several of our test cases—specifically the compression and Cook's membrane problems—are adapted from standard solid mechanics benchmarks for incompressible elasticity. While these problems traditionally examine pure solid mechanics, we extend them to fluid-structure interaction by embedding the solid bodies in an incompressible Newtonian fluid. Notably, the steady-state solutions for these FSI problems match their pure solid mechanics counterparts. The volume conservation error is quantified as the deviation of the Jacobian from unity.

We impose zero velocity boundary conditions on the computational domain to facilitate convergence to a steady state. Simulations run until a final time $T_f$, chosen such that fluid velocities effectively decay to zero. For computational efficiency, we set equal densities for both fluid and structure in all tests. While this choice does not affect steady-state solutions, it enables the use of efficient constant-coefficient linear solvers.

Unless specified, all simulations use CGS units and first-order quadrilateral ($\Qone$) elements. For grid convergence analysis, following Lee et al.\cite{lee2022lagrangian} and Wells et al.\cite{wells2023nodal}, we define the relative solid-fluid mesh factor (MFAC) as the ratio between the Lagrangian marker spacing and the Eulerian grid spacing. When evaluating volume conservation in IFED, we use the element Jacobian ($J$) as a measure of local volume change, with the error quantified by $|J-1|$.

\subsection{Pressurized Elastic Structures}

\subsubsection{A Pressurized Circular Membrane at Equilibrium}

\begin{figure}[H]
    \centering
    \includegraphics[width=0.4\linewidth]{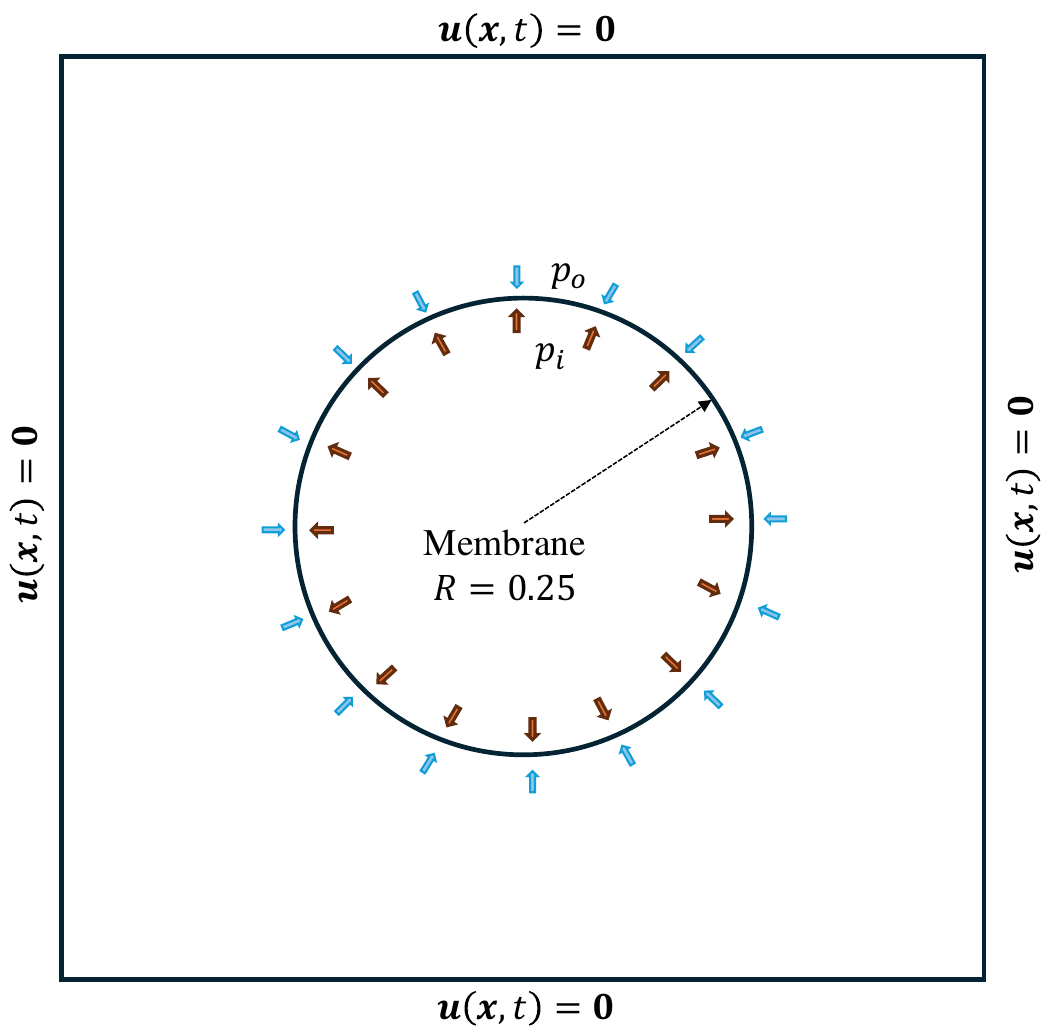}
    \caption{Pressurized membrane setup. The circular membrane (radius 0.25) at equilibrium experiences a pressure difference across its interface, with interior pressure $p_i$ and exterior pressure $p_o$.}
 
    \label{fig:pressurized_membrance_pressure}
\end{figure}

A quasi-static pressurized thin membrane was employed to assess the volume conservation properties of different computational kernels. This membrane is initialized in a circular equilibrium configuration with its center at $(1/2, 1/2)$ and radius $R = 1/4$, situated in a periodic unit square domain $[0, 1]^2$ with still background flow. The pressure difference across the boundary (Fig. \ref{fig:pressurized_membrance_pressure}) is balanced by the membrane force, with force density specified by $F(s,t) = \kappa \frac{\partial^2 X(s,t)}{\partial s^2}$, in which $X$ represents the Lagrangian coordinates and $s$ is the arclength parameter along the membrane in its reference configuration. 

In the simulation, the spring constant is set to $\kappa = 1$, and the time step $\Delta t = h / 8$, where $h = 1/128$ is the grid spacing of the computational domain. The Lagrangian mesh factor $\MFAC = 0.5$. 

Because the immersed membrane is initialized in its equilibrium configuration, any deviation from its initial area results from discretization errors in the IB method, as the fluid's incompressibility should maintain a constant area when advecting the closed membrane. We track the enclosed area using passive tracer points distributed along the immersed boundary to evaluate volume conservation accuracy. We use 10,000 Lagrangian markers as tracer points to ensure accurate area calculations. The volume conservation error is quantified through the relative area change at each timestep, defined as $\Delta A(t, X) = |A-A_0|/|A_0|$, where $A_0$ is the initial area.

\begin{figure}[H]
    \centering

    \includegraphics[width=0.4\linewidth]{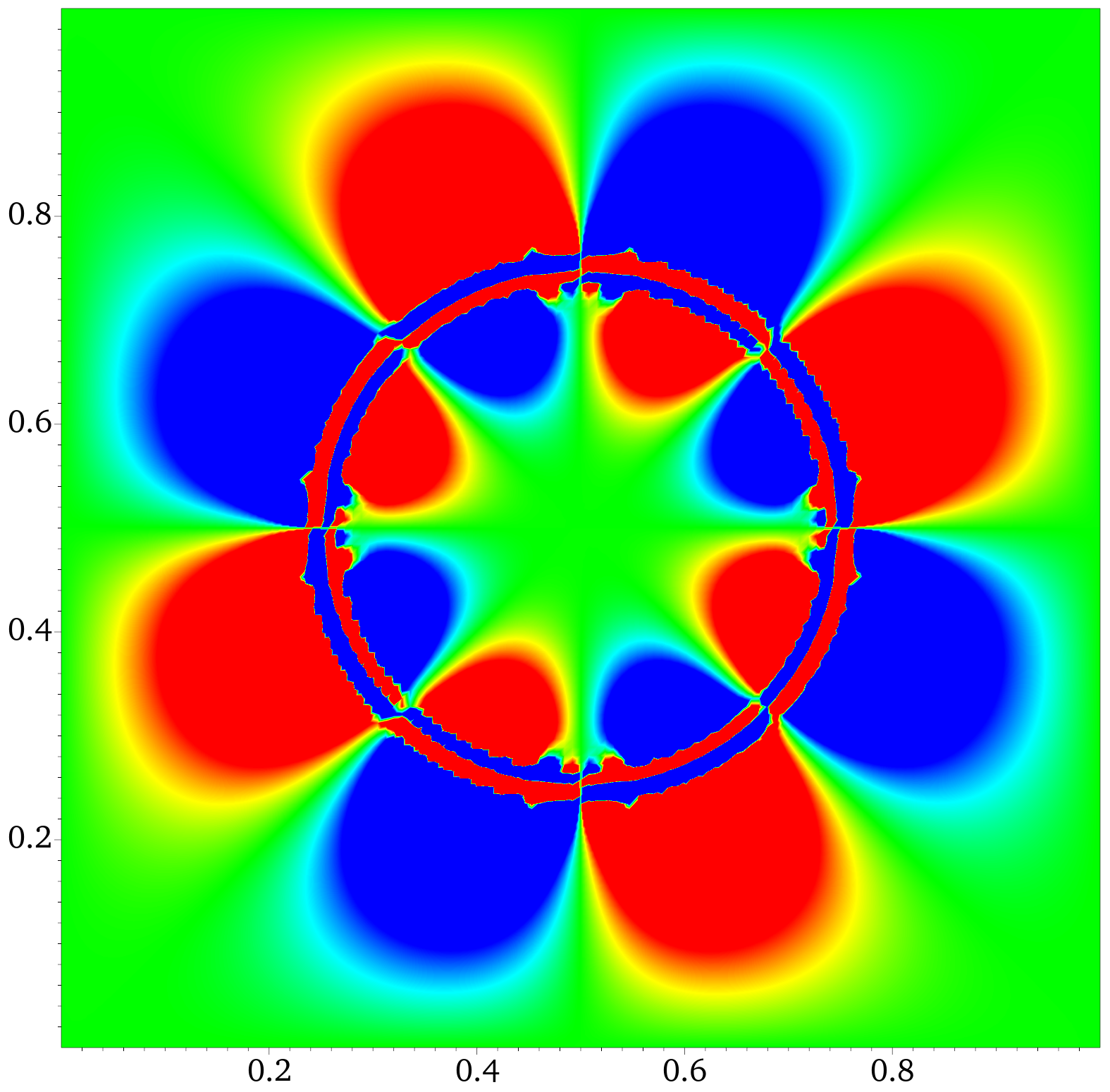}
    \includegraphics[width=0.4\linewidth]{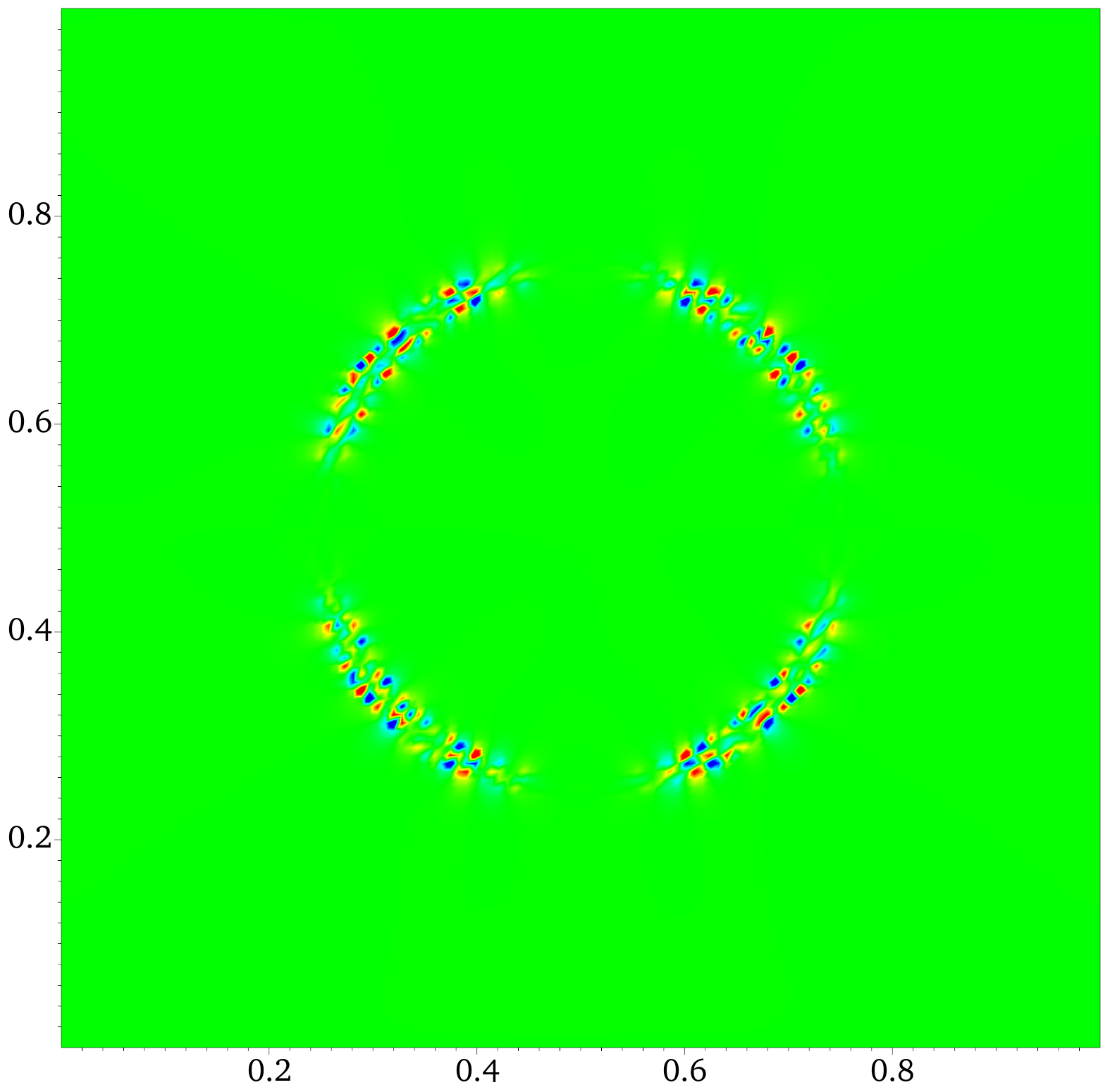}
    \includegraphics[width=0.08\linewidth]{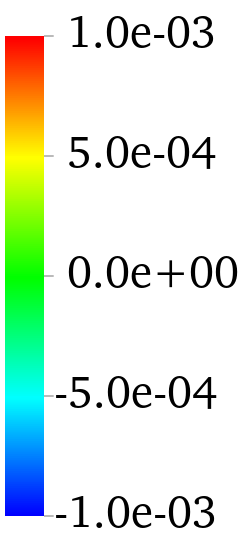}
    \caption{Vorticity comparison between $\text{IB}_{4}$ (left) and $\text{CBS}_{43}$ (right) of for the pressurized thin elastic membrane. The vorticity should ideally be zero everywhere. The CBS kernel exhibits a noticeably smaller error induced by non-zero velocity due to errors in the discrete force spreading operator.}
    \label{fig:pressurized_membrance_vorticity}
\end{figure}

Theoretically, the vorticity should be zero throughout the domain\cite{gruninger2024local}. The CBS kernels demonstrate minimal vorticity error arising from spurious non-zero velocities, which occur due to errors associated with the approximation of the force spreading operator (Fig. \ref{fig:pressurized_membrance_vorticity}). 

\begin{figure}[H]
    \centering
    \includegraphics[width=0.5\linewidth]{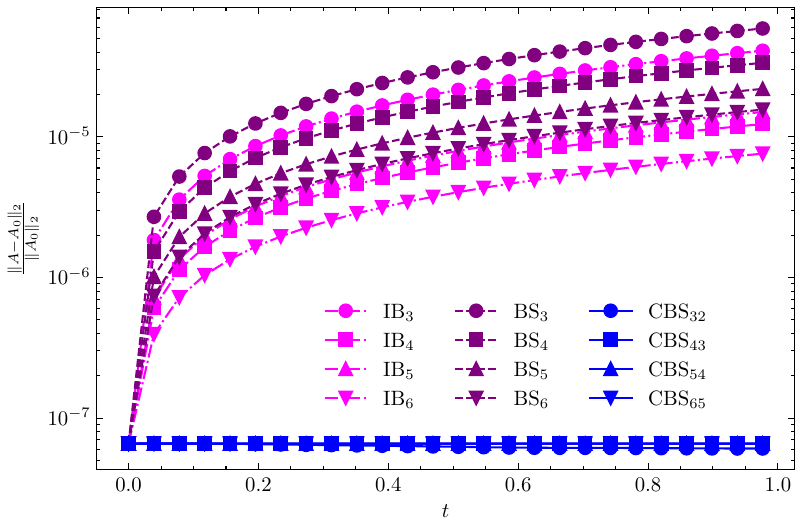}
    \caption{Volume conservation comparison of pressurized thin elastic membrane for different kernels. }
    \label{fig:pressurized_membrance_area_conservation}
\end{figure}

As illustrated in Figure \ref{fig:pressurized_membrance_area_conservation}, CBS kernels achieve superior volume conservation compared to conventional IB and BS kernels, reducing volume errors by approximately two orders of magnitude. Although wider kernels improve volume conservation for both IB and BS methods, the CBS kernels maintain consistently excellent performance regardless of kernel width, with negligible differences between variants.

\begin{figure}[H]
    \centering
    \includegraphics[width=0.5\linewidth]{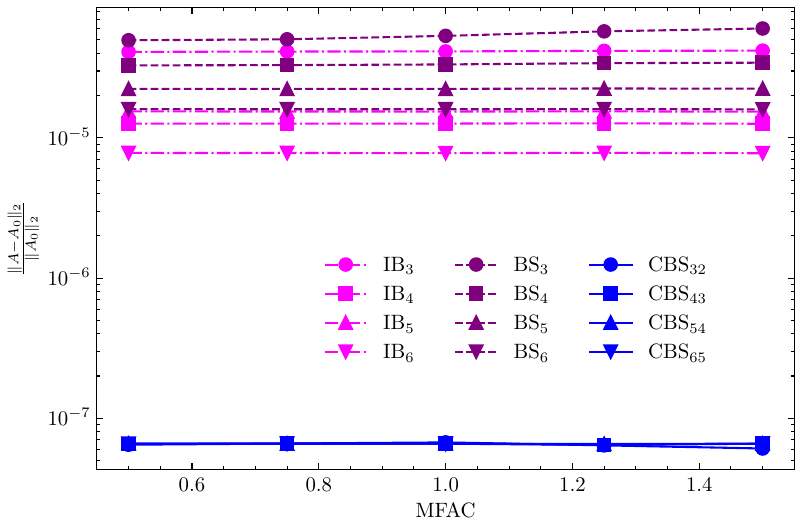}
    \caption{The effects of MFAC on the volume conservation error of pressurized thin elastic membrane for different kernels. }
    \label{fig:pressurized_membrance_area_error_convergence_mfac}
\end{figure}

To examine the sensitivity of grid spacing ratios between the solid and fluid meshes (MFAC), we fix the fluid grid with a high resolution ($N=128$) and vary MFAC as 0.5, 0.75, 1.0, 1.25, and 1.5. Fig. \ref{fig:pressurized_membrance_area_error_convergence_mfac} presents grid convergence results for different MFAC values based on the volume error $t$ = 1 s. For this test case, all kernels are insensitive to MFAC.

\subsubsection{Two-dimensional Pressure-loaded Elastic Band}
This benchmark examines the deformation of a thin immersed elastic band under pressure loading, which presents several numerical challenges. First, maintaining accurate volume conservation is particularly challenging when the structure thickness is comparable to or smaller than the fluid grid spacing. Additionally, the added-mass effect becomes more pronounced for thin structures, especially in fluids of similar density as in this case. Further, standard constitutive models may become unstable under the large deformations experienced by thin structures, necessitating special stabilization techniques \cite{vadalaroth2020}.

Following prior studies \cite{wells2023nodal, vadalaroth2020, lee2022lagrangian}, we use a $2L \times L$ computational domain ($L = 1$ cm) discretized by a uniform Cartesian grid with $N = 128$ cells in each direction on a single level. The elastic band measures $0.8L \times h$, where $h = L/10$ is its thickness. As shown in Fig. \ref{fig:elastic-band-mesh}, the band is fixed by rigid blocks at its top and bottom ends and discretized using $\Qone$ elements.

The material properties include a fluid density of 1 g/cm$^3$ and a dynamic viscosity of 0.01 g/cm$\cdot$s. The band is modeled as an incompressible neo-Hookean material with a shear modulus of 200 dyn/cm$^2$.

The time step size is initially set to $\Delta t = 10^{-3}\euleriandx$ and adaptively reduced to maintain numerical stability throughout the simulation.
Boundary conditions are prescribed as follows: fluid traction boundary conditions $\bbsigma (\xb, t)\nb(\xb) = -\boldsymbol{\tau}$ and $\bbsigma (\xb, t)\nb(\xb) = \boldsymbol{\tau}$ are imposed on the left and right domain boundaries respectively, where $\boldsymbol{\tau}=(5,0)^T$. Zero velocity conditions are enforced along the top and bottom boundaries. Using an unmodified invariants material model without volumetric energy terms, the simulation runs until reaching a steady state at $t = 10$ s.

\begin{figure}[H]
    \centering
    \includegraphics[width=0.5\linewidth]{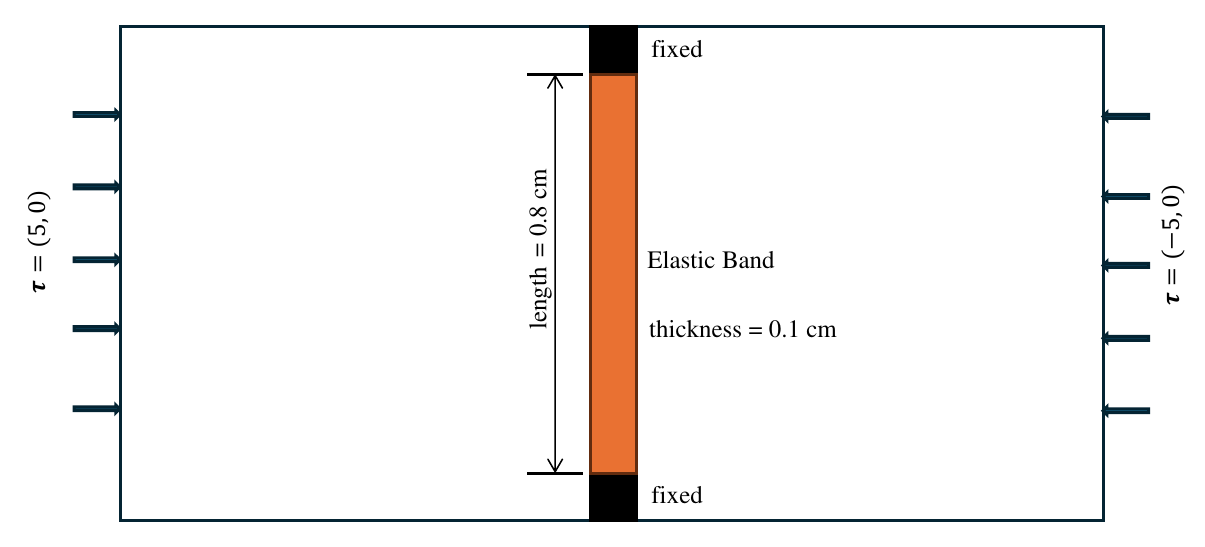}
    \caption{Schematic of a two-dimensional pressure-loaded elastic band. The elastic band is fixed at the top and bottom and experiences pressure differences across the band. }
    \label{fig:elastic-band-mesh}
\end{figure}

Table \ref{tab:thick_band_x_disp} shows the maximum $x$-displacement of the elastic band, measured at its center ($y$ = 0.5). All kernels demonstrate consistent convergence behavior as MFAC decreases. However, stability issues emerge with larger MFAC values: simulations fail to reach steady-state when MFAC exceeds 1.0. At MFAC = 1.0, wider kernels generally exhibit better stability, though this relationship is not strictly monotonic -- increased kernel width does not always guarantee improved stability.

\begin{table}[H]
  \centering
  \caption{The maximum horizontal displacement (cm) of the thick elastic band at steady-state of different kernels with different MFACs (missing value represents the simulation fails to reach the steady-state).}
    \begin{tabular}{lcccccc}
    \toprule
    \multicolumn{1}{c}{\multirow{2}[4]{*}{Kernel}} & \multicolumn{6}{c}{MFAC} \\
\cmidrule{2-7}      & \multicolumn{1}{r}{0.25} & \multicolumn{1}{r}{0.5} & \multicolumn{1}{r}{0.75} & \multicolumn{1}{r}{1.0} & \multicolumn{1}{r}{1.25} & \multicolumn{1}{r}{1.5} \\
    \midrule
    $\text{IB}_3$ & 0.11935 & 0.11787 & 0.11230 & \textbackslash{} & \textbackslash{} & \textbackslash{} \\
    $\text{IB}_4$ & 0.11943 & 0.11862 & 0.11510 & \textbackslash{} & \textbackslash{} & \textbackslash{} \\
    $\text{IB}_5$ & 0.12056 & 0.11905 & 0.11536 & 0.10880 & \textbackslash{} & \textbackslash{} \\
    $\text{IB}_6$ & 0.12073 & 0.11998 & 0.11802 & \textbackslash{} & \textbackslash{} & \textbackslash{} \\
    $\text{BS}_3$ & 0.11931 & 0.11662 & 0.11153 & \textbackslash{} & \textbackslash{} & \textbackslash{} \\
    $\text{BS}_4$ & 0.11999 & 0.11817 & 0.11339 & \textbackslash{} & \textbackslash{} & \textbackslash{} \\
    $\text{BS}_5$ & 0.12045 & 0.11883 & 0.11372 & \textbackslash{} & \textbackslash{} & \textbackslash{} \\
    $\text{BS}_6$ & 0.12060 & 0.11905 & 0.11473 & 0.10918 & \textbackslash{} & \textbackslash{} \\
    $\text{CBS}_{32}$ & 0.12145 & 0.12197 & 0.12135 & 0.11742 & \textbackslash{} & \textbackslash{} \\
    $\text{CBS}_{43}$ & 0.12231 & 0.12247 & 0.12282 & \textbackslash{} & \textbackslash{} & \textbackslash{} \\
    $\text{CBS}_{54}$ & 0.12267 & 0.12285 & 0.12319 & 0.12798 & \textbackslash{} & \textbackslash{} \\
    $\text{CBS}_{65}$ & 0.12277 & 0.12299 & 0.12338 & 0.12437 & \textbackslash{} & \textbackslash{} \\
    \bottomrule
    \end{tabular}%
  \label{tab:thick_band_x_disp}%
\end{table}%

Figs. \ref{fig:elastic-band-thick-mfac-isotropic} and \ref{fig:elastic-band-thick-mfac-cbs32} show color maps of the Jacobian values in the deformed structure at steady state, illustrating incompressibility errors for different kernels and MFACs. The most significant errors occur along the fluid-structure interface, with all kernels showing reduced volumetric errors as MFAC decreases. $\text{IB}_5$ and $\text{BS}_6$ exhibit similar performance at the same MFAC values. However, CBS kernels demonstrate significantly smaller errors compared to isotropic kernels when MFAC is less than 1.0.

\begin{figure}[H]
    \centering
    \begin{subfigure}[b]{0.9\textwidth} \ContinuedFloat
        \centering
        \includegraphics[width=0.15\textwidth]{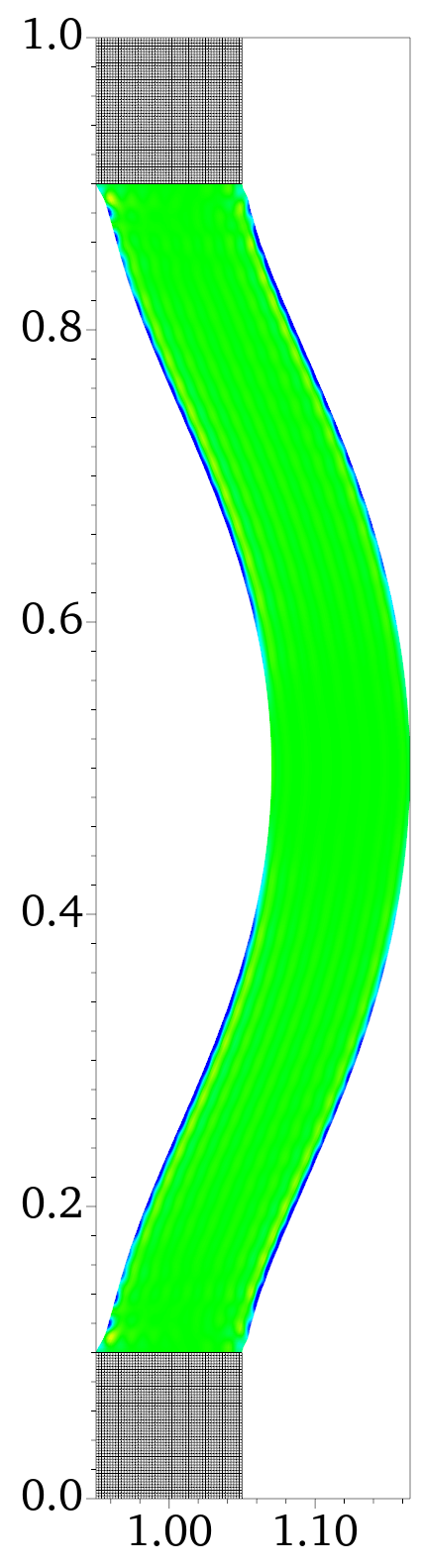}\hspace{0.04\textwidth} 
        \includegraphics[width=0.15\textwidth]{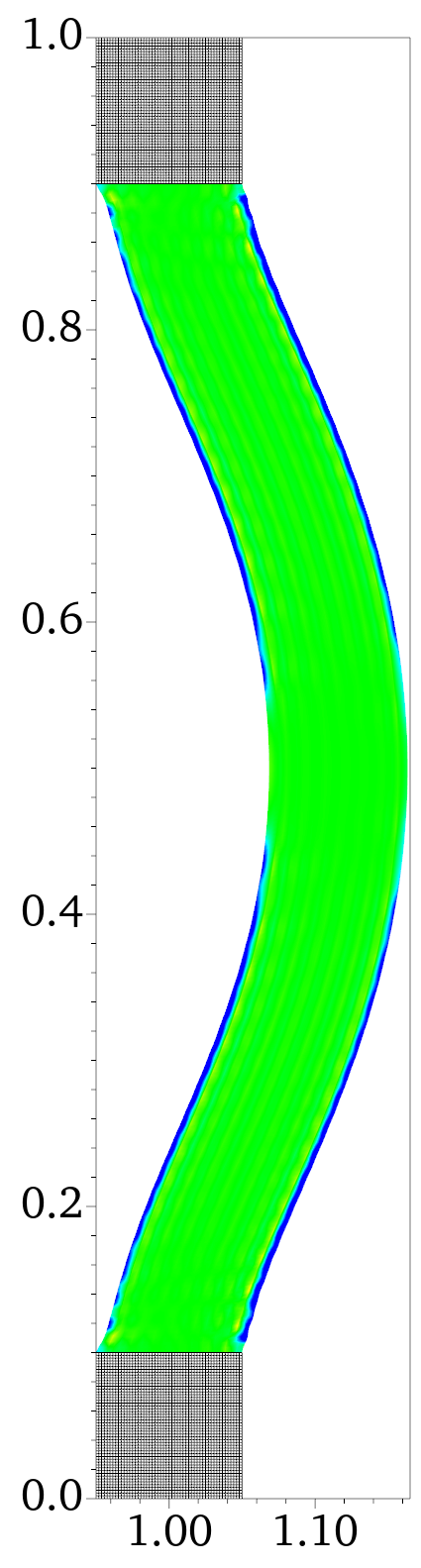}\hspace{0.04\textwidth} 
        \includegraphics[width=0.15\textwidth]{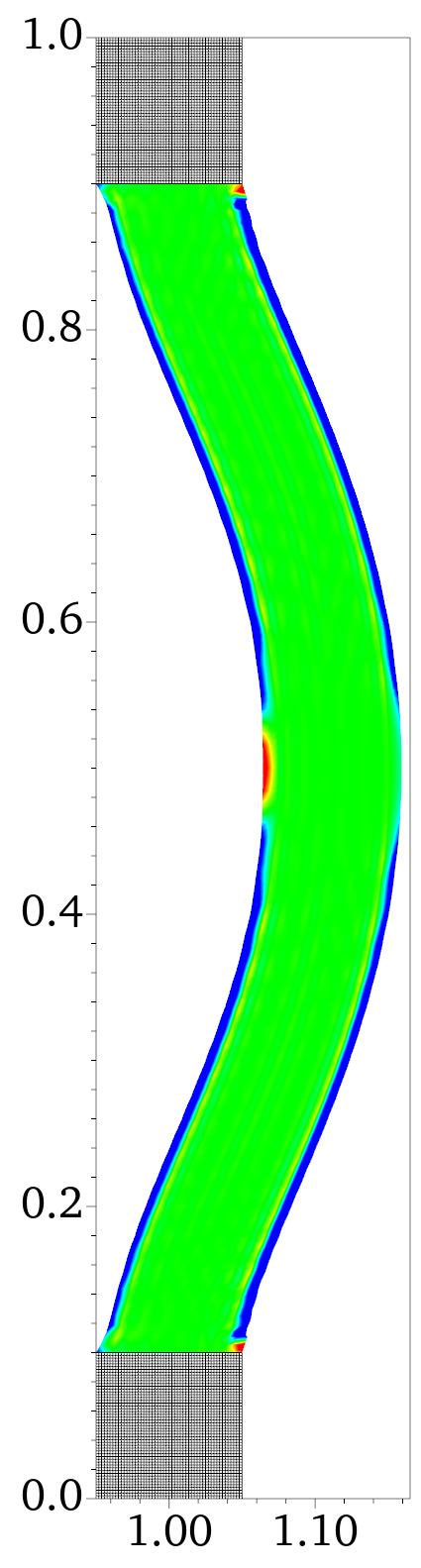} \hspace{0.04\textwidth} 
        \includegraphics[width=0.15\textwidth]{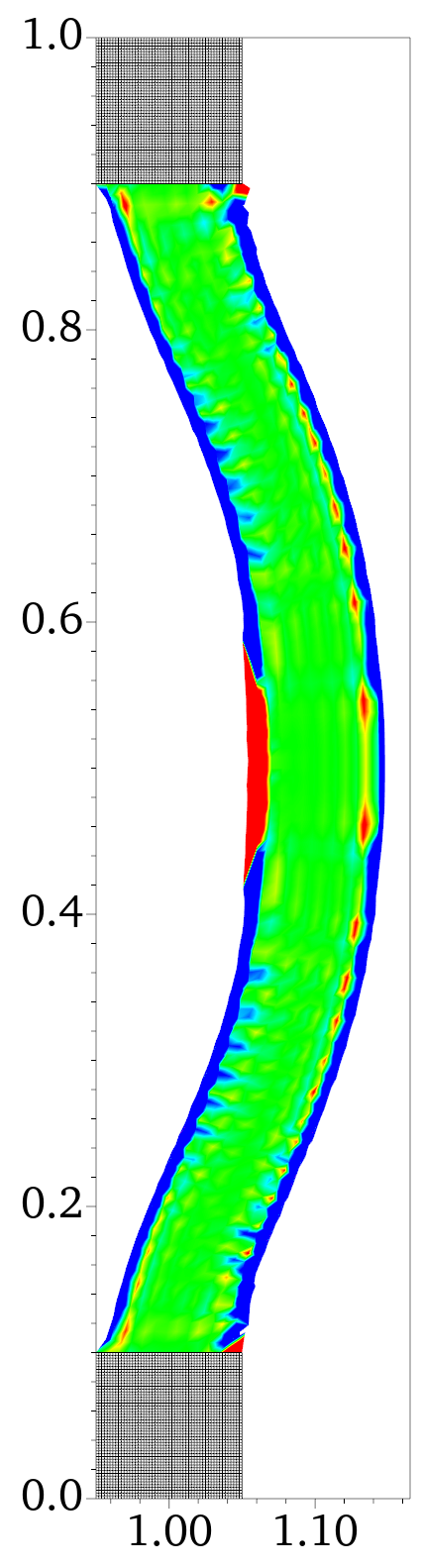}    
         \includegraphics[width=0.07\linewidth]{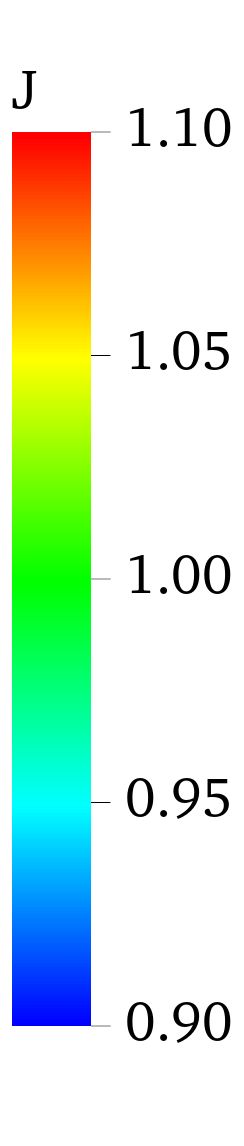}
        \caption{$\text{IB}_5$ with MFAC = 0.2, 0.5, 0.75 and 1.0 (left to right)}
    \end{subfigure}

    \begin{subfigure}[b]{0.9\textwidth} \ContinuedFloat
        \centering
        \includegraphics[width=0.14\textwidth]{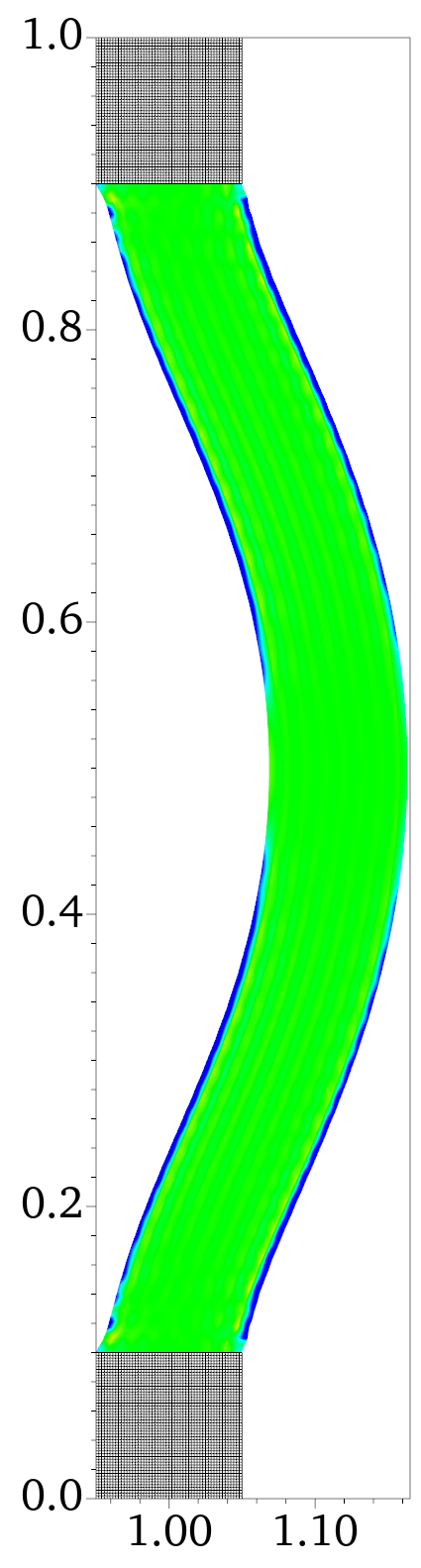}\hspace{0.04\textwidth} 
        \includegraphics[width=0.14\textwidth]{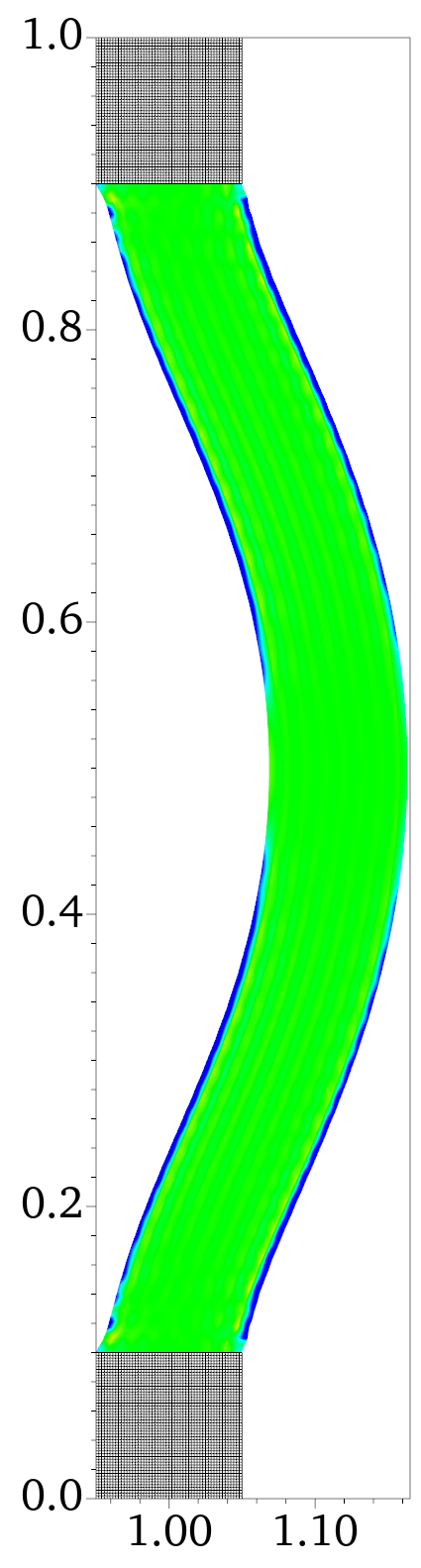}\hspace{0.04\textwidth} 
        \includegraphics[width=0.14\textwidth]{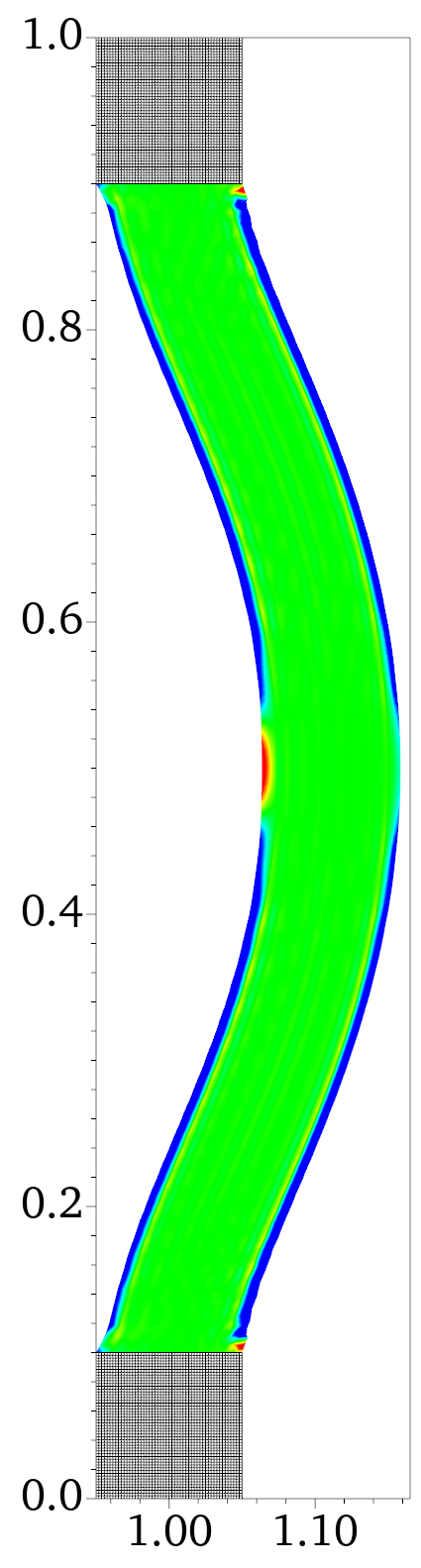}\hspace{0.04\textwidth}  
        \includegraphics[width=0.14\textwidth]{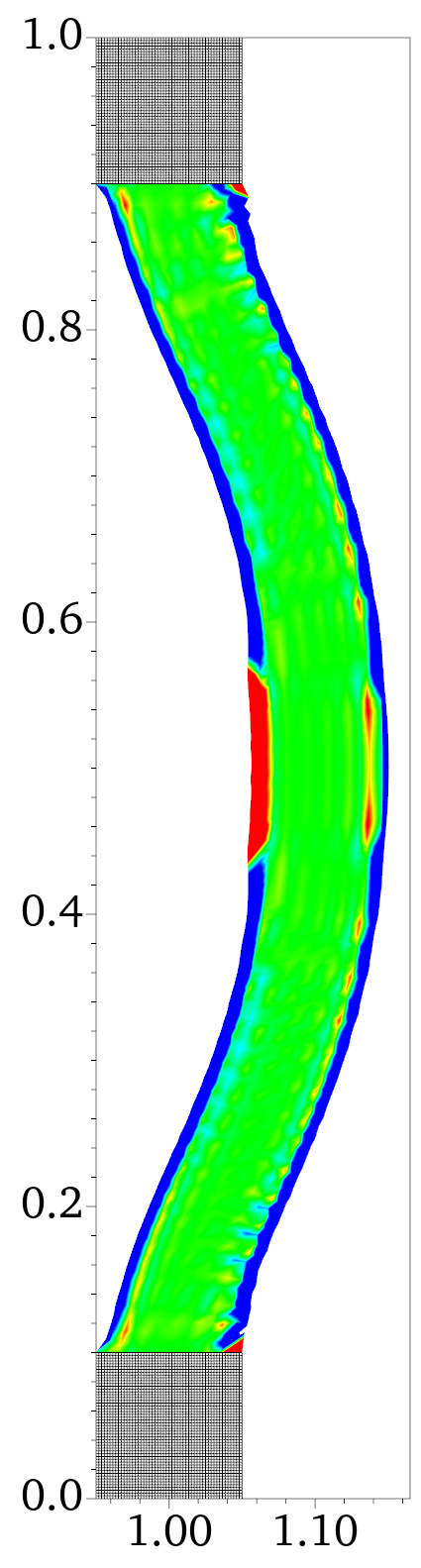}    
         \includegraphics[width=0.07\linewidth]{images/elastic_band/J_legend.png}
        \caption{$\text{BS}_6$ with MFAC = 0.2, 0.5, 0.75, and 1.0 (left to right)}
    \end{subfigure}

    \caption{Comparison of the influence of MFAC on the volume conservation of the pressure-loaded elastic band to isotropic kernels. (a) $\text{IB}_5$ and (b) $\text{BS}_6$. Generally, smaller MFAC yields better results regarding the volume conservation for all kernels. There are indistinct differences between $\text{IB}_5$ and $\text{BS}_6$ with the same MFAC.}
    \label{fig:elastic-band-thick-mfac-isotropic}
\end{figure}

\begin{figure}
\centering
    \begin{subfigure}[b]{0.9\textwidth}
        \centering
        \includegraphics[width=0.15\textwidth]{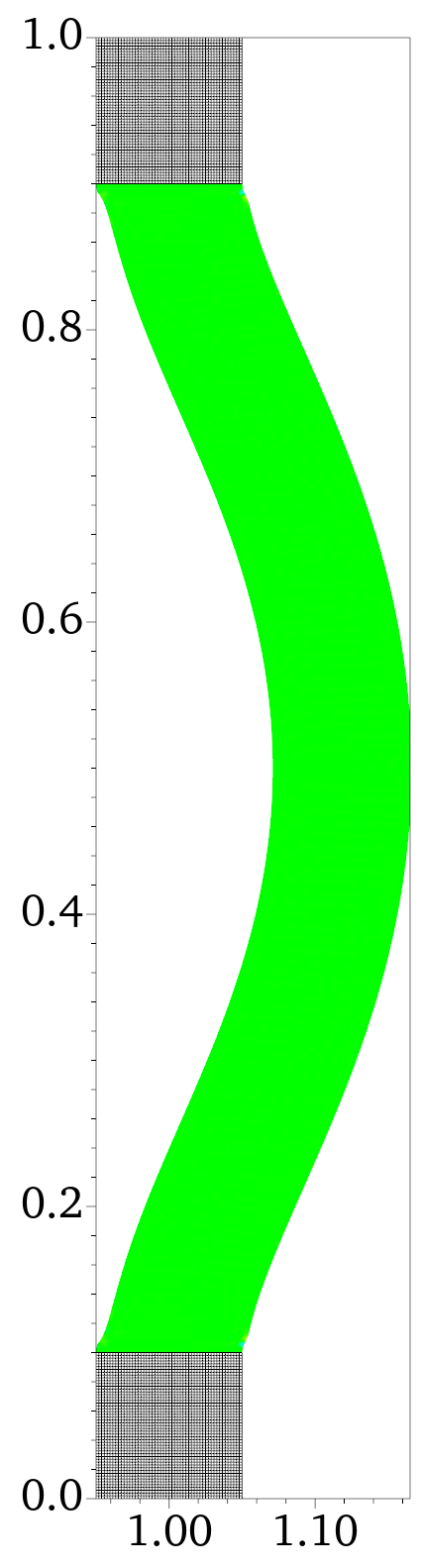}\hspace{0.04\textwidth} 
        \includegraphics[width=0.15\textwidth]{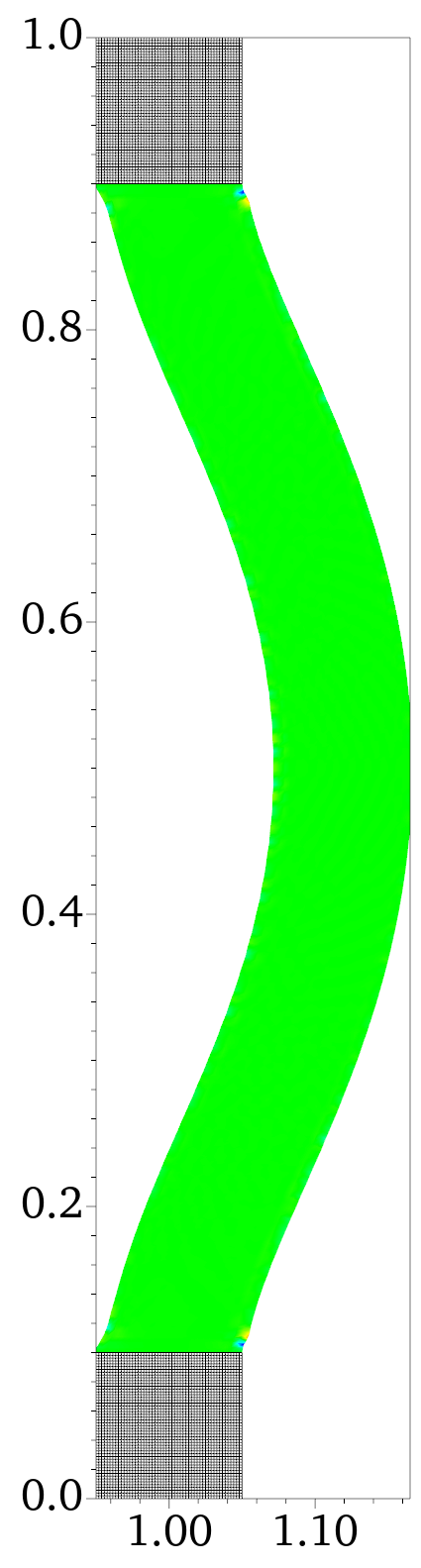}\hspace{0.04\textwidth} 
        \includegraphics[width=0.15\textwidth]{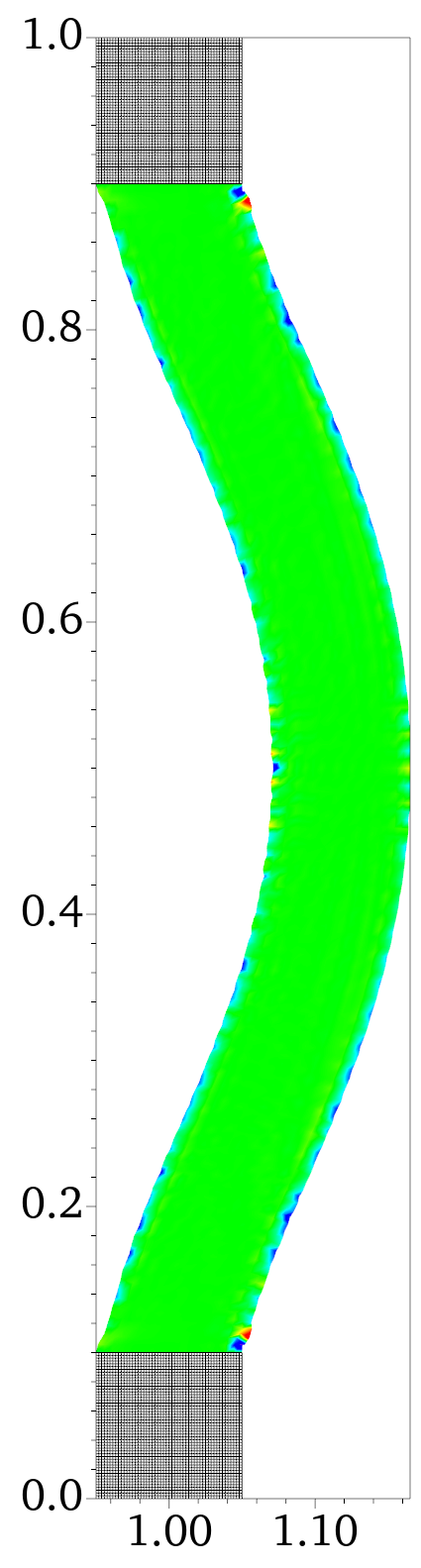}\hspace{0.04\textwidth}  
        \includegraphics[width=0.15\textwidth]{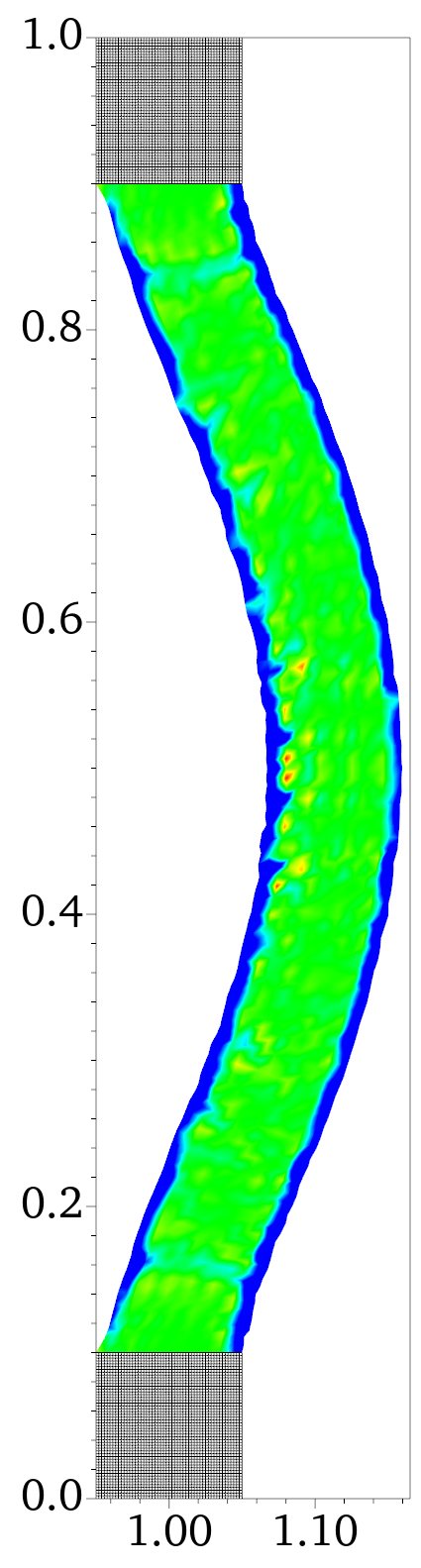}    
         \includegraphics[width=0.07\linewidth]{images/elastic_band/J_legend.png}
    \end{subfigure}
    \caption{Comparison of the influence of MFAC on the volume conservation of the pressure-loaded elastic band using  $\text{CBS}_{32}$ with MFAC = 0.2, 0.5, 0.75, and 1.0 (left to right). Compared with isotropic kernels (Fig. \ref{fig:elastic-band-thick-mfac-isotropic}), CBS kernels show a significantly smaller error when MFAC is less than 1.0.}
    \label{fig:elastic-band-thick-mfac-cbs32}   
\end{figure}

This is more prominent when the elastic band is thinner. We reduce the thickness of the elastic band to $L/32$ with only 2 elements across the thickness direction.

\begin{figure}[H]
    \centering
    \begin{subfigure}[b]{0.18\textwidth} 
    \centering
    \includegraphics[width=\textwidth]{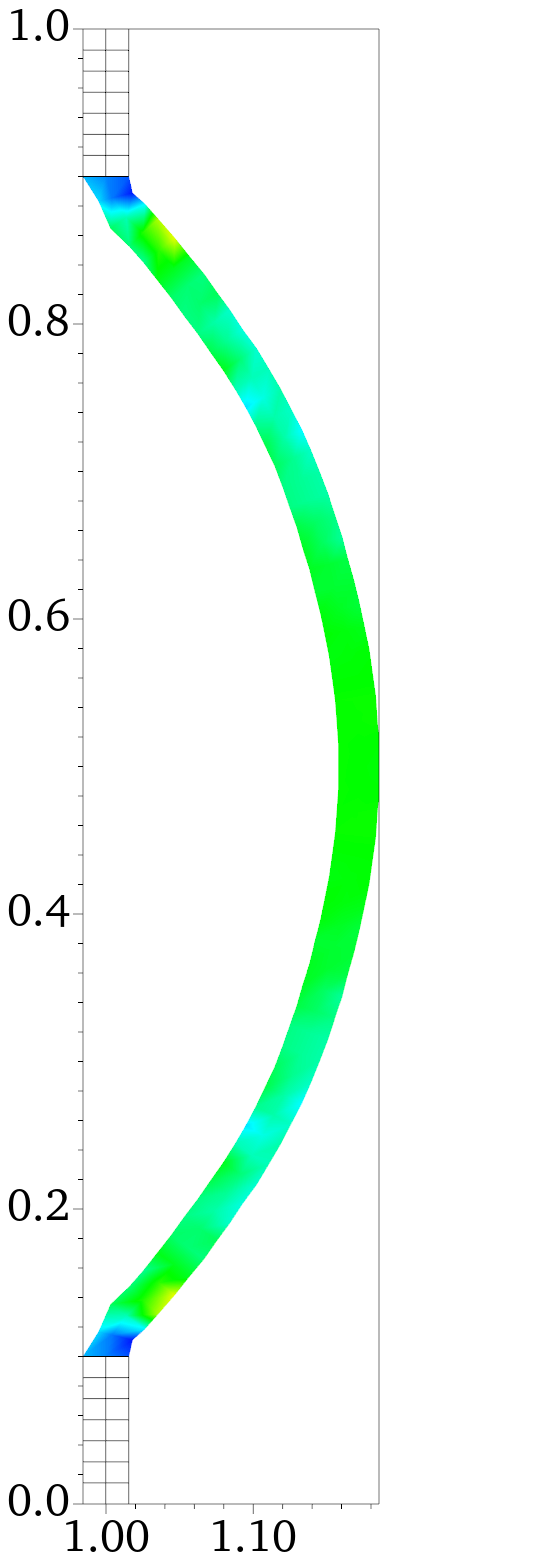}
    \caption{MFAC = 0.5}
    \end{subfigure}
    \hspace{0.02\textwidth} 
    \begin{subfigure}[b]{0.18\textwidth}   
    \centering
    \includegraphics[width=\textwidth]{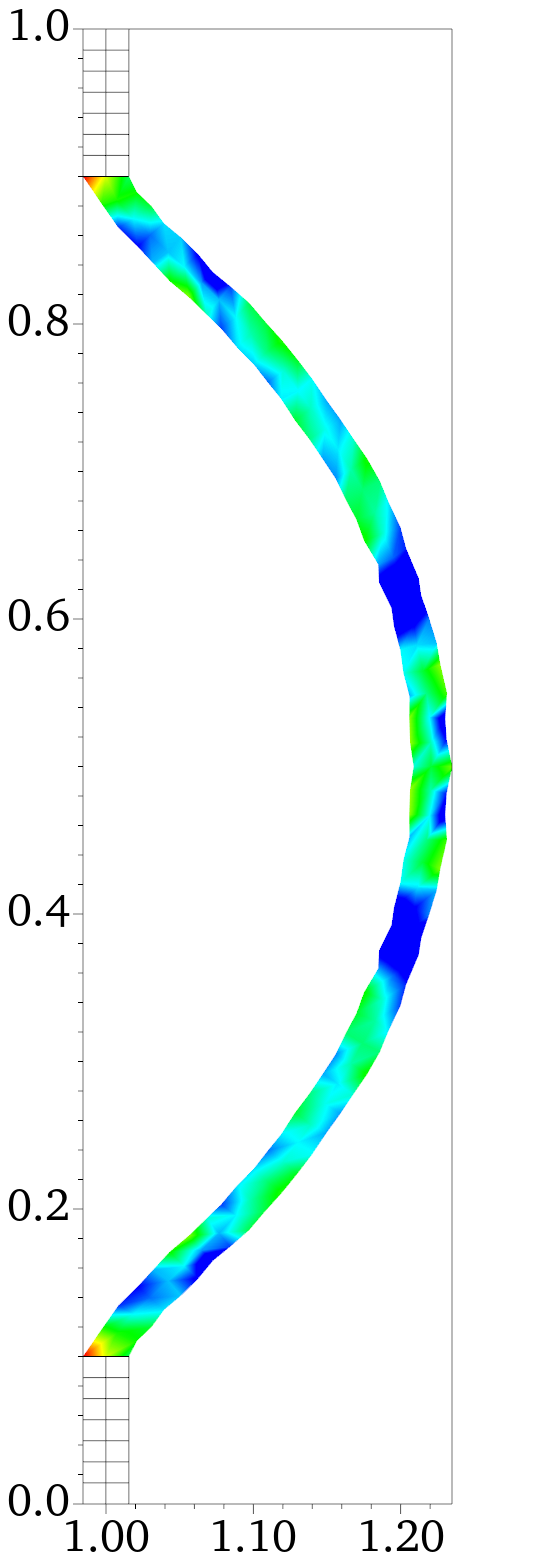}
    \caption{MFAC = 1.0}
    \end{subfigure}  
    \hspace{0.02\textwidth}
    \begin{subfigure}[b]{0.153\textwidth}
    \centering
    \includegraphics[width=\textwidth]{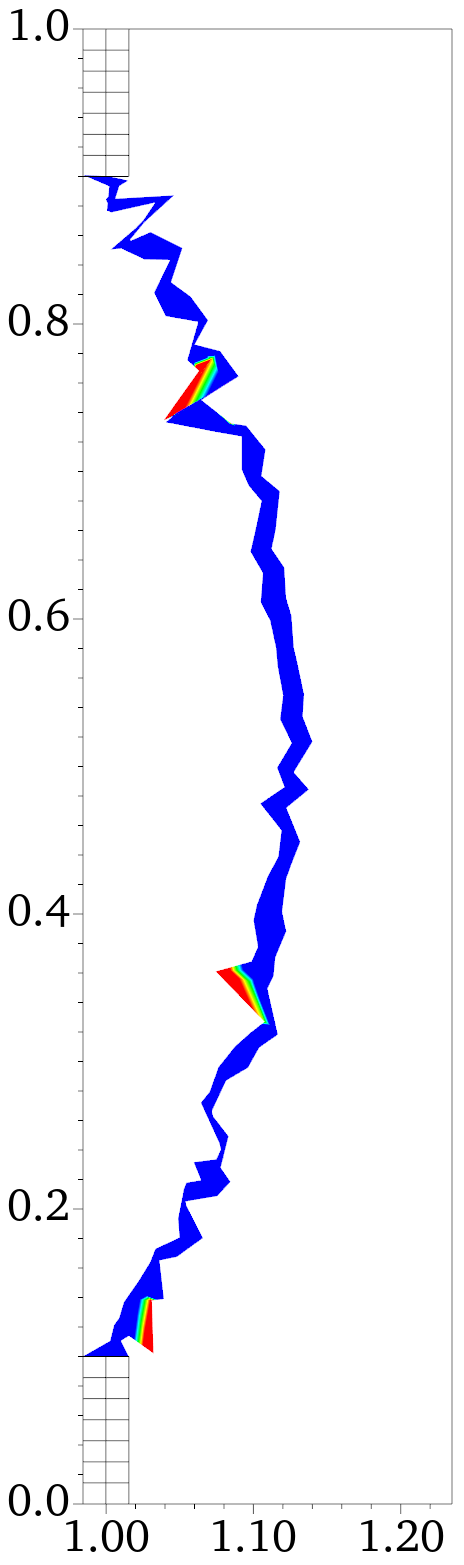}
    \caption{MFAC = 1.5}
    \end{subfigure}
    \begin{subfigure}[b]{0.06\textwidth}
    \centering
    \includegraphics[width=\linewidth]{images/elastic_band/J_legend.png}
    \end{subfigure}

    \caption{The deformation of the elastic band for $\text{CBS}_{32}$ at $t$ = 0.2 s for different MFAC values. The color map represents the elemental Jacobian values. At $\text{MFAC}=0.5$, the elastic band deforms smoothly and the elemental Jacobians remain close to unity. As MFAC increases to 1.0 and 1.5, the elemental Jacobians deviate increasingly from unity. For $\text{MFAC} \ge 1$, as the simulation proceeds, eventually the band's deformation becomes spurious and nonphysical.}
    \label{fig:elastic-band-mfac}
\end{figure}

Fig. \ref{fig:elastic-band-mfac} demonstrates the sensitivity of the $\text{CBS}_{32}$ kernel to MFAC. The results indicate that CBS kernels require a smaller MFAC for stability than traditional kernels. When using neither volume energy term nor modified invariants, three distinct behaviors are observed: At MFAC = 0.5, the simulation remains stable and completes successfully; at MFAC = 1.0, the structure develops crinkling instabilities at $t$ = 0.2 s, ultimately failing due to excessive element distortion at MFAC = 1.5, numerical instabilities emerge even earlier, leading to premature simulation failure. These findings align with Wells et al. \cite{wells2023nodal}, despite differences in element order. While they used second-order elements and we employed first-order elements, both studies demonstrate that nodal coupling experiences significant accuracy loss with relatively coarse structural meshes (MFAC $\geq$ 1) for this case.

\begin{figure}[H]
    \centering
    \begin{subfigure}[b]{0.2\textwidth}
        \centering
        \includegraphics[width=\textwidth]{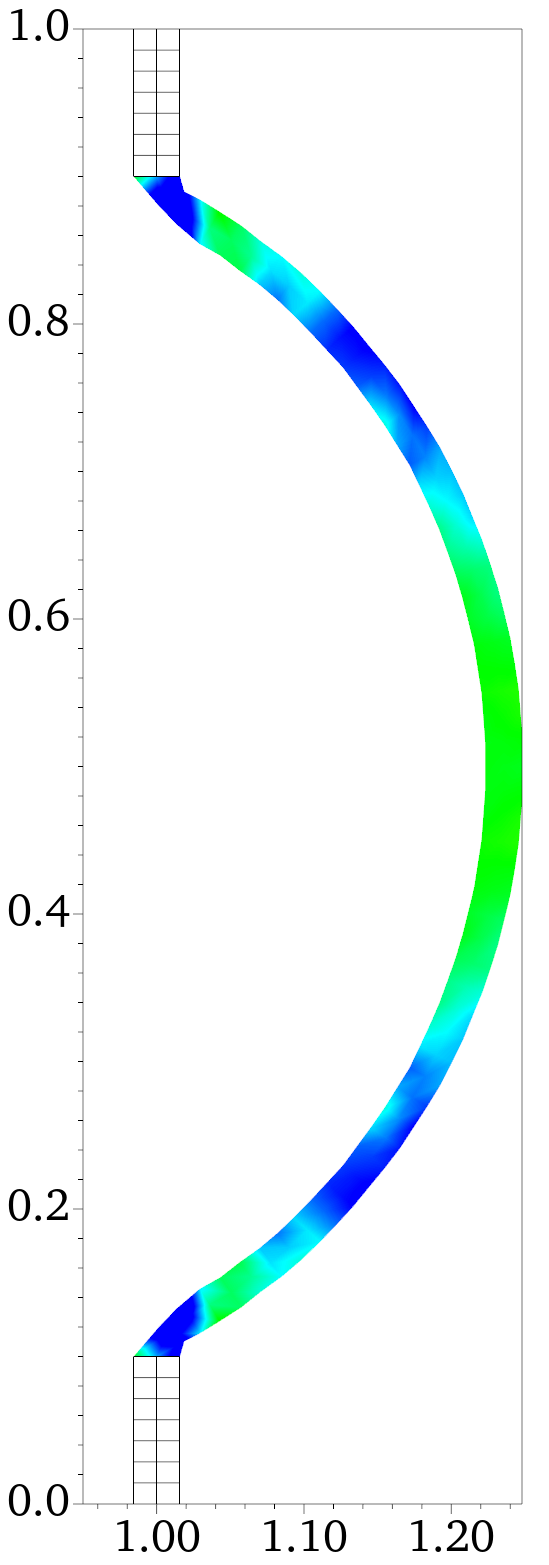}
        \caption{}
    \end{subfigure}
    \hspace{0.02\textwidth}
    \begin{subfigure}[b]{0.2\textwidth}
        \centering
        \includegraphics[width=\textwidth]{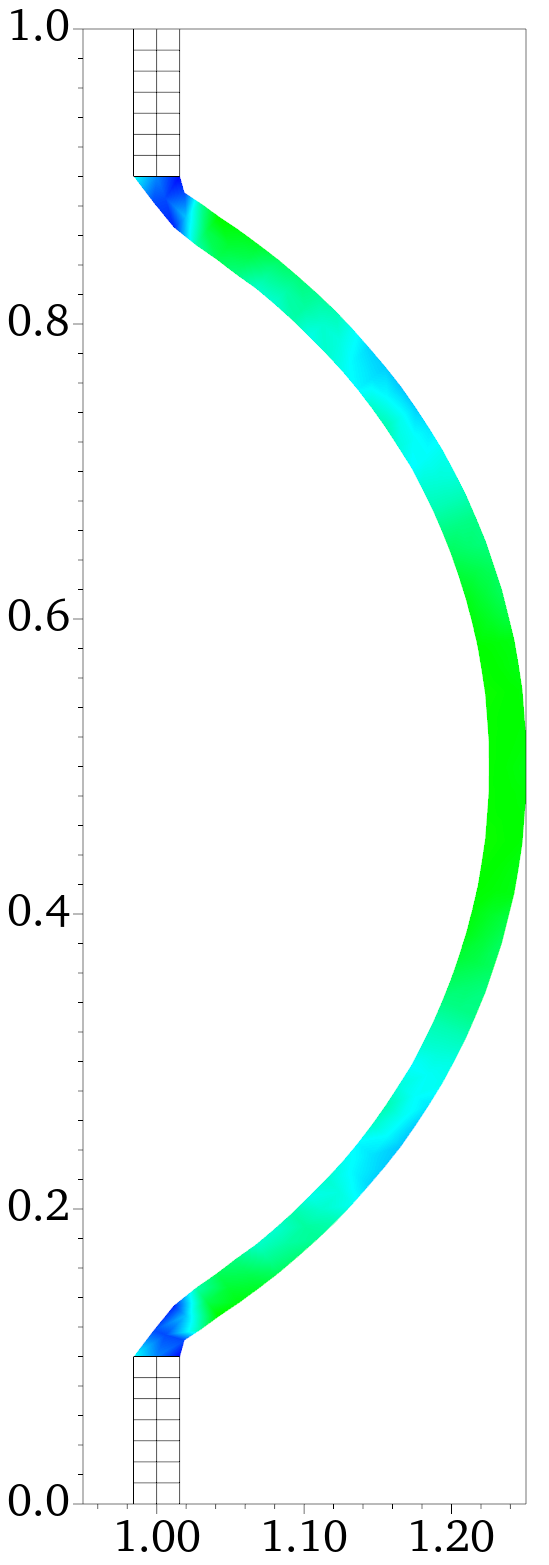}
        \caption{}
    \end{subfigure}
    \hspace{0.02\textwidth}
    \begin{subfigure}[b]{0.2\textwidth}
        \centering
        \includegraphics[width=\textwidth]{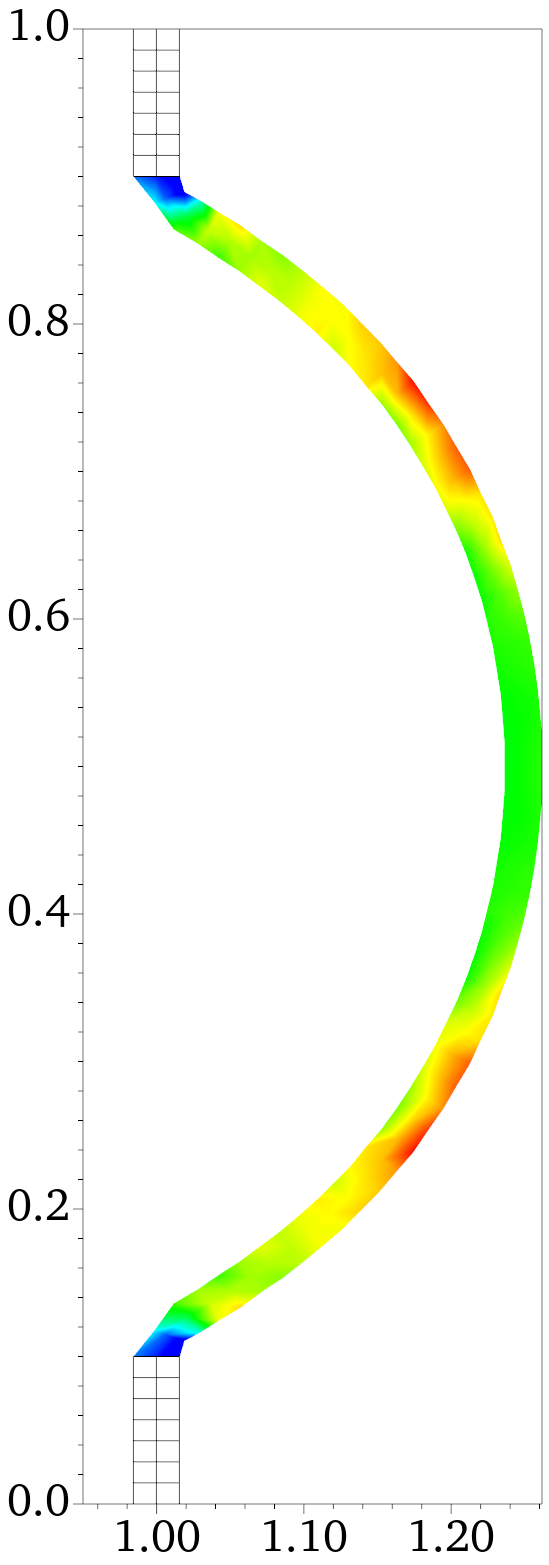}    
        \caption{}
    \end{subfigure}
    \hspace{0.02\textwidth}
    \begin{subfigure}[b]{0.2\textwidth}
        \centering
        \includegraphics[width=\textwidth]{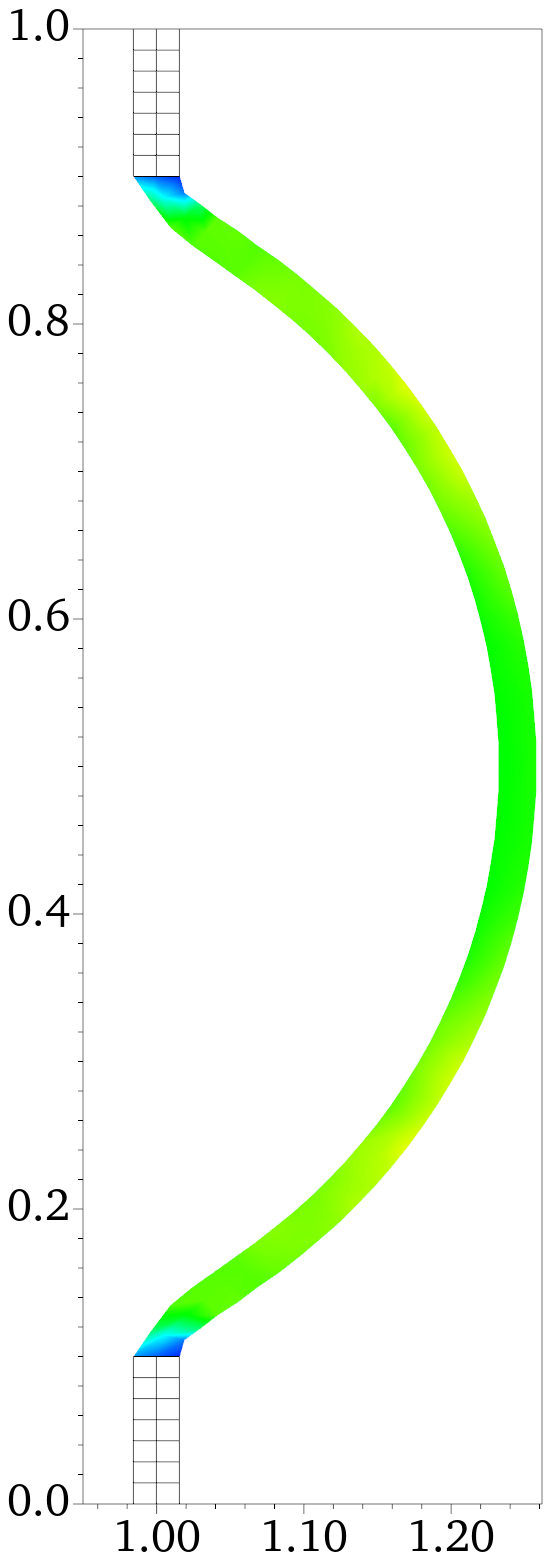}    
        \caption{}
    \end{subfigure}
    \begin{subfigure}[b]{0.07\textwidth}
    \centering
    \includegraphics[width=\linewidth]{images/elastic_band/J_legend.png}
    \end{subfigure}    
\caption{Comparison of stabilization treatments for $\text{CBS}_{32}$ with MFAC = 0.5 ($t$ = 10 s), visualized through element Jacobian distributions: (a) without a volumetric penalty term and using unmodified invariants, (b) including only a volumetric penalty term, (c) using modified invariants only, (d) using both a volumetric penalty and modified invariants.}
    \label{fig:elastic-band-stab}
\end{figure}

To evaluate the effectiveness of both volume energy terms and modified invariant treatments, we conducted simulations at MFAC = 0.5 using different combinations of these approaches. Fig. \ref{fig:elastic-band-stab} compares the Jacobian values at $t = 10$ s across these cases. The results show a hierarchy of improvements: the volume energy term provides modest enhancement, modified invariants yield better results, and combining both treatments significantly improves solution quality.

\subsection{Prescribed Deformation of Elastic Bodies}
\subsubsection{Compression Test}
\label{subsec:compression-test}

\begin{figure}[H]
    \centering
    \includegraphics[width=0.5\linewidth]{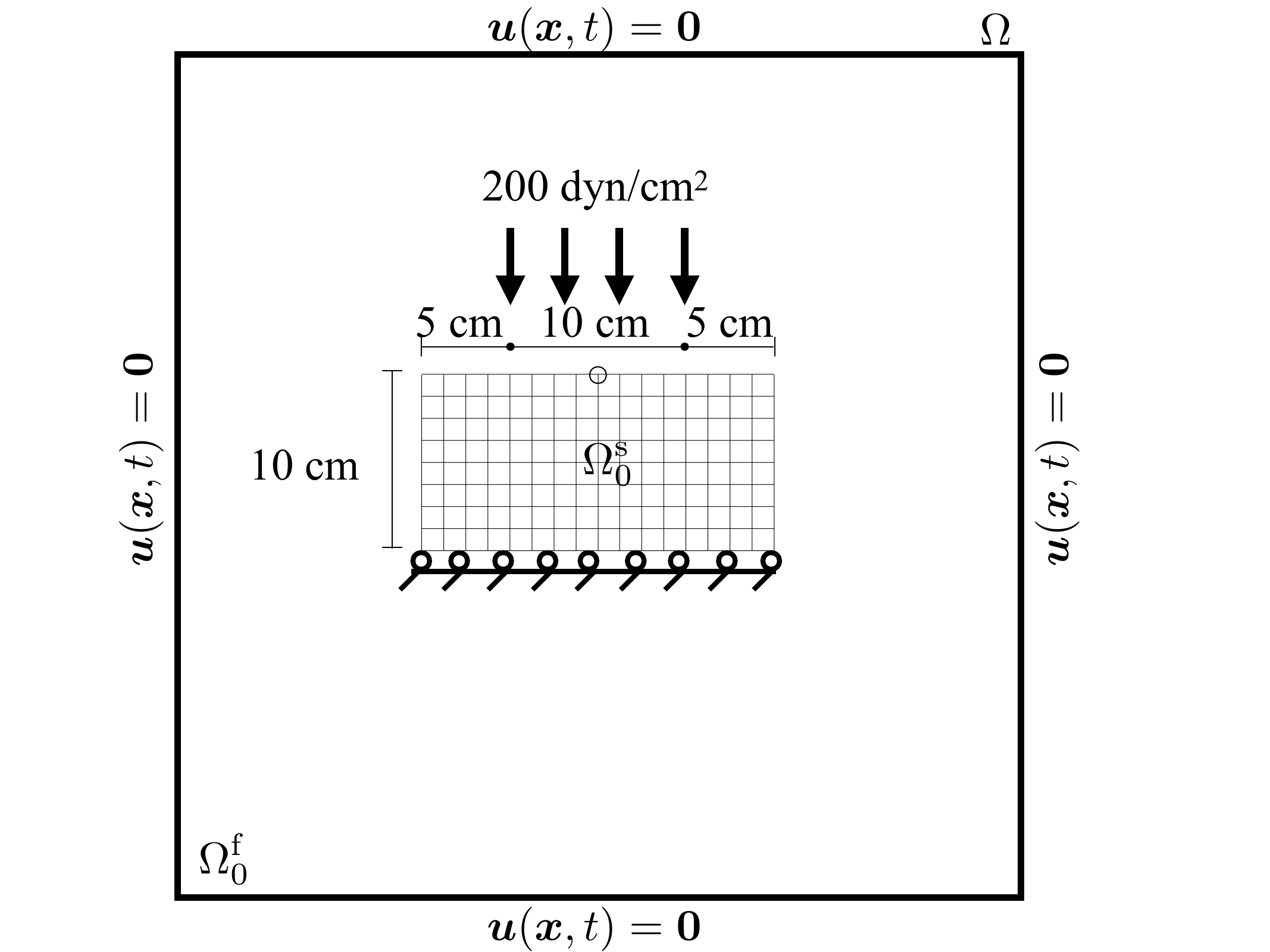}
    \caption{Schematic of the compressed block benchmark}
    \label{fig:compresed_block_setup}
\end{figure}

Following the configuration presented in the previous study \cite{wells2023nodal}, this benchmark examines a plane strain quasi-static problem of a rectangular elastic block subjected to a central downward traction. The computational domain is $40 \ \text{cm} \times 40 \ \text{cm}$, within which a block of dimensions $20 \ \text{cm} \times 10 \ \text{cm}$ is immersed centrally in the fluid domain. A downward traction of magnitude $200\ \text{dyn}/\text{cm}$ is applied to the central $10$ cm of the block's top surface. Figure~\ref{fig:compresed_block_setup} illustrates this loading configuration and the structural dimensions.

The discretization uses a uniform Cartesian grid with $N = \lceil M \cdot \mfac \rceil$ cells per direction, where $M$ is the number of $\Qone$ elements along the longest edge of the Lagrangian mesh. We examine cases with $M = 4, 8, 16, 32, 48,$ and $64$ elements.

The structure follows a neo-Hookean material model, with physical and numerical parameters specified in Table~\ref{tb:cb-param}. The traction load increases linearly in time until reaching full magnitude at $T_{\text{l}} = 40.0$ s, using a time step size of $\Delta t = 0.001 \euleriandx\ \text{s}$. The simulation runs until $T_\text{f} = 100.0$ s to ensure equilibrium is achieved.

\begin{table}[H]
\centering
\caption{Parameters for the compressed block benchmark }
\begin{tabular}{ l  l  l  c  }
\hline
Quantity & Symbol & Value & Unit\\
\hline
Density & $\rho$ & $1.0$ & $\frac{\text{g}}{\text{cm}^3}$\\
Viscosity & $\mu$ & $0.16$ & $\frac{\text{dyn} \cdot \text{s}}{\text{cm}^2}$ \\
Shear modulus & $G$ & $80.194$ & $\frac{\text{dyn}}{\text{cm}^2}$  \\
Numerical bulk modulus & $\kappas$ & $374.239$ &
$\frac{\text{dyn}}{\text{cm}^2}$\\
Final time & $T_\text{f}$ & $100.0$  & s\\
Load time & $T_{\text{l}}$ & $40.0$ & s \\
\hline
\end{tabular}
\label{tb:cb-param}
\end{table}

The boundary conditions for the structure include zero vertical displacements along the bottom boundary and zero horizontal displacements along the top boundary, with zero traction imposed on all other boundaries. To enforce the bottom boundary constraint, we use a penalty parameter $\kappa_{\text{S}} = 2.5 \cdot \frac{2.5\Delta x}{\Delta t} \frac{\text{dyn}}{\text{cm}^3}$.

The primary quantities of interest are the element Jacobians ($J$), which characterize the volumes of the deformed elements, and the vertical displacement ($\Delta Y$) at the center of the top surface.

%

\begin{figure}[H]
  \centering
  \begin{subfigure}[b]{0.45\textwidth}
    \centering
    \includegraphics[width=\textwidth]{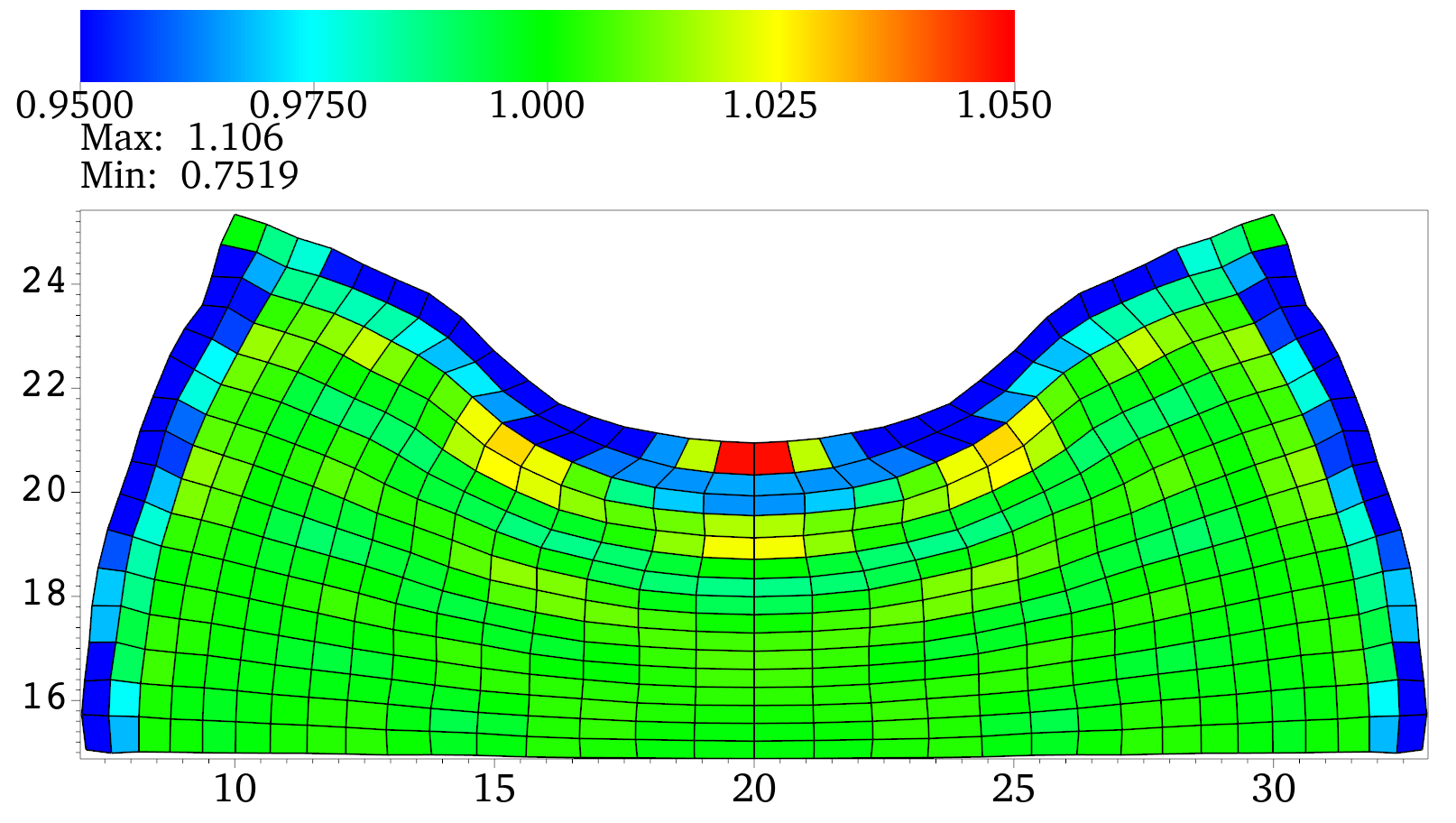}
    \caption{$\text{IB}_3$}
    \label{fig:compressed_block_IB3_no_stab}
  \end{subfigure}
  \hfill
  \begin{subfigure}[b]{0.45\textwidth}
    \centering
    \includegraphics[width=\textwidth]{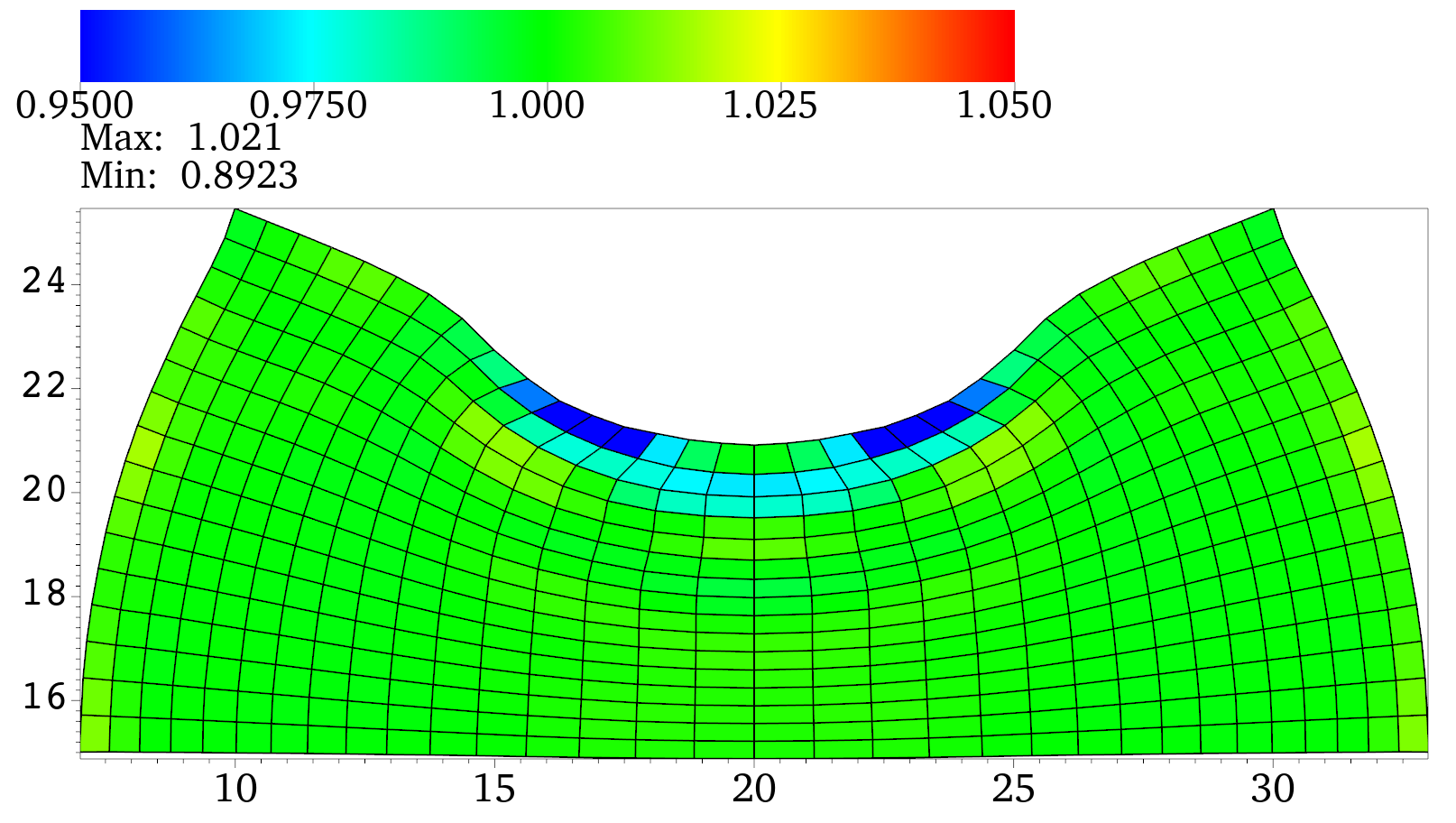}
    \caption{$\text{IB}_3$, $\nu=0.4$, modified invariants}
    \label{fig:compressed_block_IB3_stab}
  \end{subfigure}
  
  \begin{subfigure}[b]{0.45\textwidth}
    \centering
    \includegraphics[width=\textwidth]{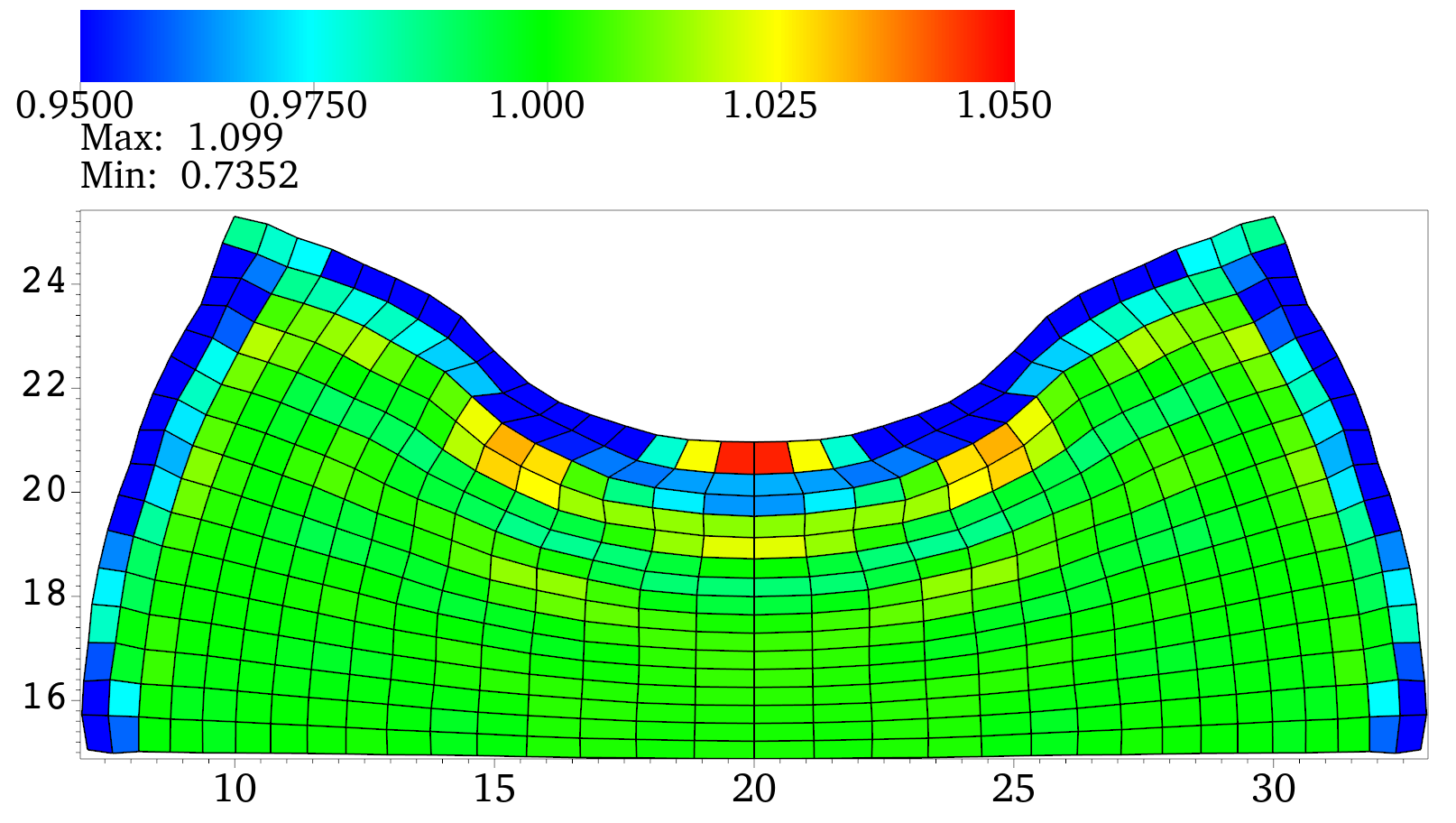}
    \caption{$\text{BS}_3$}
    \label{fig:compressed_block_BS3_no_stab}
  \end{subfigure}
  \hfill
  \begin{subfigure}[b]{0.45\textwidth}
    \centering
    \includegraphics[width=\textwidth]{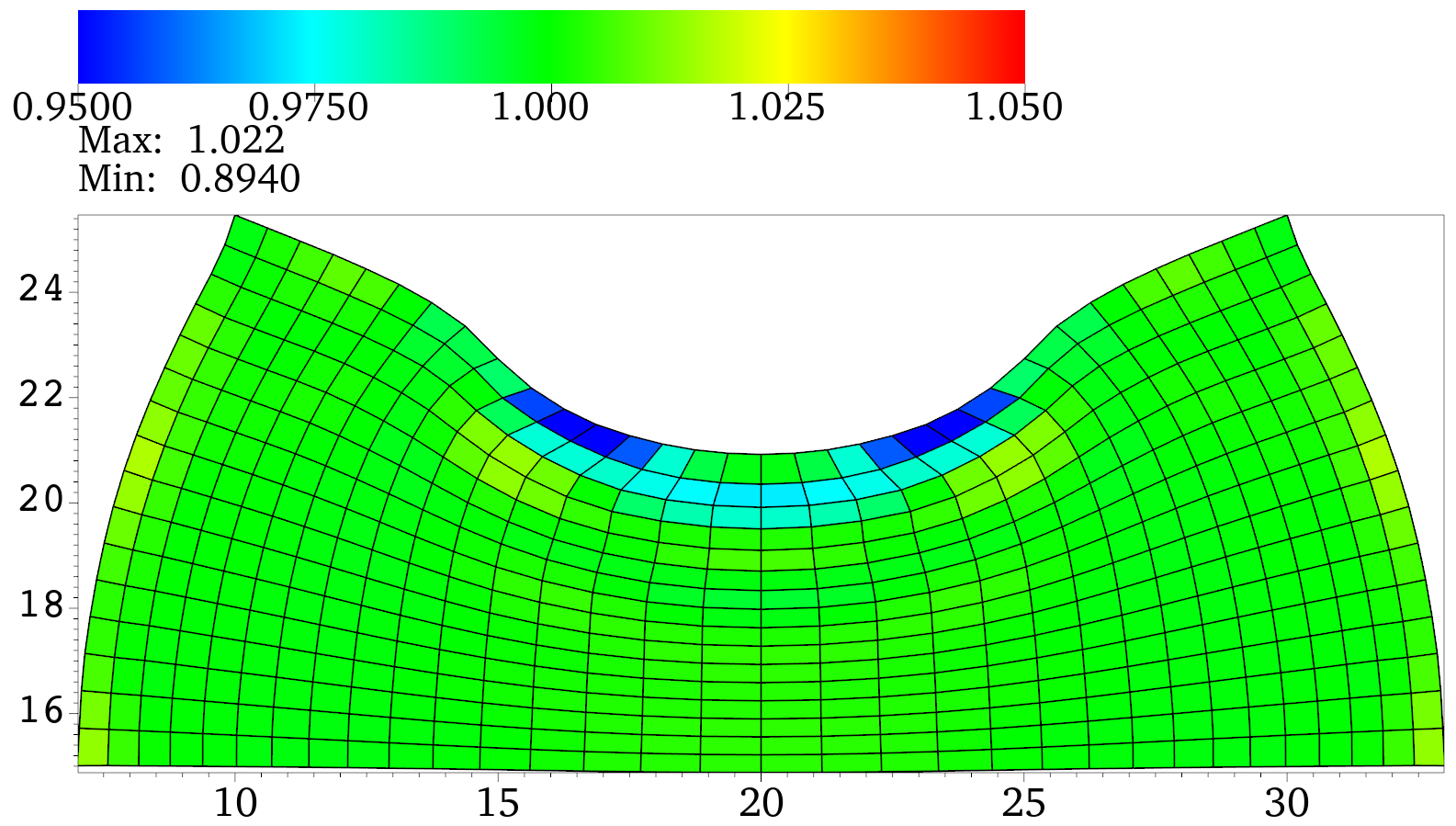}
    \caption{$\text{BS}_3$, $\nu=0.4$, modified invariants}
    \label{fig:compressed_block_BS3_stab}
  \end{subfigure}  
  
  \begin{subfigure}[b]{0.45\textwidth}
    \centering
    \includegraphics[width=\textwidth]{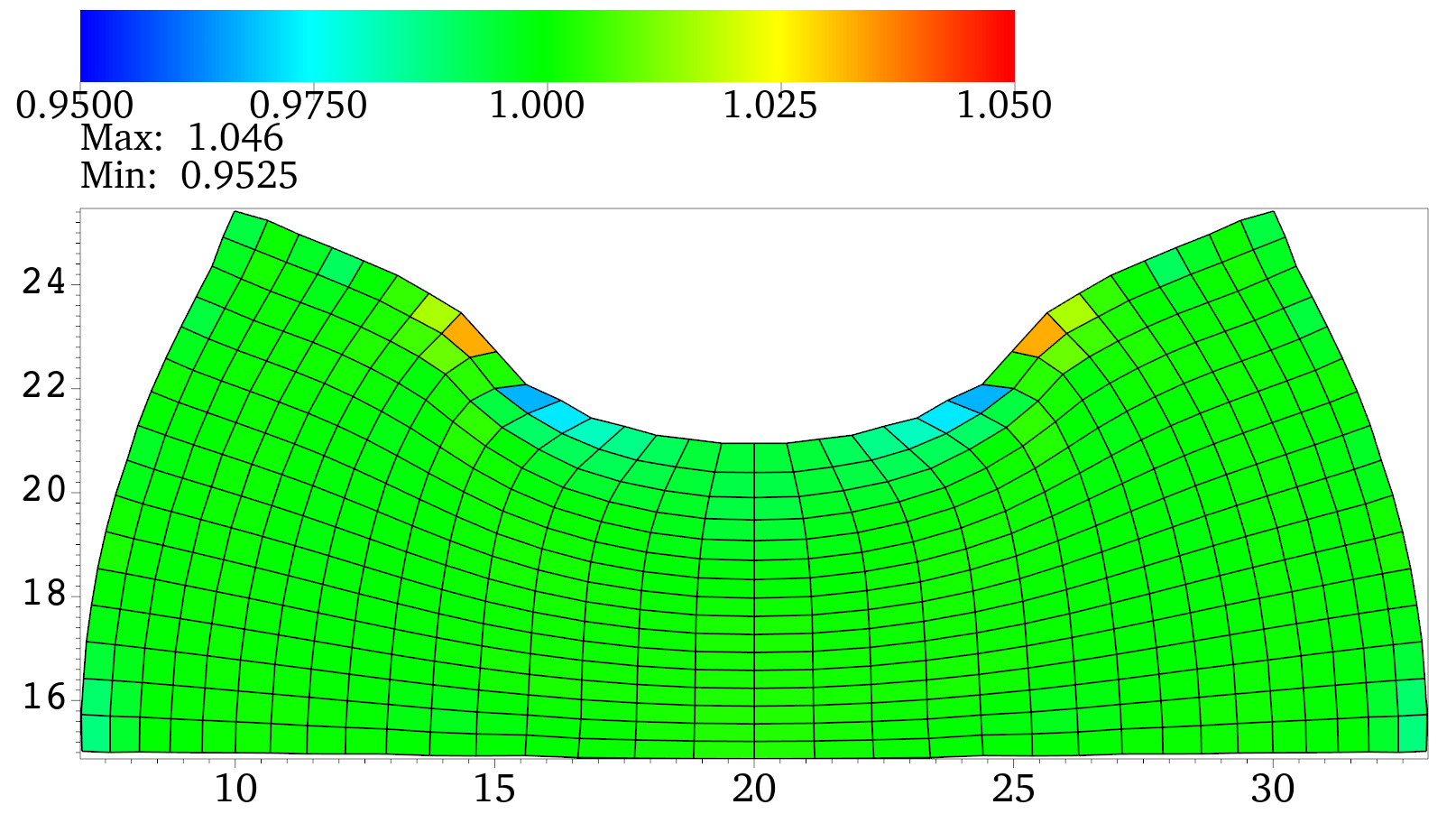}
    \caption{$\text{CBS}_{32}$}
    \label{fig:compressed_block_CBS32_no_stab}
  \end{subfigure}
  \hfill
  \begin{subfigure}[b]{0.45\textwidth}
    \centering
    \includegraphics[width=\textwidth]{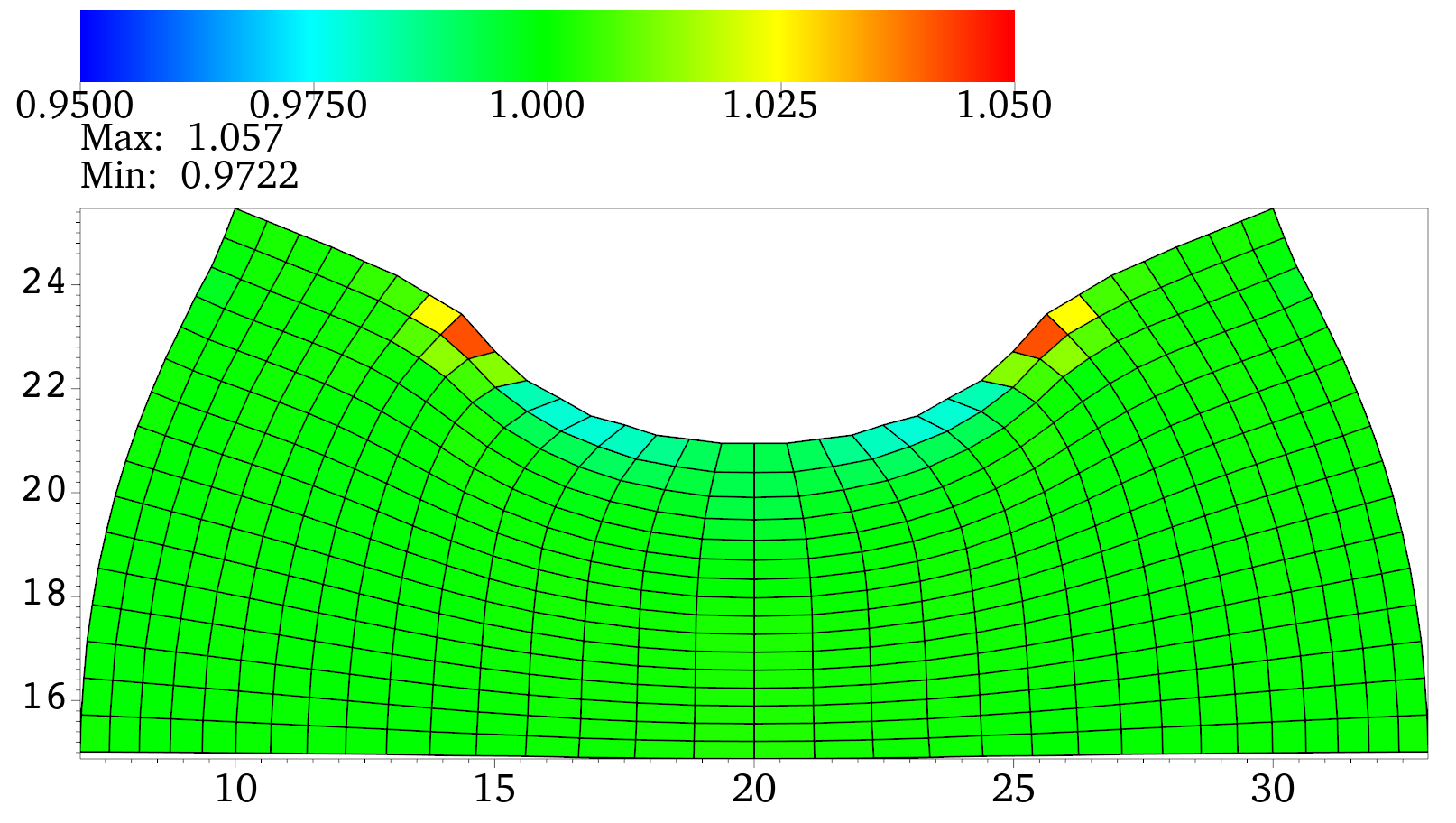}
    \caption{$\text{CBS}_{32}$, $\nu=0.4$, modified invariants}
    \label{fig:compressed_block_CBS32_stab}
  \end{subfigure}    
  
  \begin{subfigure}[b]{0.45\textwidth}
    \centering
    \includegraphics[width=\textwidth]{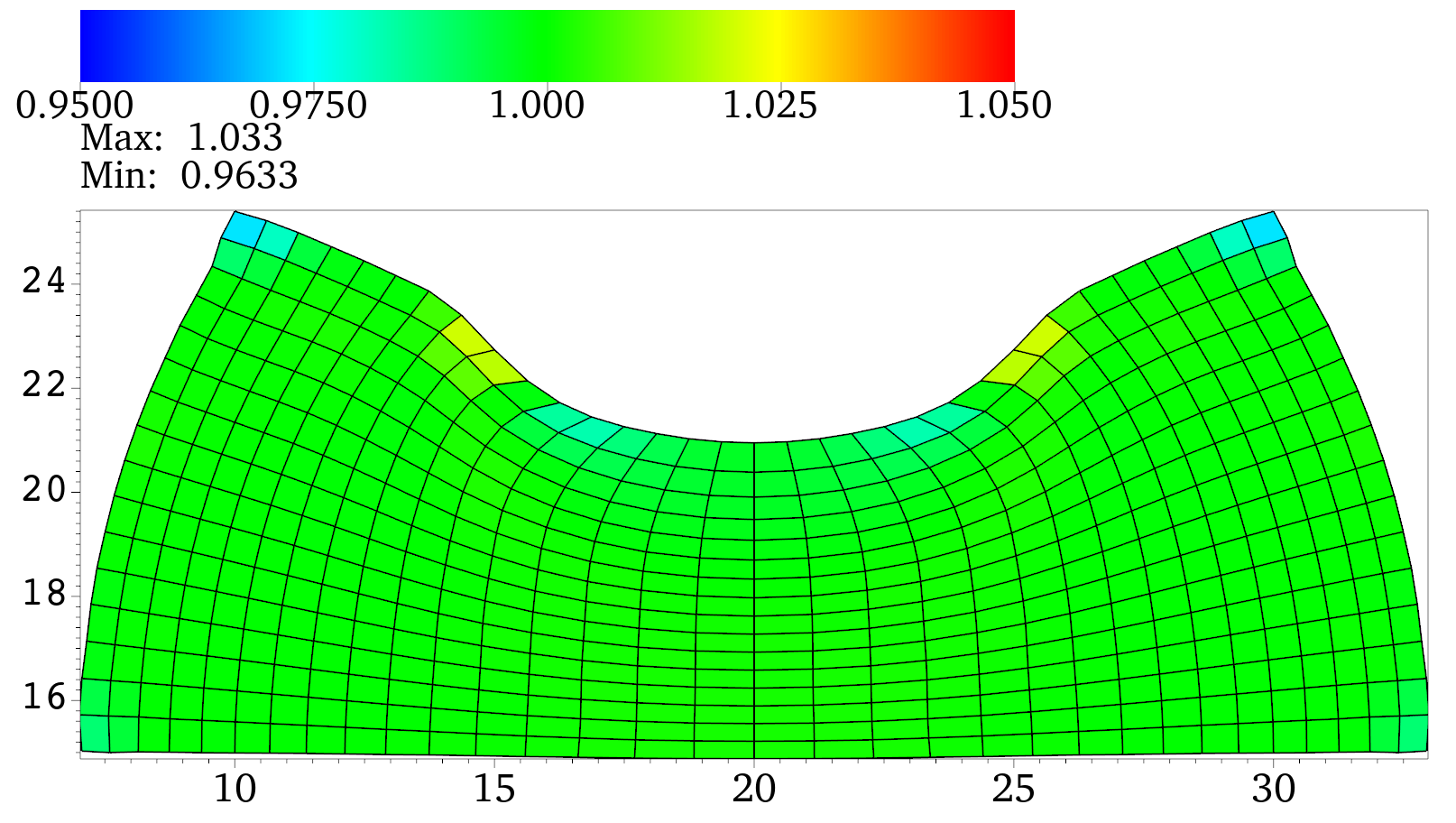}
    \caption{$\text{CBS}_{43}$}
    \label{fig:compressed_block_CBS43_no_stab}
  \end{subfigure}
  \hfill
  \begin{subfigure}[b]{0.45\textwidth}
    \centering
    \includegraphics[width=\textwidth]{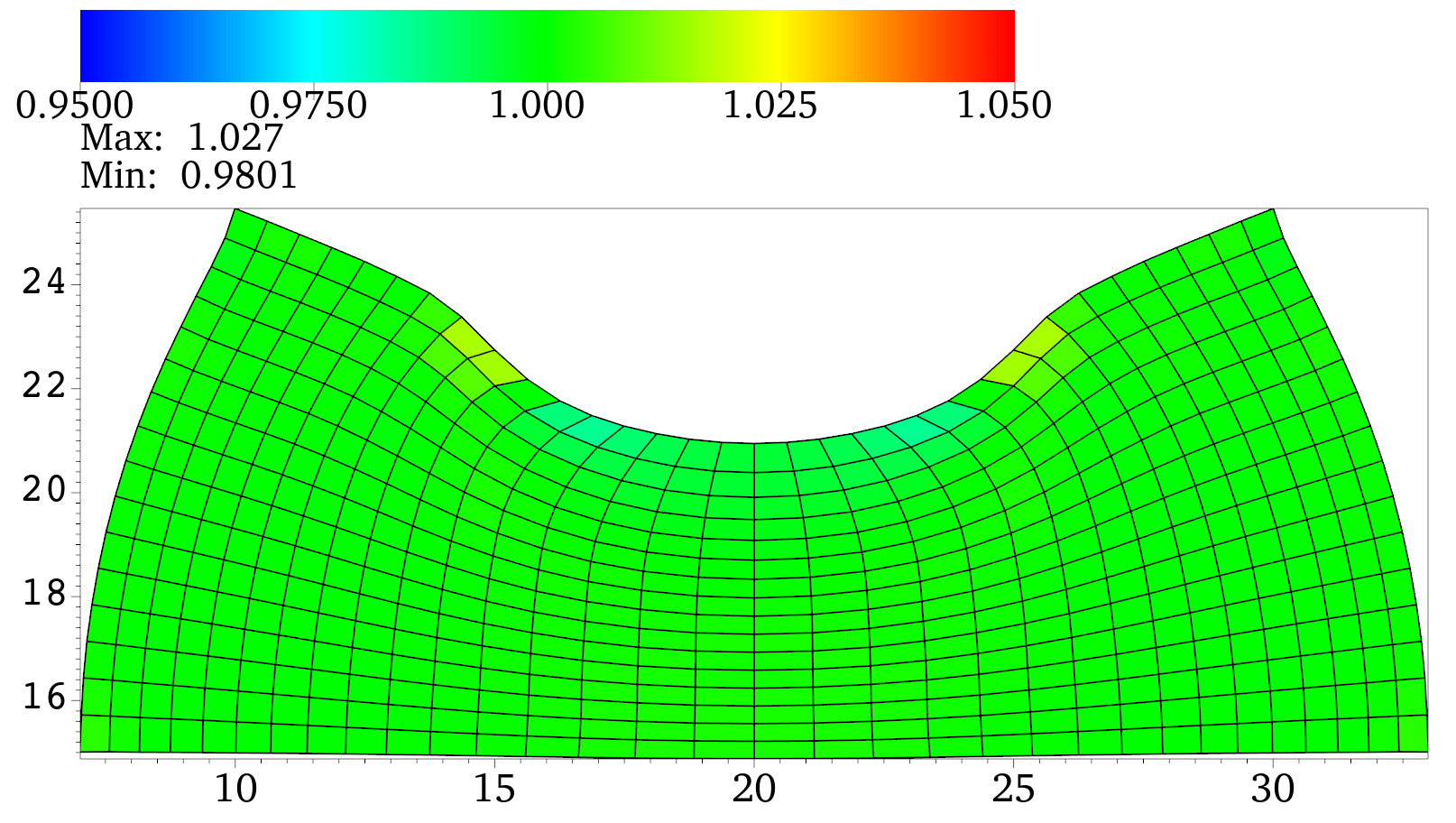}
    \caption{$\text{CBS}_{43}$, $\nu=0.4$, modified invariants}
    \label{fig:compressed_block_CBS43_stab}
  \end{subfigure}    
  
  \caption{The comparison of the Jocobians between the unmodified invariants with omitting the volumetric energy (left column) and with both treatments (right column) for different kernels ($t$ = 100 s). The special treatments improve $\text{IB}_3$ and $\text{BS}_3$ a lot for volume conservation but have little effect on $\text{CBS}_{32}$ and $\text{CBS}_{43}$.}
  \label{fig:compressed_block_jacobian_stab_no_stab}
\end{figure}

Fig. \ref{fig:compressed_block_jacobian_stab_no_stab} compares the Jacobians for simulations with and without stabilization treatments (unmodified invariants and no volumetric energy versus both treatments) across different kernels. While stabilization treatments substantially improve volume conservation for IB ($\text{IB}_3$) and BS ($\text{BS}_3$) kernels, they have minimal effect on CBS kernels ($\text{CBS}_{32}$ and $\text{CBS}_{43}$). This demonstrates that CBS kernels inherently maintain incompressibility without requiring volumetric energy terms or modified invariants. Although the modified invariants and volumetric penalty treatments provide slight improvements in the reduction of spurious volume errors for CBS kernels, particularly at the upper corners, these effects are much less significant compared to their impact on isotropic kernels.

\begin{figure}[H]
  \centering
  \begin{subfigure}[b]{0.49\textwidth}
    \centering
    \includegraphics[width=\textwidth]{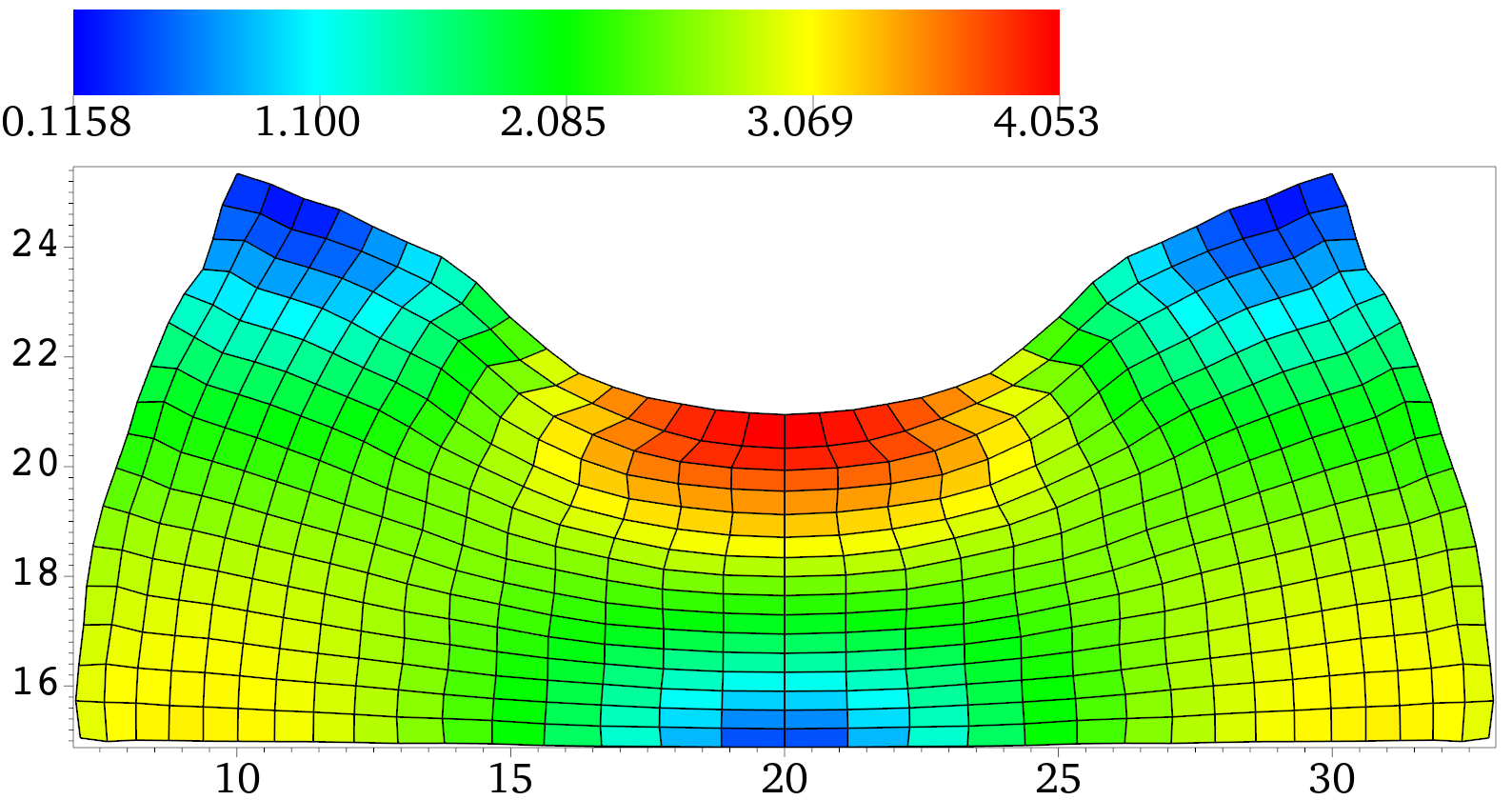}
    \caption{$\text{IB}_3$}
    \label{fig:compressed_block_disp_IB3_no_stab}
  \end{subfigure}
  \hfill
  \begin{subfigure}[b]{0.49\textwidth}
    \centering
    \includegraphics[width=\textwidth]{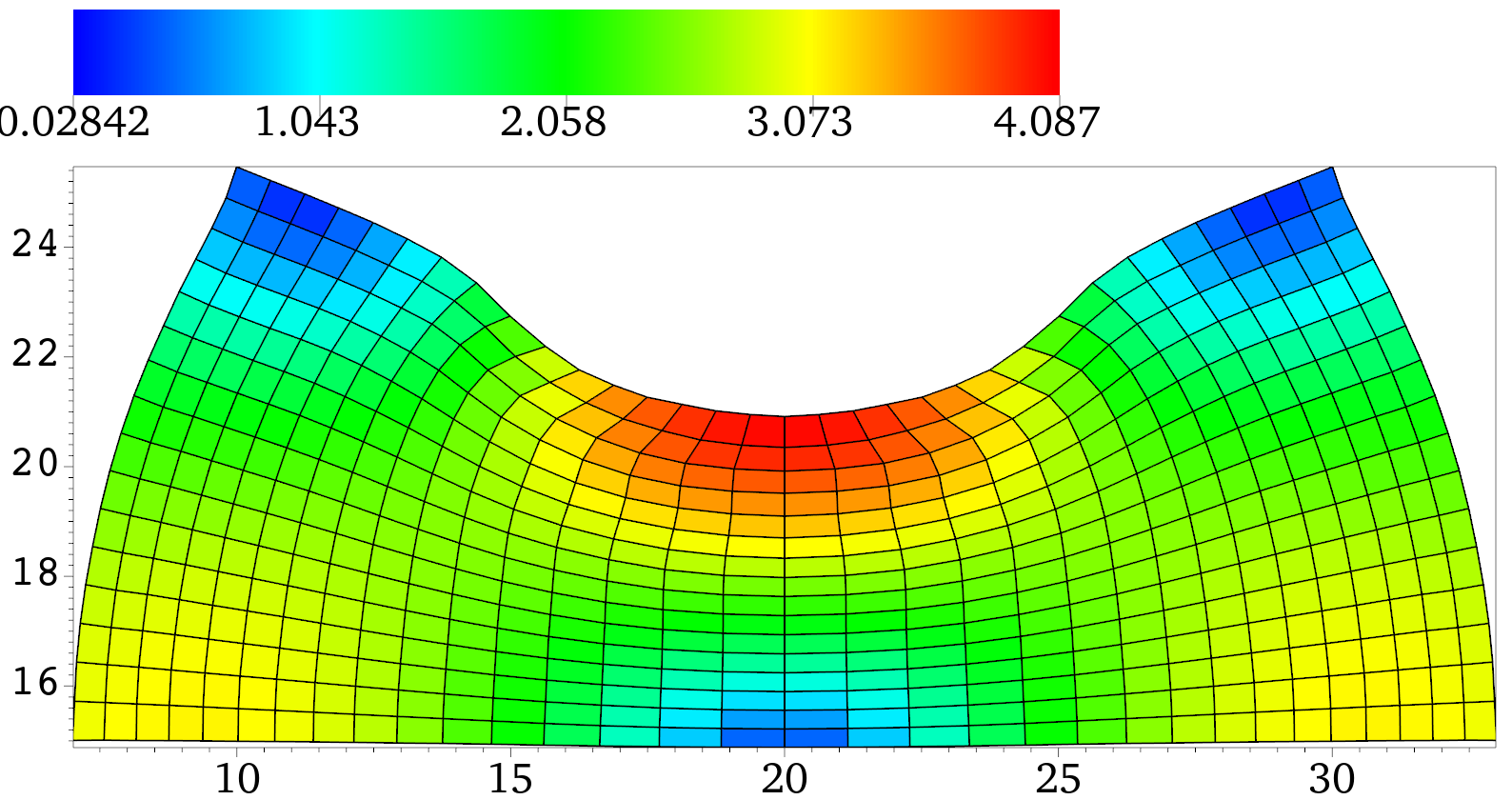}
    \caption{$\text{IB}_3$, $\nu=0.4$, modified invariants}
    \label{fig:compressed_block_disp_IB3_stab}
  \end{subfigure}
  
  \begin{subfigure}[b]{0.49\textwidth}
    \centering
    \includegraphics[width=\textwidth]{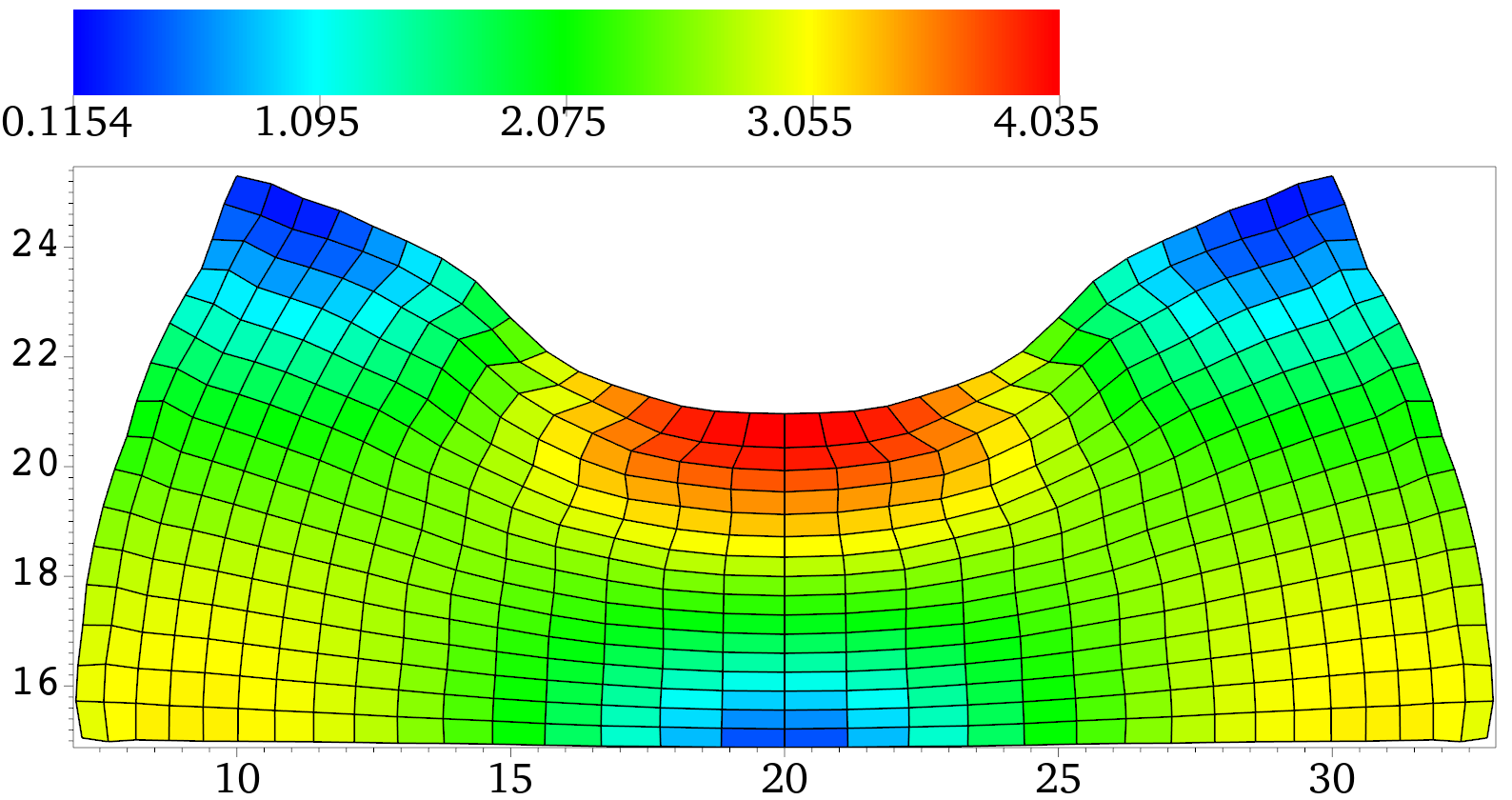}
    \caption{$\text{BS}_3$}
    \label{fig:compressed_block_disp_BS3_no_stab}
  \end{subfigure}
  \hfill
  \begin{subfigure}[b]{0.49\textwidth}
    \centering
    \includegraphics[width=\textwidth]{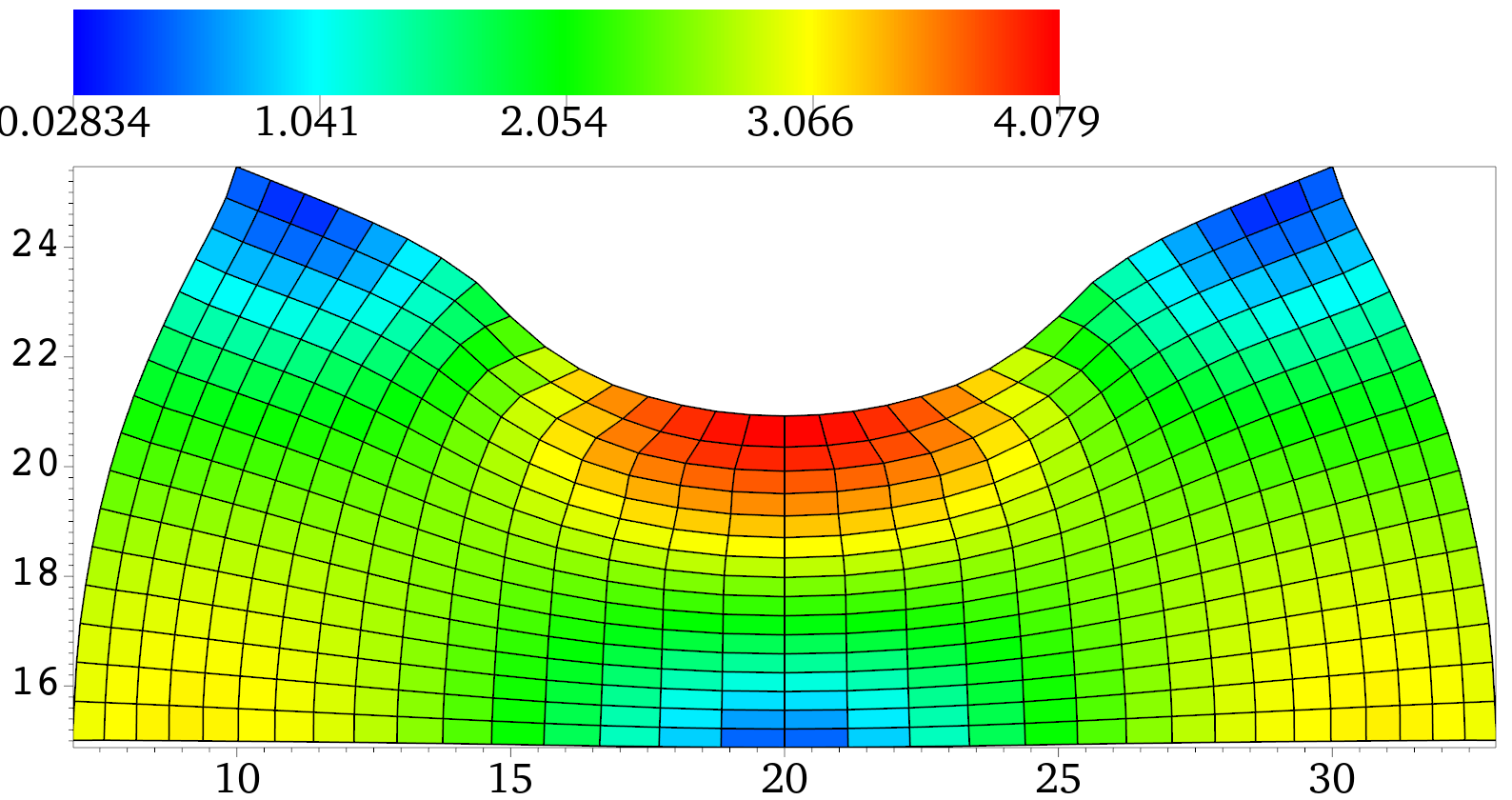}
    \caption{$\text{BS}_3$, $\nu=0.4$, modified invariants}
    \label{fig:compressed_block_disp_BS3_stab}
  \end{subfigure}  
  
  \begin{subfigure}[b]{0.49\textwidth}
    \centering
    \includegraphics[width=\textwidth]{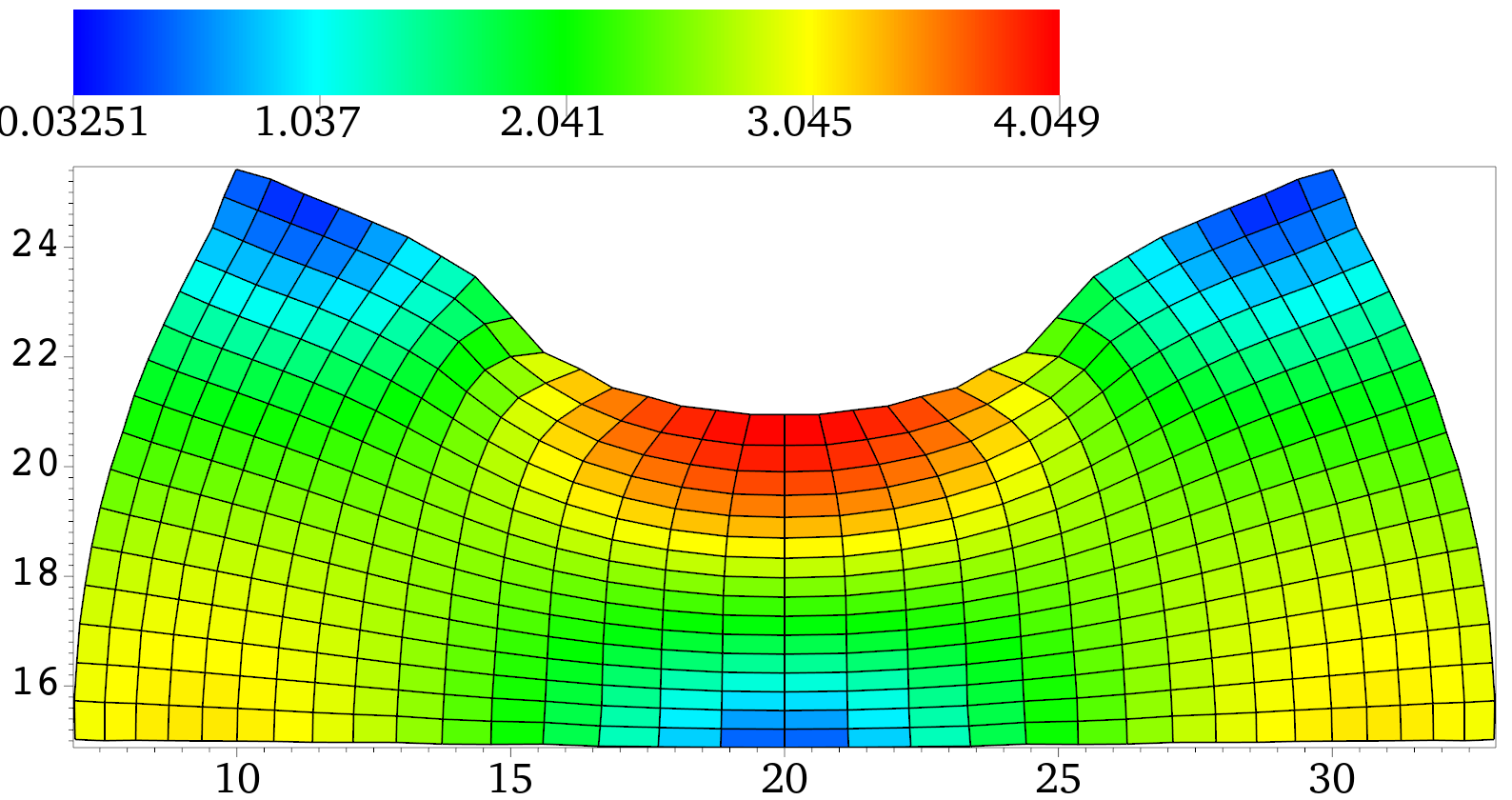}
    \caption{$\text{CBS}_{32}$}
    \label{fig:compressed_block_disp_CBS32_no_stab}
  \end{subfigure}
  \hfill
  \begin{subfigure}[b]{0.49\textwidth}
    \centering
    \includegraphics[width=\textwidth]{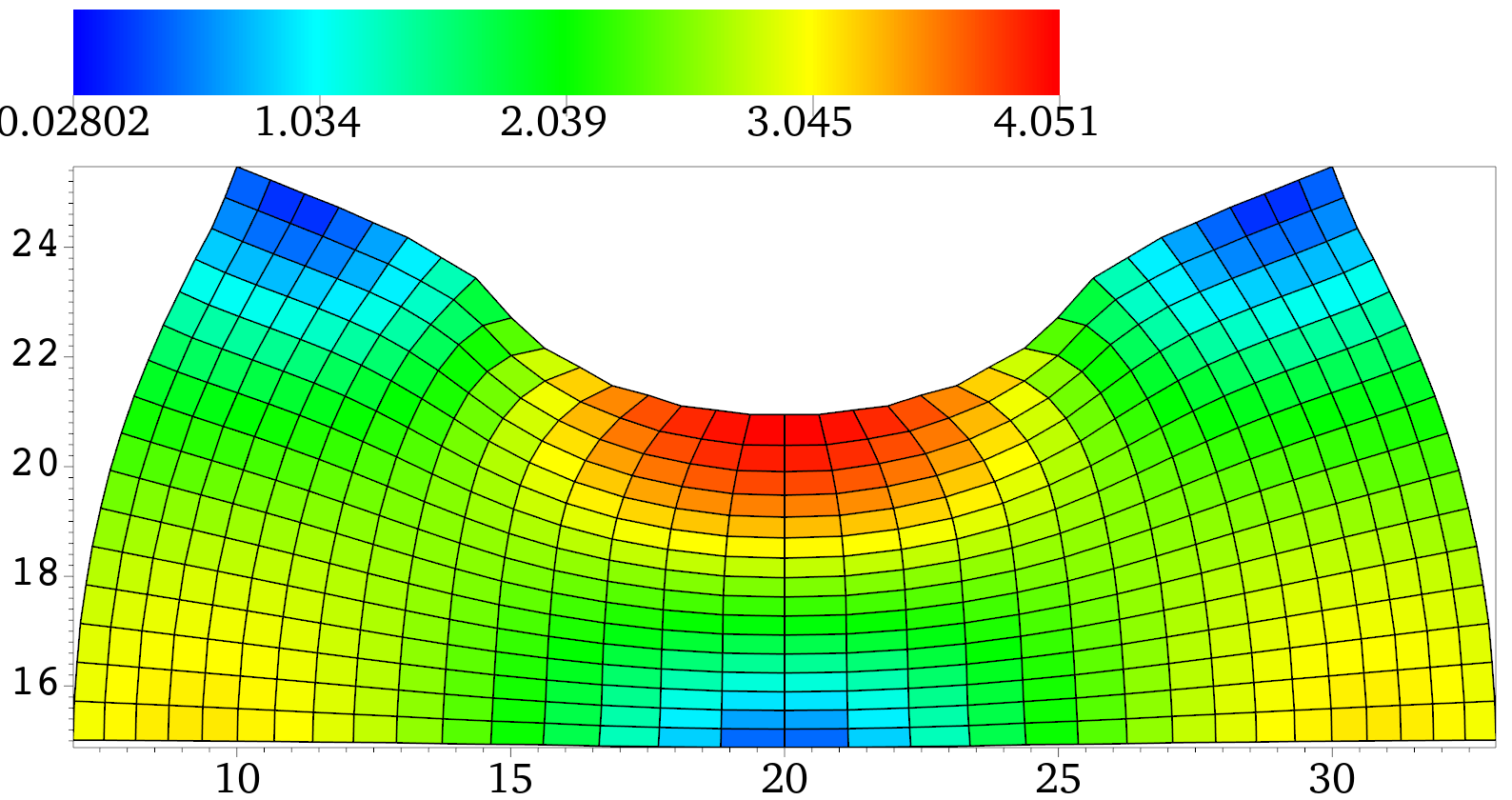}
    \caption{$\text{CBS}_{32}$, $\nu=0.4$, modified invariants}
    \label{fig:compressed_block_disp_CBS32_stab}
  \end{subfigure}    
  
  \begin{subfigure}[b]{0.49\textwidth}
    \centering
    \includegraphics[width=\textwidth]{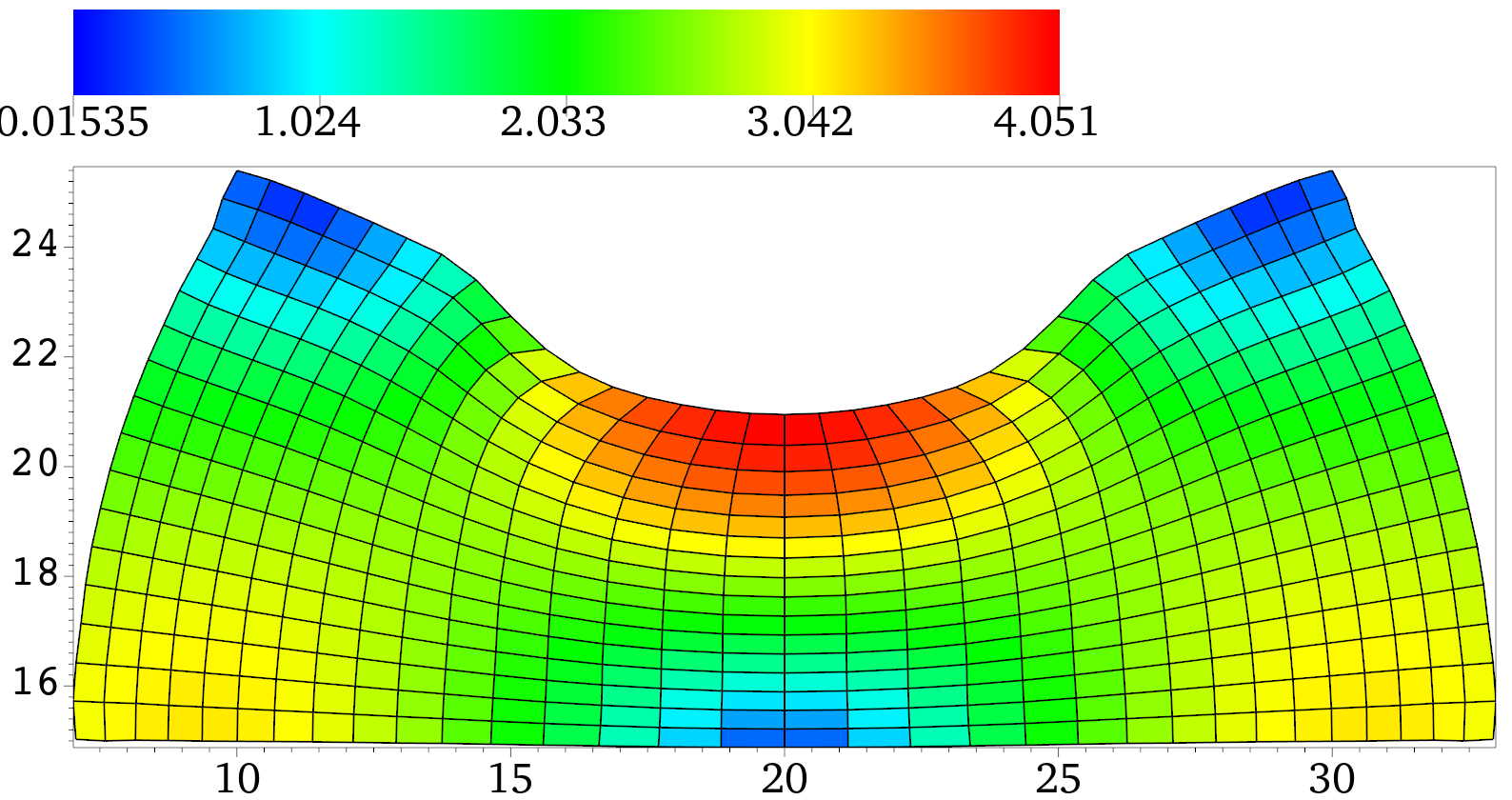}
    \caption{$\text{CBS}_{43}$}
    \label{fig:compressed_block_disp_CBS43_no_stab}
  \end{subfigure}
  \hfill
  \begin{subfigure}[b]{0.49\textwidth}
    \centering
    \includegraphics[width=\textwidth]{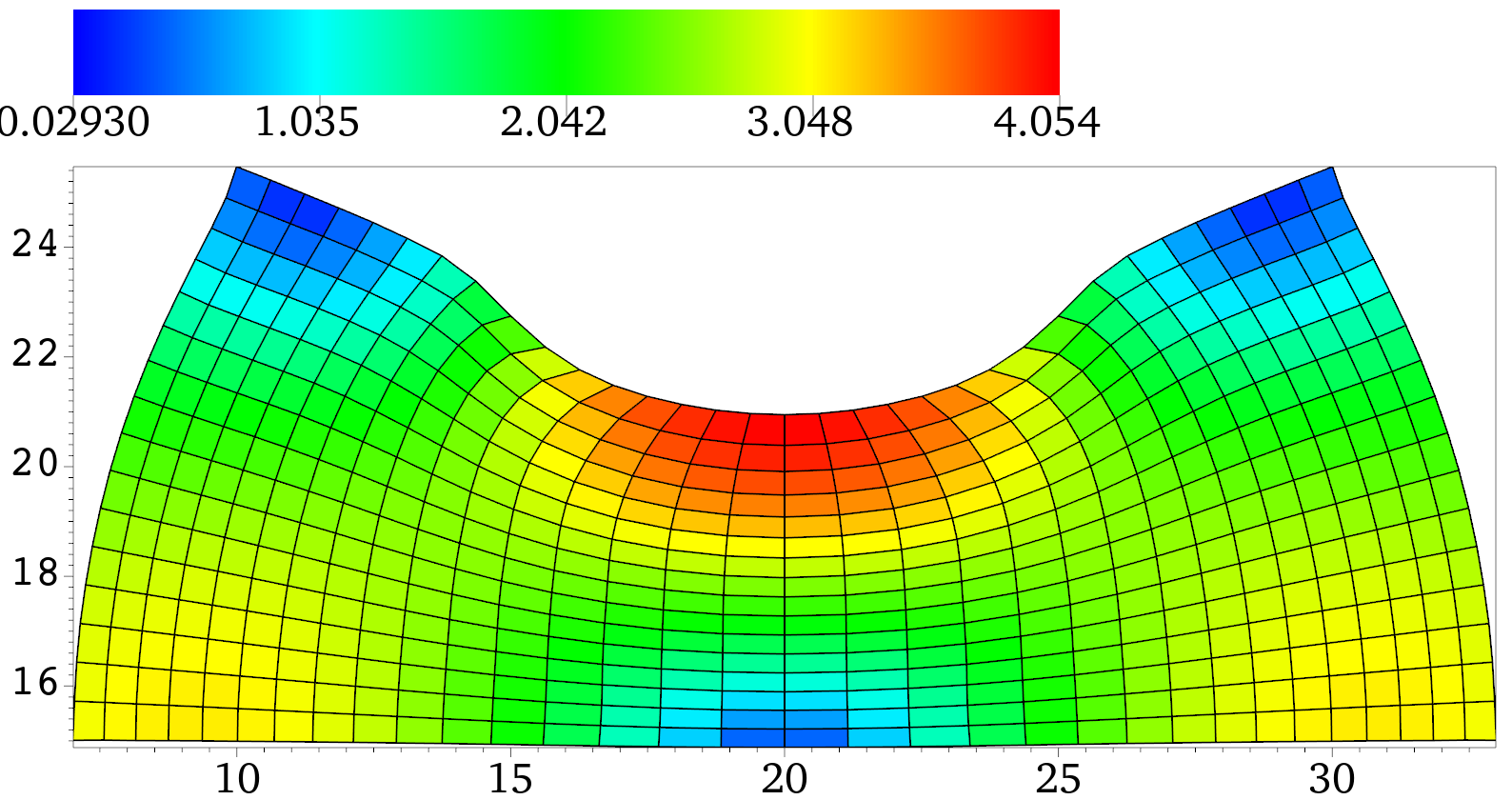}
    \caption{$\text{CBS}_{43}$, $\nu=0.4$, modified invariants}
    \label{fig:compressed_block_disp_CBS43_stab}
  \end{subfigure}    
  
  \caption{The comparison of the displacement between the unmodified invariants with omitting the volumetric energy (left column) and with both treatments (right column) for different kernels. The special treatments improve $\text{IB}_3$ and $\text{BS}_3$ a lot for volume conservation but have little effect on $\text{CBS}_{32}$ and $\text{CBS}_{43}$.}
  \label{fig:compressed_block_disp_jacobian_stab_no_stab}
\end{figure}

\begin{figure}[H]
    \centering
    \includegraphics[width=1\linewidth]{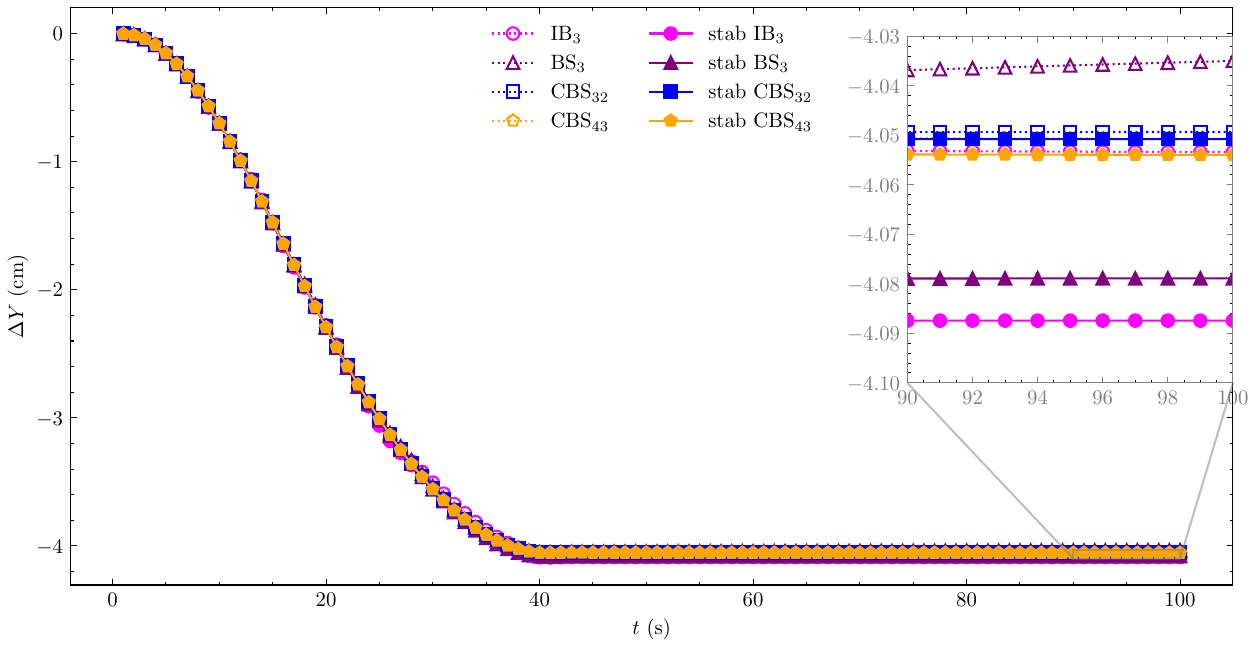}
    \caption{The comparison of the vertical displacements at the center of the top surface with and without the volumetric stabilization treatments for different kernels with $M=32$ and $\MFAC = 0.5$ for the compressed block benchmark. The displacement of the top-mid point is similar for different kernels and treatments, though CBS kernels are much less sensitive to the treatments than IB and BS.}
    \label{fig:compressed_block_point_displacement_stab_no_stab}
\end{figure}

The comparison of displacement fields in Fig.~\ref{fig:compressed_block_disp_jacobian_stab_no_stab} further illustrates these findings. CBS kernels naturally produce smooth displacement fields without requiring volumetric energy or modified invariants for stabilization. In contrast, IB and BS kernels exhibit irregular boundaries when these stabilization treatments are omitted. As shown in Fig.~\ref{fig:compressed_block_point_displacement_stab_no_stab}, whereas CBS kernels show negligible differences between modified and unmodified versions, IB and BS kernels demonstrate significant changes in the vertical displacement at the center of the top surface. With the addition of modified invariants and volumetric stabilization, IB and BS kernels achieve results comparable to those obtained with CBS kernels.

%
%

\begin{figure}[H]
  \centering
  \begin{subfigure}[b]{0.45\textwidth}
    \centering
    \includegraphics[width=\textwidth]{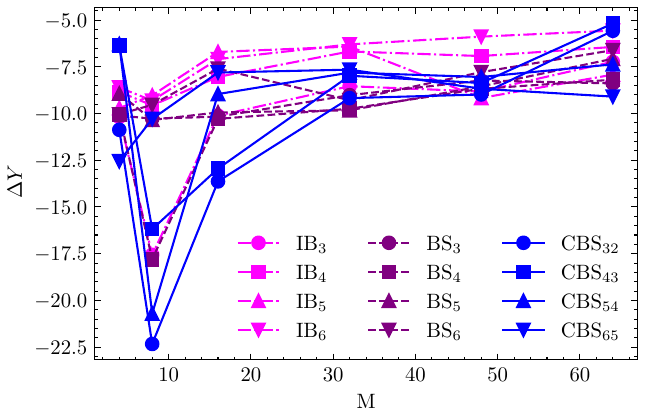}
    \caption{$\MFAC=1.5$}
    \label{fig:compressed_block_mfac=1.5}
  \end{subfigure}
  \hfill
  \begin{subfigure}[b]{0.45\textwidth}
    \centering
    \includegraphics[width=\textwidth]{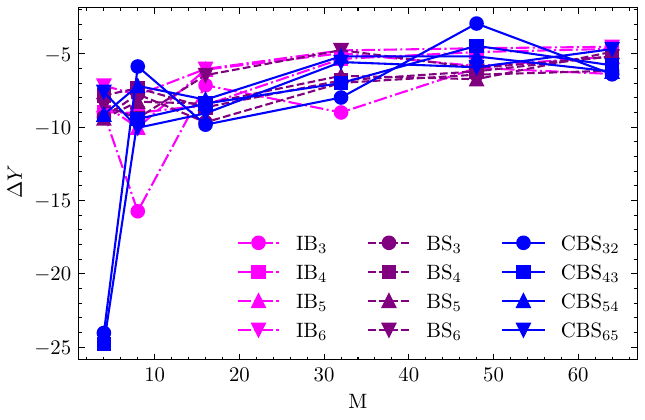}
    \caption{$\MFAC=1.25$}
    \label{fig:compressed_block_mfac=1.25}
  \end{subfigure}
  
  \begin{subfigure}[b]{0.45\textwidth}
    \centering
    \includegraphics[width=\textwidth]{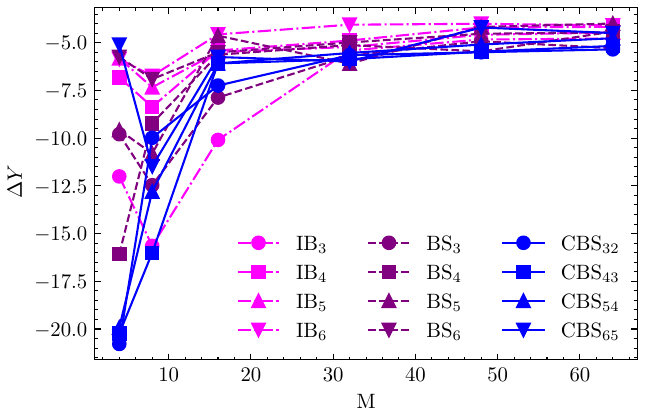}
    \caption{$\MFAC=1.0$}
    \label{fig:compressed_block_mfac=1.0}
  \end{subfigure}
  \hfill
  \begin{subfigure}[b]{0.45\textwidth}
    \centering
    \includegraphics[width=\textwidth]{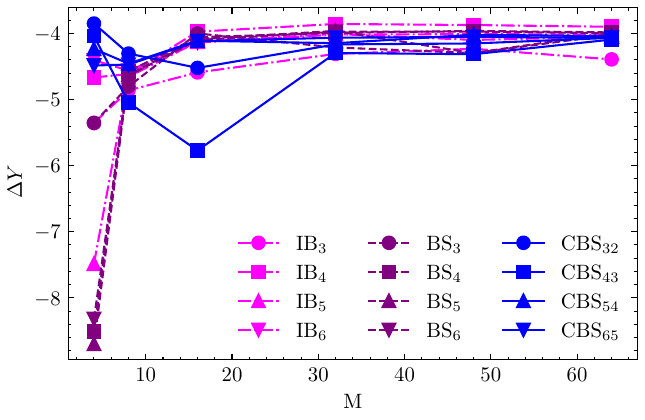}
    \caption{$\MFAC=0.75$}
    \label{fig:compressed_block_mfac=0.75}
  \end{subfigure}  
  
  \begin{subfigure}[b]{0.45\textwidth}
    \centering
    \includegraphics[width=\textwidth]{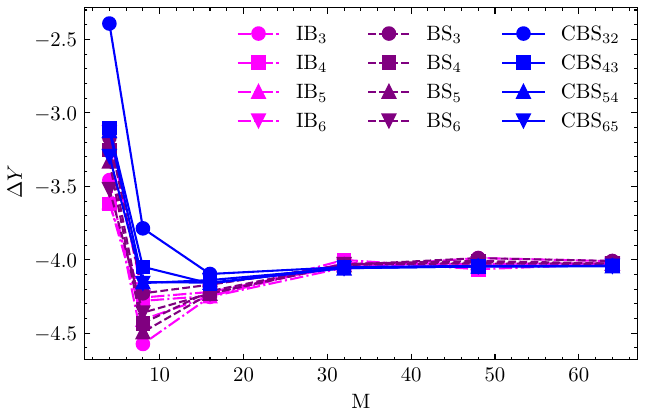}
    \caption{$\MFAC=0.5$}
    \label{fig:compressed_block_mfac=0.5}
  \end{subfigure}
  \hfill
  
  \caption{Grid convergence with different MFAC in terms of the displacement of the probed location. These are the results without modified invariants or volumetric energies. (1) Grid convergence is better for smaller MFACs for all kernels, but CBS kernels converge at a coarse mesh level. For large MFACs, CBS kernels do not converge. (2) Wider kernels give a smaller volume error with the same MFAC.}
  \label{fig:compressed_block/jacobian_stab_no_stab}
\end{figure}

To examine the sensitivity of grid spacing ratios between the solid and fluid meshes (MFAC), Fig. \ref{fig:compressed_block/jacobian_stab_no_stab} presents grid convergence results for different MFAC values based on the displacement at a probed location. Observations indicate that: (1) grid convergence improves with smaller MFAC values across all kernels, with CBS achieving convergence at a coarser mesh level, while larger MFAC values prevent CBS kernels from converging; and (2) wider kernels yield more minor volume errors for the same MFAC.

\begin{figure}[H]
\centering
\includegraphics[width=0.55\textwidth]{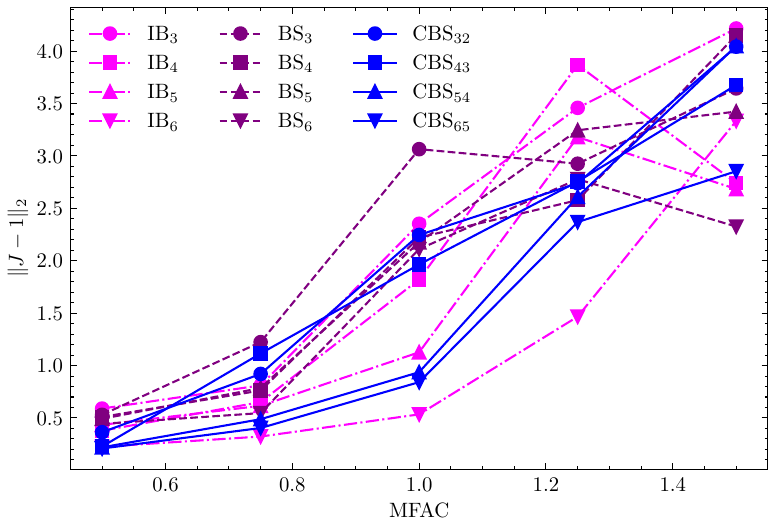}
\caption{Error norms of Jacobian against MFAC with $N$ = 90. All kernel types show convergence under solid mesh refinement, with wider kernels yielding smaller errors. CBS kernels generally demonstrate superior volume conservation across different $\mfac$ values, except for $\text{IB}_6$ which achieves the smallest overall error. CBS kernels show improved performance with decreasing $\mfac$, while IB and BS kernels often perform better at larger $\mfac$ values.}
\label{fig:compressed_block_mfac_J}
\end{figure}

Keeping the fluid mesh resolution fixed at $N=90$, we examine the effect of solid mesh refinement by varying the mesh factor $\mfac$ from $0.5$ to $1.5$ in increments of $0.25$. Figure~\ref{fig:compressed_block_mfac_J} shows the Jacobian error $|J - 1|_2$ for different kernels. The results demonstrate convergence for all kernel types as the solid mesh is refined, with wider kernels consistently producing smaller errors. CBS kernels generally exhibit smaller volume errors across different values of $\mfac$, with $\text{IB}_6$ being the notable exception in achieving the smallest error overall. The relationship between $\mfac$ and volume conservation error shows distinct patterns: CBS kernels consistently decrease errors with decreasing $\mfac$, whereas IB and BS kernels may achieve better volume conservation with larger values of $\mfac$. Specifically, $\text{IB}_4$ and $\text{IB}_5$ exhibit smaller errors at $\mfac = 1.5$ compared to $\mfac = 1.25$, $\text{BS}_6$ shows improved performance at $\mfac = 1.5$ versus $\mfac = 1.25$, and $\text{BS}_3$ achieves better accuracy at $\mfac = 1.25$ than at $\mfac = 1.0$. 

\subsubsection{Cook's Membrane}

\begin{figure}[H]
    \centering
    \includegraphics[width=0.5\linewidth]{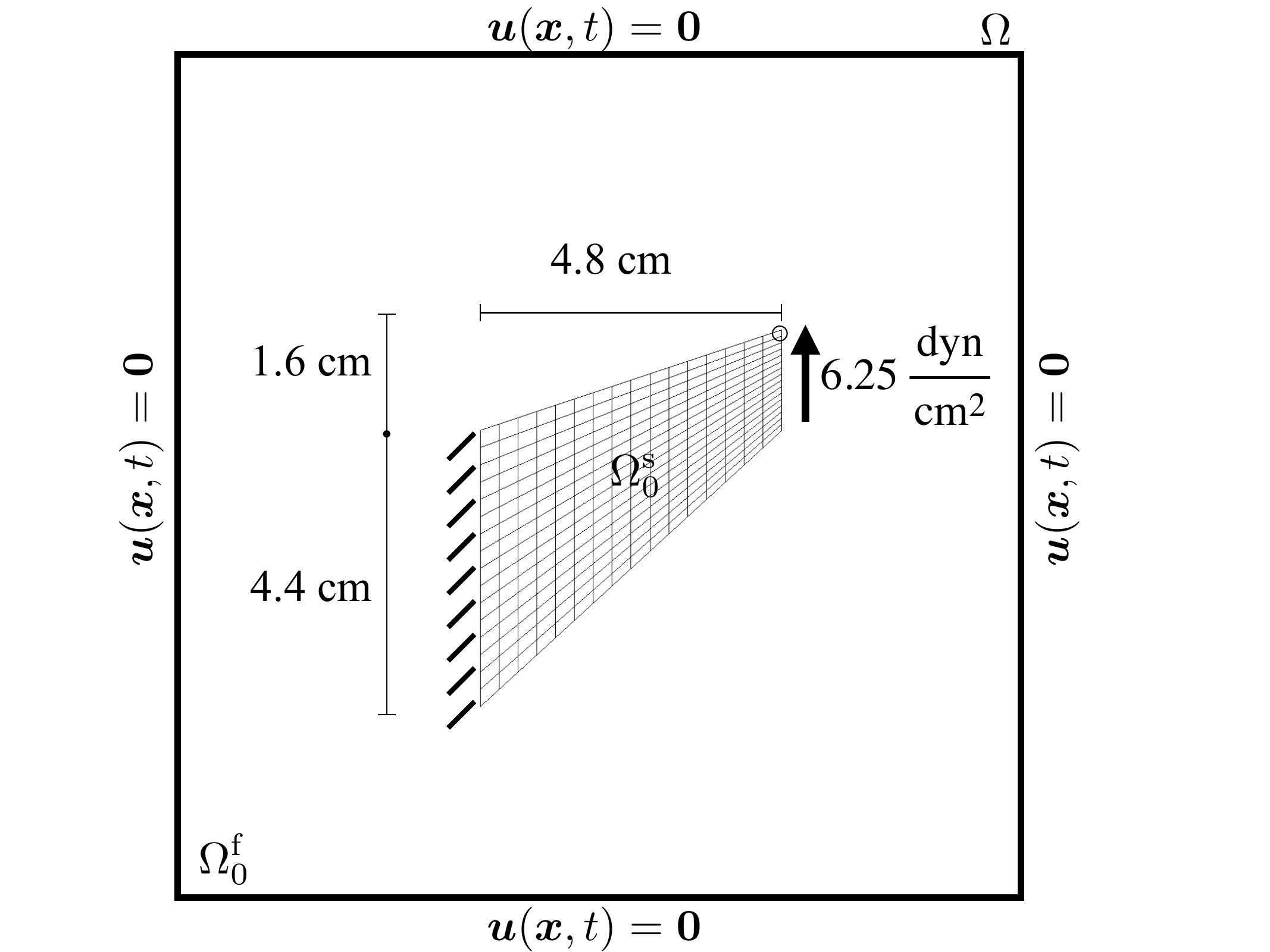}
    \caption{Setup of the Cook's membrane benchmark problem. The $y$-displacement at the upper-right corner (indicated by the circle) is monitored for subsequent analysis. The initial configuration of the structure, denoted by $\Omega_0^s$, is immersed within the fluid domain $\Omega_0^f$. The entire computational domain is $\Omega = \Omega_0^s \cup \Omega_0^f$, with zero fluid velocity enforced on $\partial \Omega$.}
    \label{fig:cook_membrane_problem_setup}
\end{figure}

The Cook's membrane problem is a classical plane strain benchmark widely used to evaluate numerical methods for incompressible elasticity. Following the configuration in \cite{wells2023nodal}, but with a modified computational domain of $13 \text{ cm} \times 13 \text{ cm}$ instead of the original $10 \text{ cm} \times 10 \text{ cm}$ to ensure integer grid resolution during mesh convergence studies, we examine this problem with the structural dimensions and loading conditions illustrated in Fig.~\ref{fig:cook_membrane_problem_setup}.

The discretization employs a Cartesian grid with $N = \lceil M \cdot \mfac \cdot \frac{10}{6.5} \rceil$ cells in each coordinate direction, where $M$ is the number of $\Qone$ elements per edge in the Lagrangian mesh, and the factor $\frac{10}{6.5}$ accounts for the ratio between the computational domain length and the structure's longest side. We examine cases with $M$ = 4, 8, 16, 32, 48, and $64$ elements per side and explore mesh factors $\mfac = 0.5, 0.75, 1.0, 1.25,$ and $1.5$, using a time step size of $\Delta t = 0.001 \euleriandx\ \text{s}$.

The boundary conditions include a fixed left side, enforced by a penalty parameter $\kappa_{\text{S}} = 0.125 \cdot \frac{\Delta x}{\Delta t}\frac{\text{dyn}}{\text{cm}^3}$, and an upward traction of density 6.25 dyn/cm applied to the right side with stress-free conditions on all other boundaries. The traction increases linearly in time until reaching full magnitude at $T_{\text{l}} = 20.0$ s, and the simulation continues until $T_\text{f} = 50.0$ s to ensure equilibrium. Table~\ref{tb:cm-param} lists all relevant physical and numerical parameters.

The primary quantities of interest are the element Jacobians ($J$), and the vertical displacement ($\Delta Y$) at the membrane's upper right corner (located at $(8.05\ \text{cm}, 9.5\ \text{cm})$).

\begin{table}[H]
\centering
\caption{Parameters for the Cook's membrane benchmark}
\begin{tabular}{ l  l  l  c  }
\hline
Quantity & Symbol & Value & Unit\\
\hline
Density & $\rho$ & $1.0$ & $\frac{\text{g}}{\text{cm}^3}$\\
Viscosity & $\mu$ & $0.16$ & $\frac{\text{dyn} \cdot \text{s}}{\text{cm}^2}$ \\
Material model & - & neo-Hookean & - \\
Shear modulus & $G$ & $83.333$ & $\frac{\text{dyn}}{\text{cm}^2}$  \\
Numerical bulk modulus & $\kappas$ & $388.889$
& $\frac{\text{dyn}}{\text{cm}^2}$\\
Final time & $T_\text{f}$ & $50.0$  & s\\
Load time & $T_{\text{l}}$ & $20.0$ & s \\
\hline
\end{tabular}

\label{tb:cm-param}
\end{table}

Fig. \ref{fig:cook_membrane_jacobian_stab_no_stab} compares the Jacobians for simulations using unmodified invariants without volumetric energy (left column) versus simulations with both stabilization treatments (right column) across different kernels. The volume conservation and stabilization treatments significantly enhance volume conservation for IB and BS kernels but have minimal impact on CBS kernels.

\begin{figure}[H]
    \centering
    \includegraphics[width=1\linewidth]{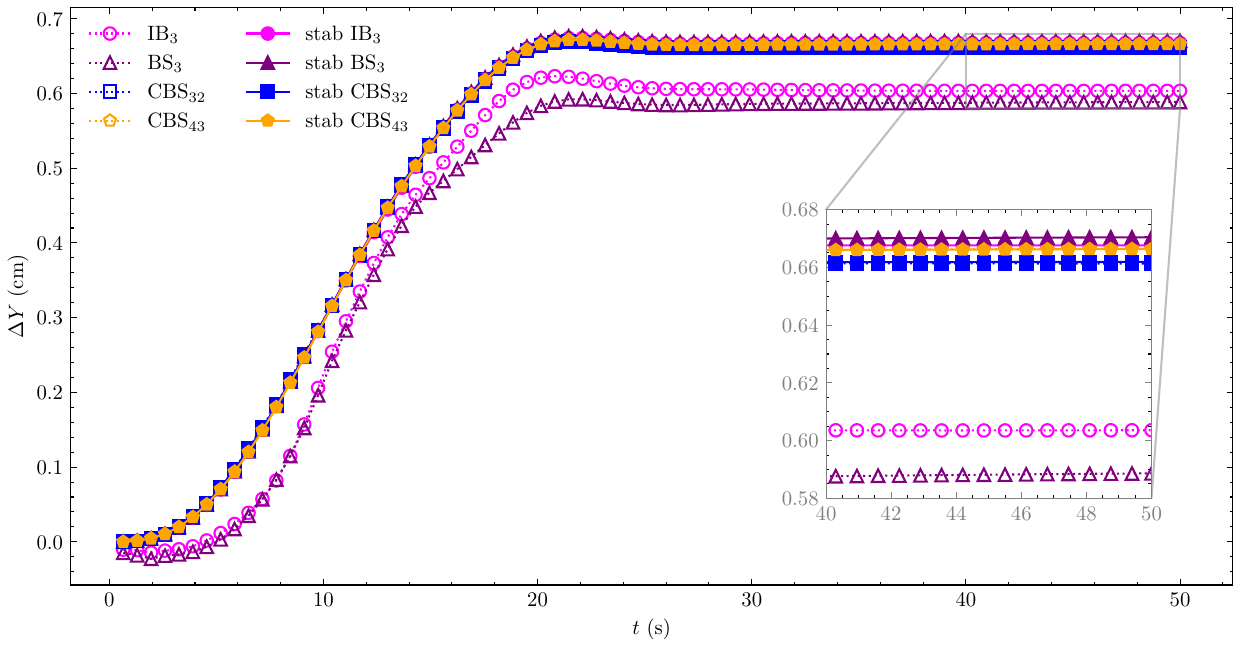}
\caption{Comparison of the displacement at the top-mid point of the compressed block with and without volumetric energy and stabilization treatments for different kernels. A clear difference is observed between cases with and without these treatments for IB and BS kernels, while CBS kernels show no distinguishable difference.}
    \label{fig:cook_membrane_tip_displacement_stab_no_stab}
\end{figure}

This is further illustrated by the comparison of displacement fields in Fig. \ref{fig:cook_membrane_disp_field_stab_no_stab}. CBS kernels produce smoother displacement
fields without needing either volumetric energy or modified invariants stabilization terms. In contrast, when these treatments are omitted, IB and BS kernels display irregular boundaries. Displacements were compared for cases without volumetric energy and with both treatments across different kernels, with $M = 32$ (see Fig. \ref{fig:cook_membrane_tip_displacement_stab_no_stab}). The displacement of the top-mid point is similar across kernels and treatments, although CBS kernels are notably less sensitive to these treatments than IB and BS kernels.

\begin{figure} [H]
  \centering
  \begin{subfigure}{0.49\textwidth}
    \centering
    \includegraphics[width=\textwidth]{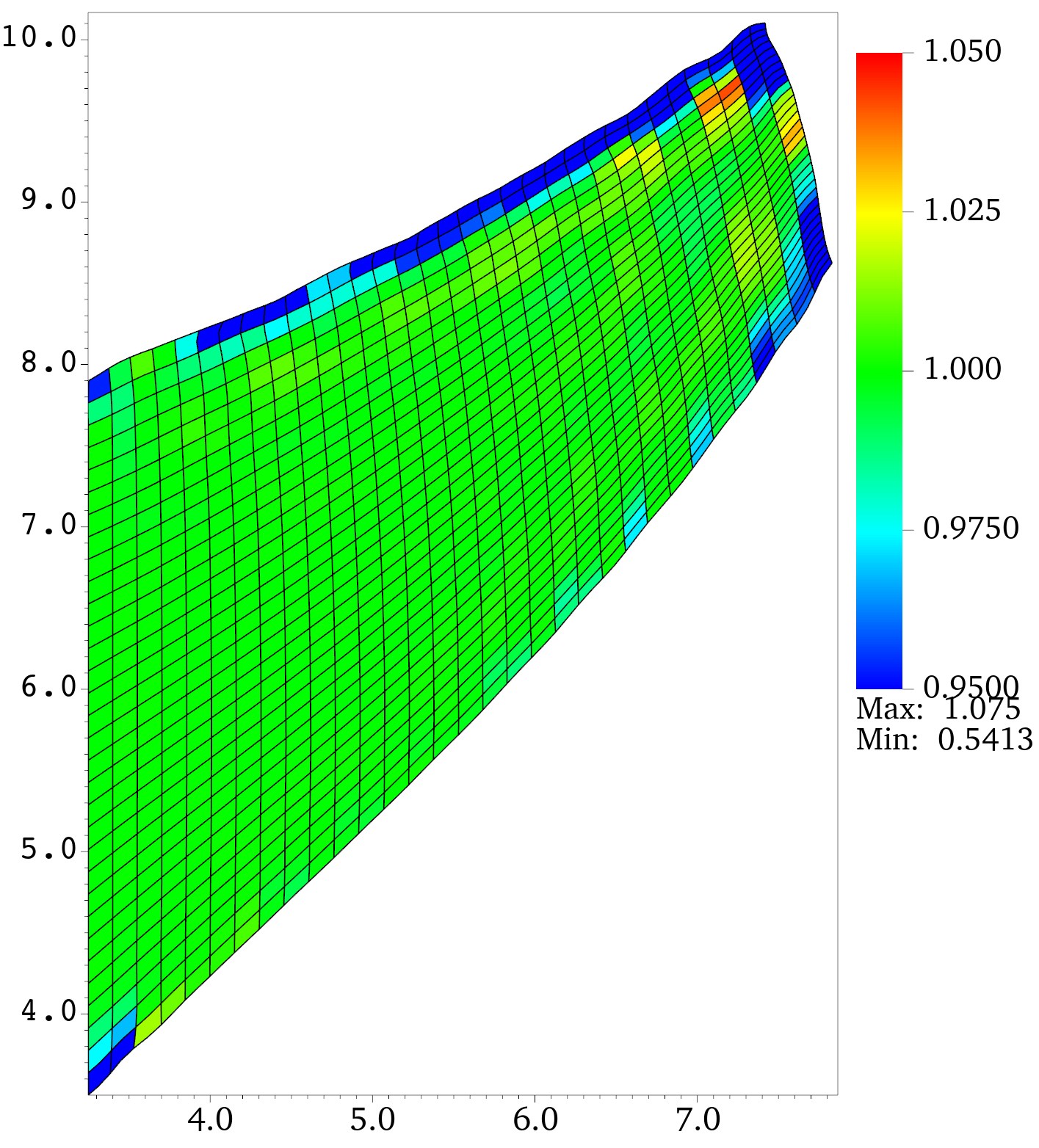}
    \caption{$\text{IB}_3$}
    \label{fig:cook_membrane_IB3_no_stab}
  \end{subfigure}
  \hfill
  \begin{subfigure}{0.49\textwidth}
    \centering
    \includegraphics[width=\textwidth]{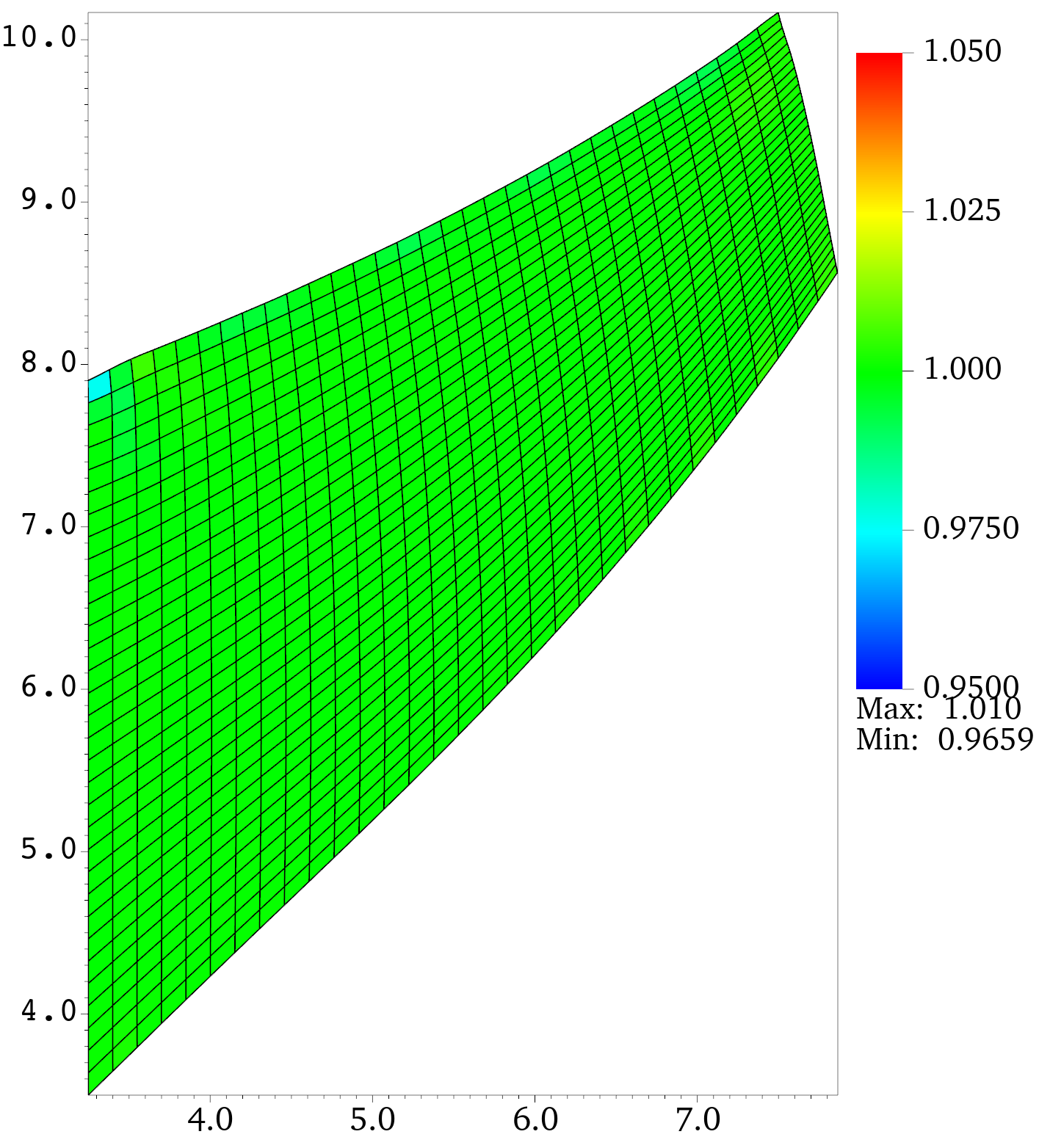}
    \caption{$\text{IB}_3$, $\nu=0.4$, modified invariants}
    \label{fig:cook_membrane_IB3_stab}
  \end{subfigure}
\end{figure}

\begin{figure}[H]\ContinuedFloat  
  \begin{subfigure}{0.49\textwidth}
    \centering
    \includegraphics[width=\textwidth]{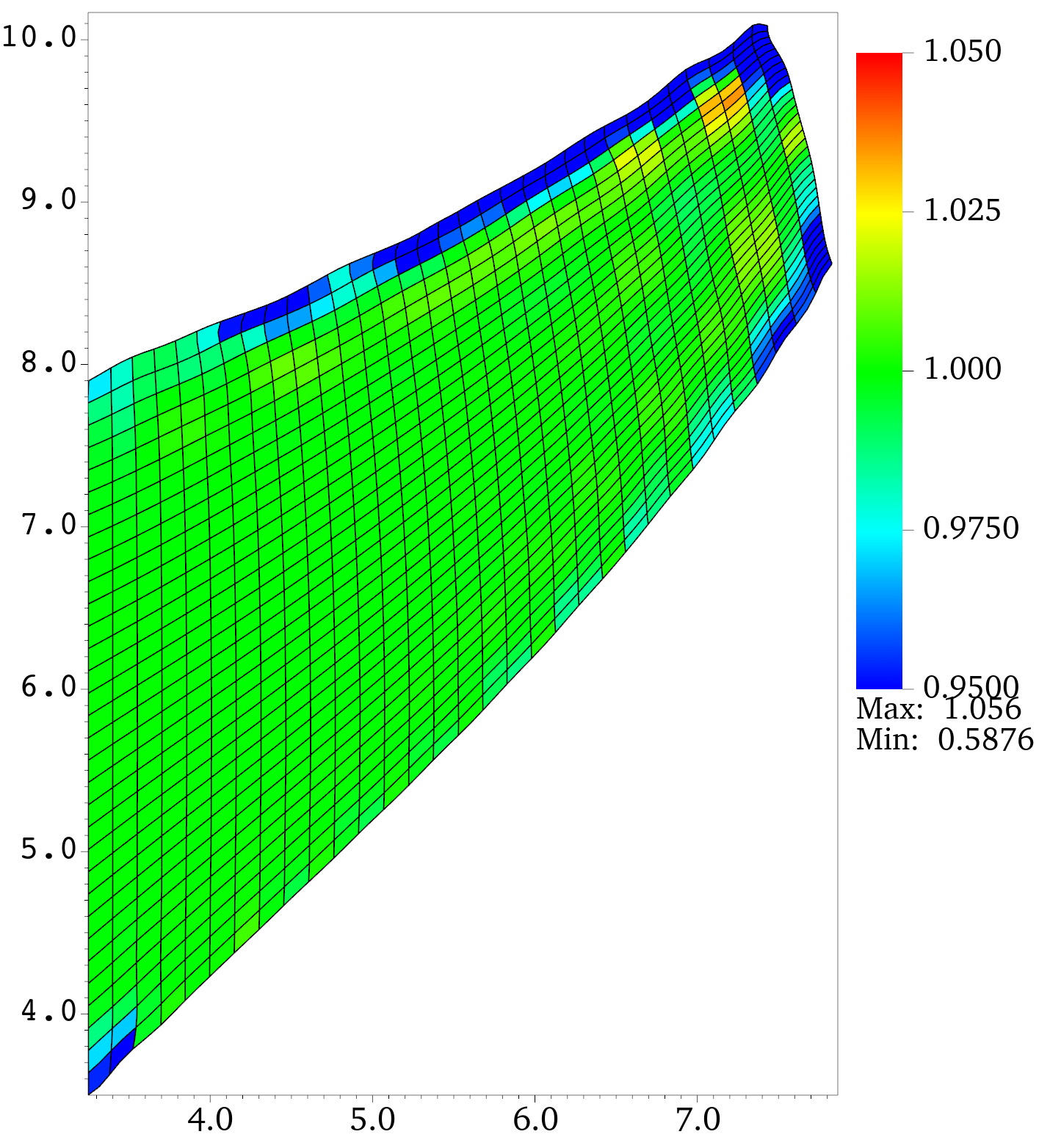}
    \caption{$\text{BS}_3$}
    \label{fig:cook_membrane_BS3_no_stab}
  \end{subfigure}
  \hfill
  \begin{subfigure}{0.49\textwidth}
    \centering
    \includegraphics[width=\textwidth]{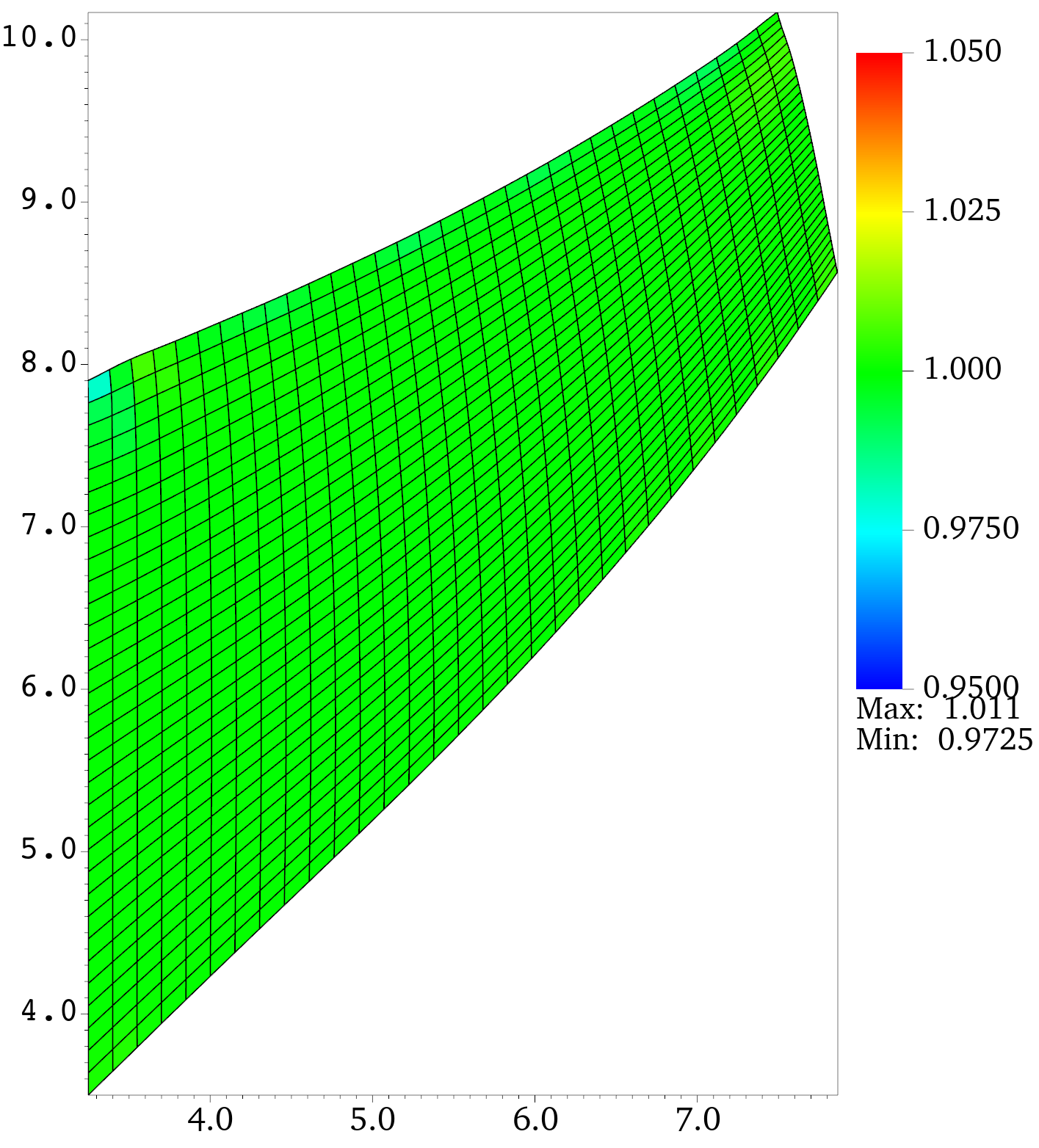}
    \caption{$\text{BS}_3$, $\nu=0.4$, modified invariants}
    \label{fig:cook_membrane_BS3_stab}
  \end{subfigure}  
\end{figure}

\begin{figure}[H]\ContinuedFloat  
  \begin{subfigure}{0.49\textwidth}
    \centering
    \includegraphics[width=\textwidth]{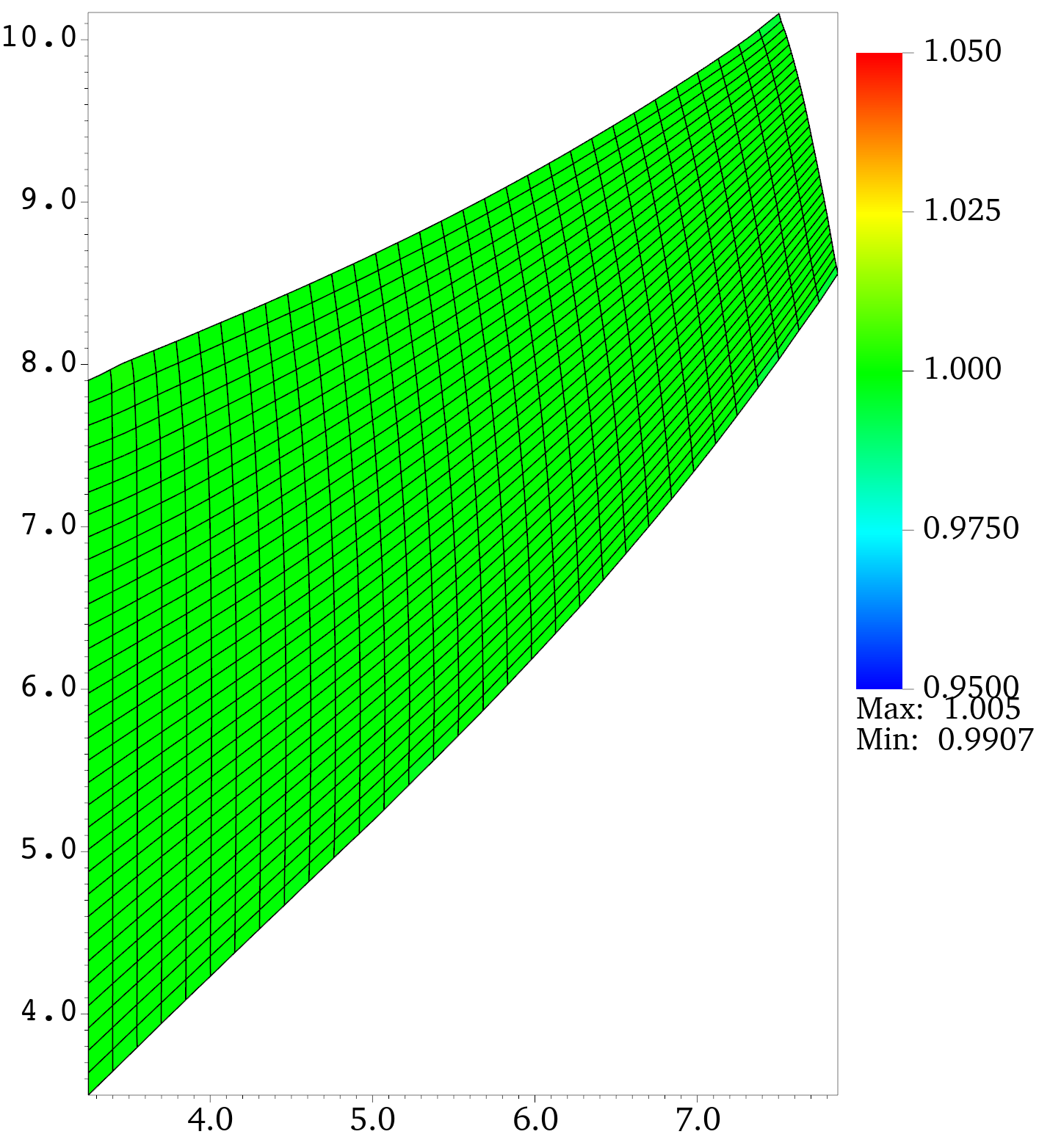}
    \caption{$\text{CBS}_{32}$}
    \label{fig:cook_membrane_CBS32_no_stab}
  \end{subfigure}
  \hfill
  \begin{subfigure}{0.49\textwidth}
    \centering
    \includegraphics[width=\textwidth]{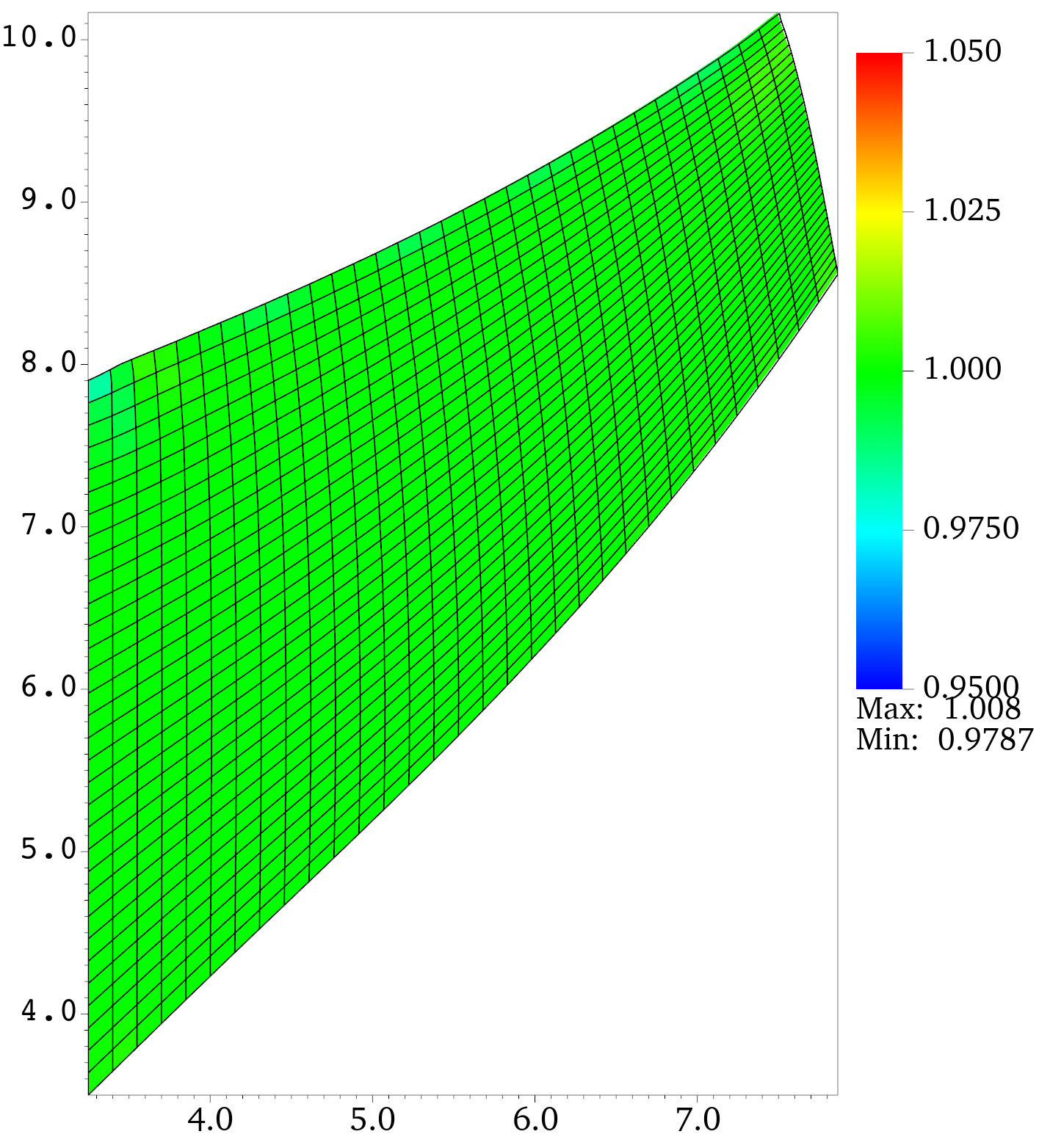}
    \caption{$\text{CBS}_{32}$, $\nu=0.4$, modified invariants}
    \label{fig:cook_membrane_CBS32_stab}
  \end{subfigure}    
\end{figure}

\begin{figure}[H]\ContinuedFloat    
  \begin{subfigure}{0.49\textwidth}
    \centering
    \includegraphics[width=\textwidth]{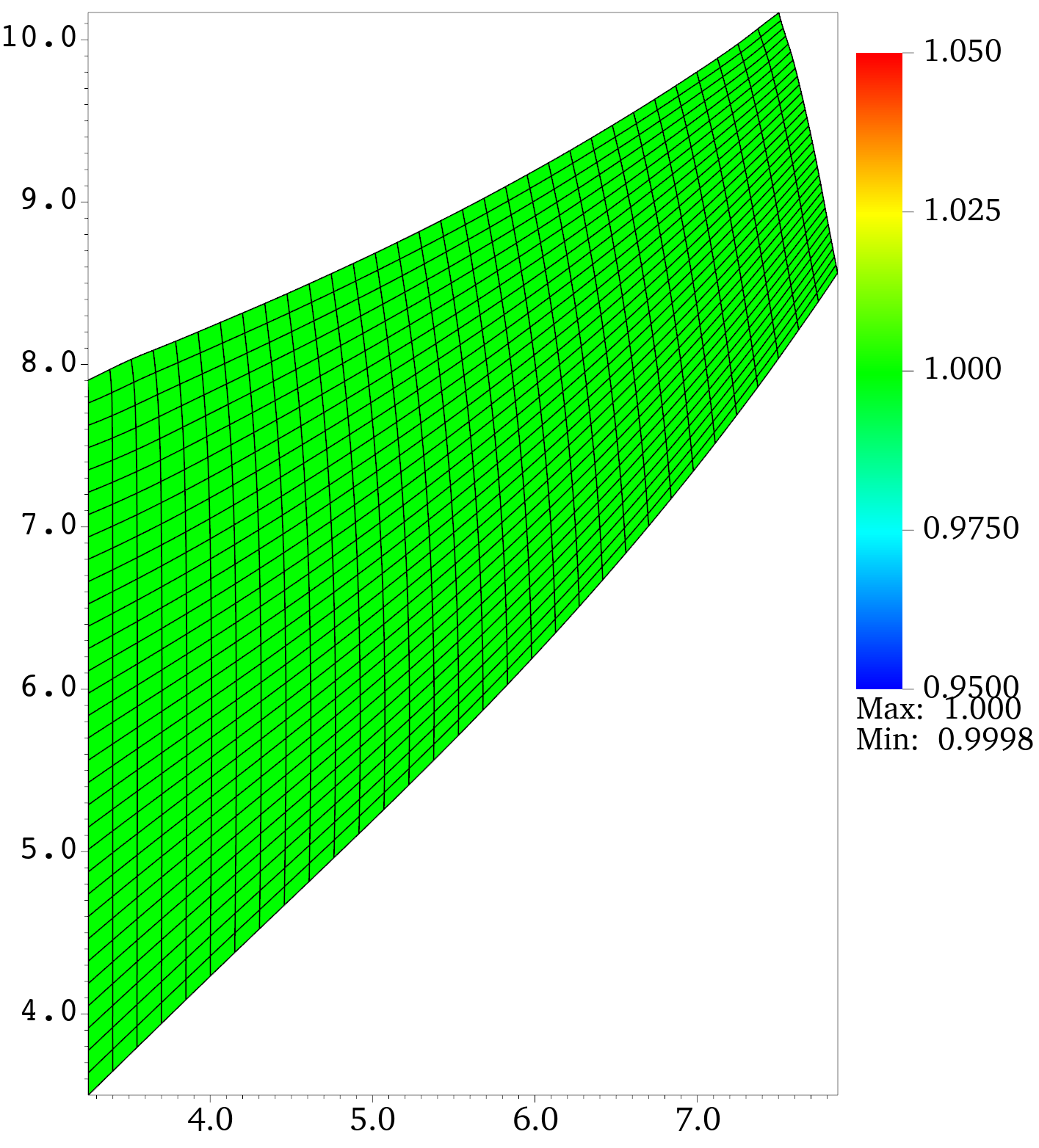}
    \caption{$\text{CBS}_{43}$}
    \label{fig:cook_membrane_CBS43_no_stab}
  \end{subfigure}
  \begin{subfigure}{0.49\textwidth}
    \centering
    \includegraphics[width=\textwidth]{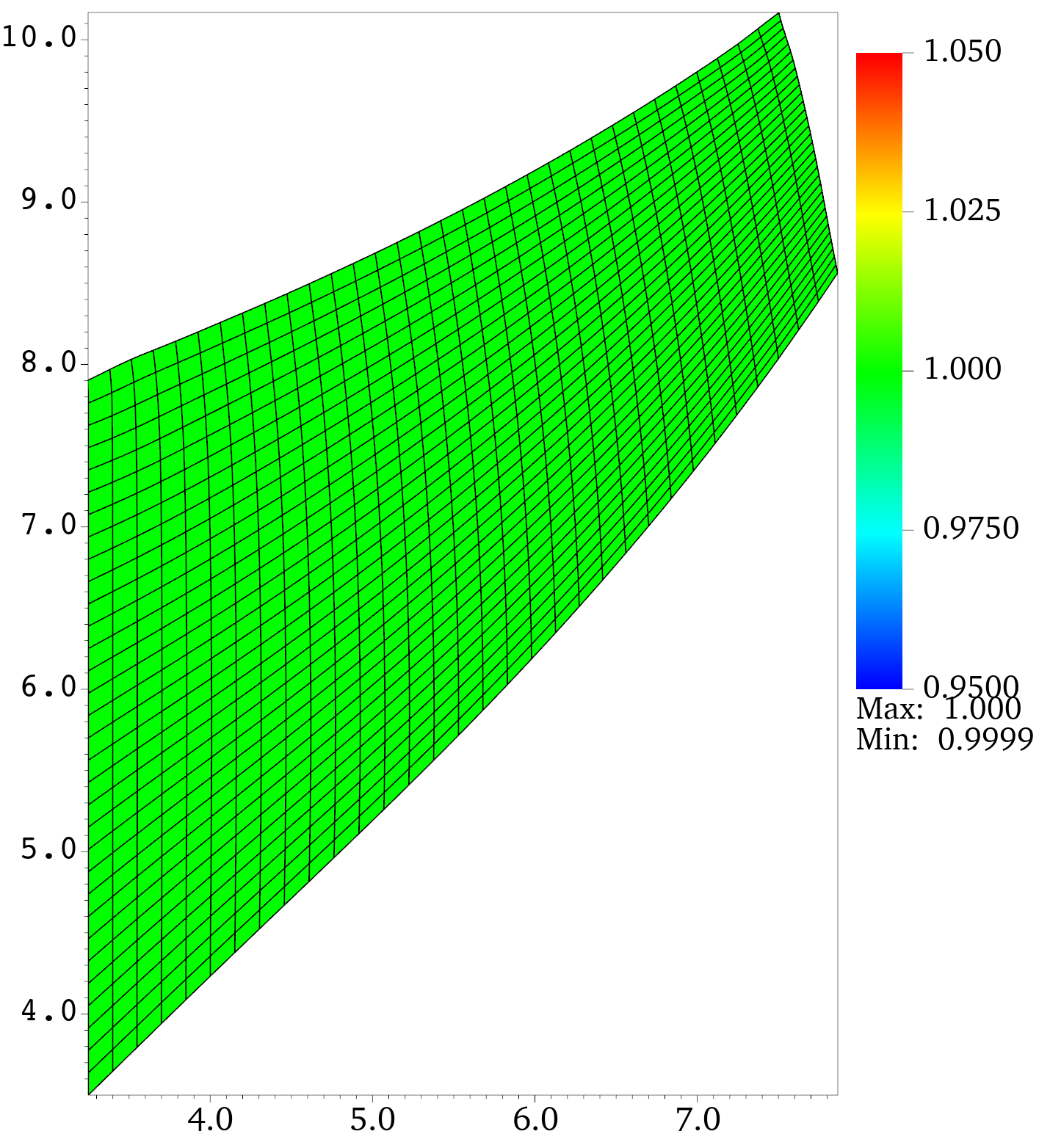}
    \caption{$\text{CBS}_{43}$, $\nu=0.4$, modified invariants}
    \label{fig:cook_membrane_CBS43_stab}
  \end{subfigure}    
  
  \caption{The comparison of the Jocobians between the unmodified invariants with no volumetric energy (left column) and with both treatments (right column) for different kernels ($\MFAC$ = 0.5, $t$ = 50 s). The special treatments improve $\text{IB}_3$ and $\text{BS}_3$ a lot for volume conservation but have little effect on $\text{CBS}_{32}$ and $\text{CBS}_{43}$.}
  \label{fig:cook_membrane_jacobian_stab_no_stab}
\end{figure}

\begin{figure}[H]
  \centering
  \begin{subfigure}{0.49\textwidth}
    \centering
    \includegraphics[width=\textwidth]{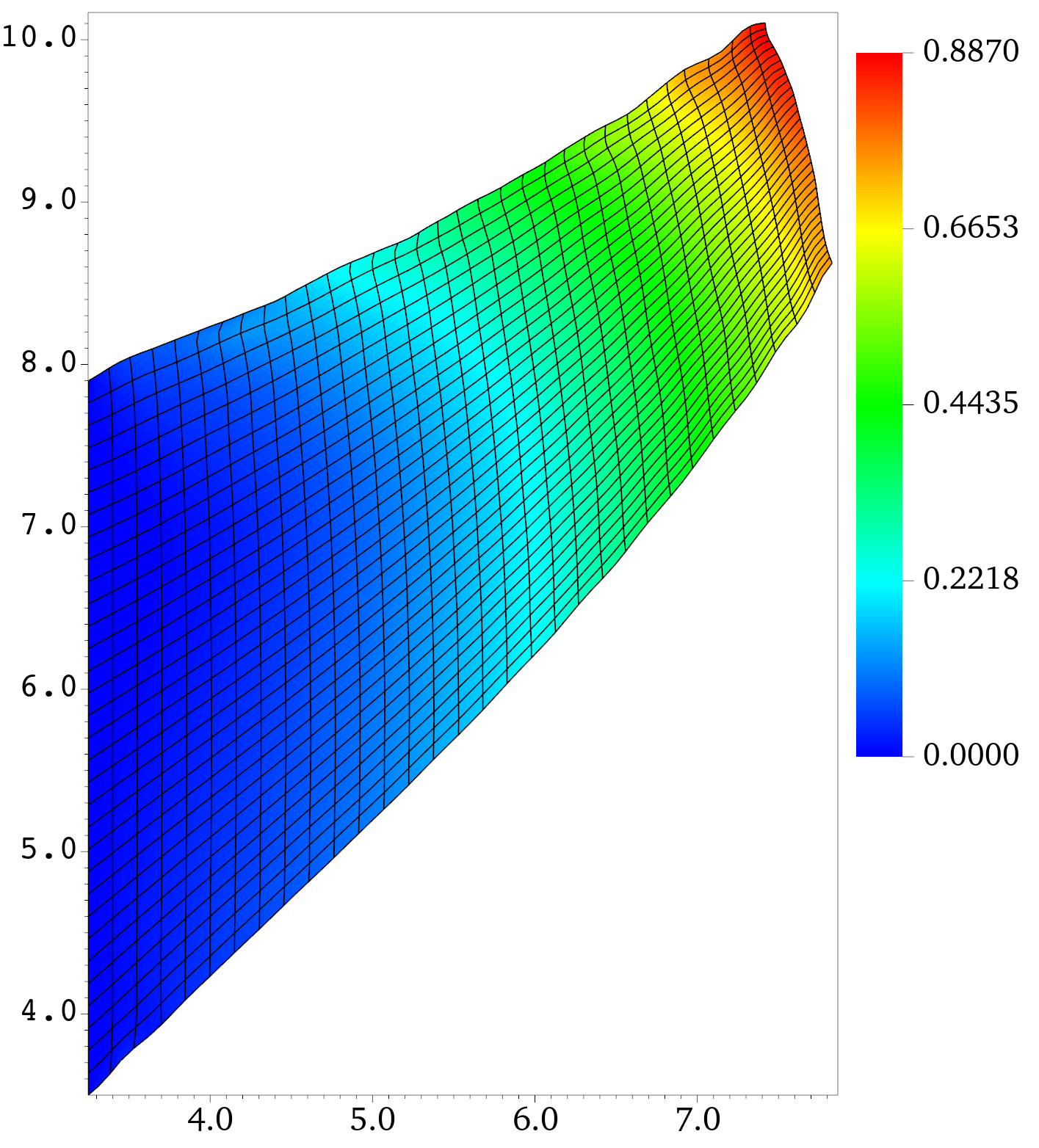}
    \caption{$\text{IB}_3$}
    \label{fig:cook_membrane_disp_IB3_no_stab}
  \end{subfigure}
  \hfill
  \begin{subfigure}{0.49\textwidth}
    \centering
    \includegraphics[width=\textwidth]{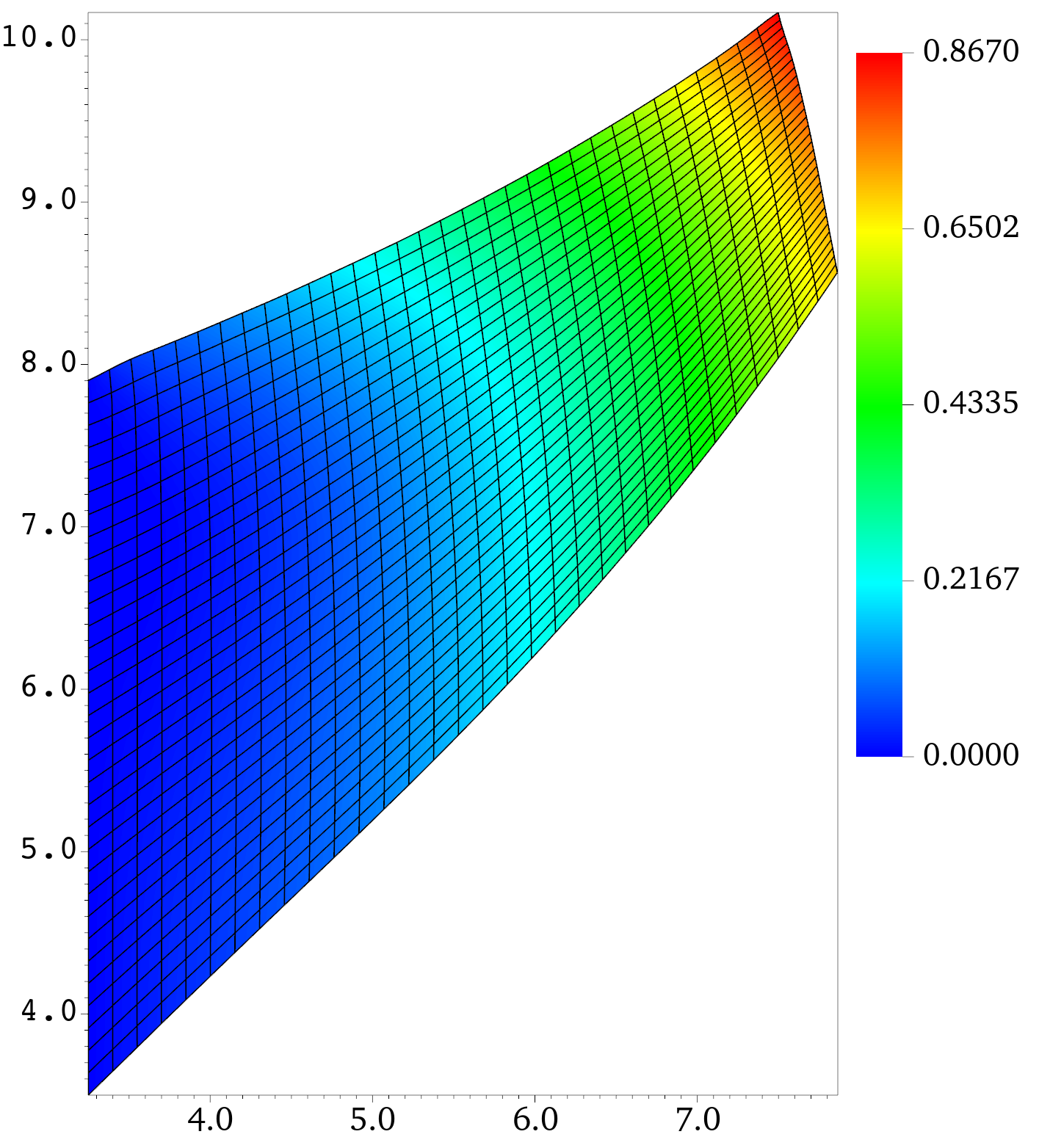}
    \caption{$\text{IB}_3$, $\nu=0.4$, modified invariants}
    \label{fig:cook_membrane_disp_IB3_stab}
  \end{subfigure}
\end{figure}

\begin{figure}[H]\ContinuedFloat    
  \begin{subfigure}{0.49\textwidth}
    \centering
    \includegraphics[width=\textwidth]{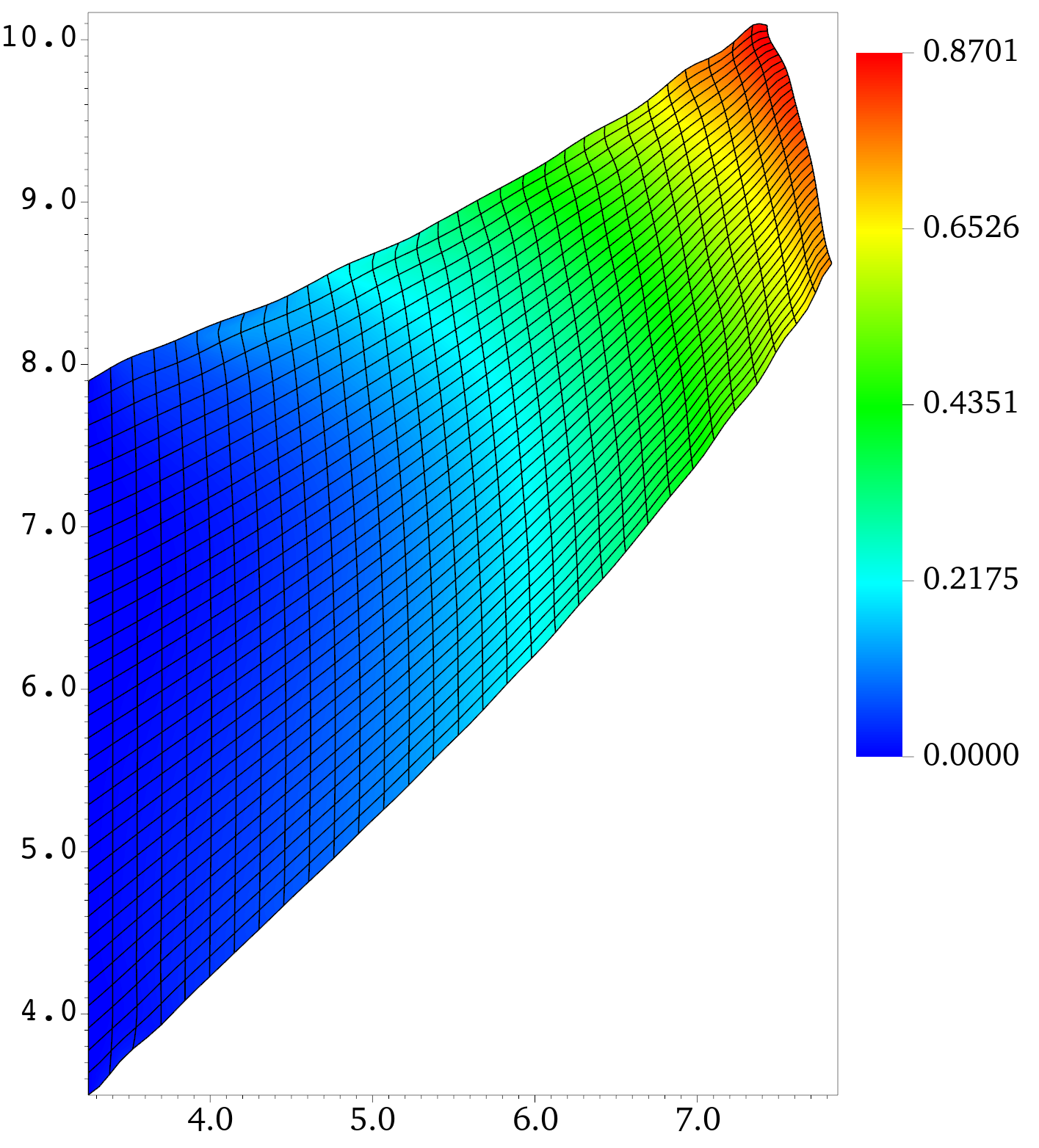}
    \caption{$\text{BS}_3$}
    \label{fig:cook_membrane_disp_BS3_no_stab}
  \end{subfigure}
  \hfill
  \begin{subfigure}{0.49\textwidth}
    \centering
    \includegraphics[width=\textwidth]{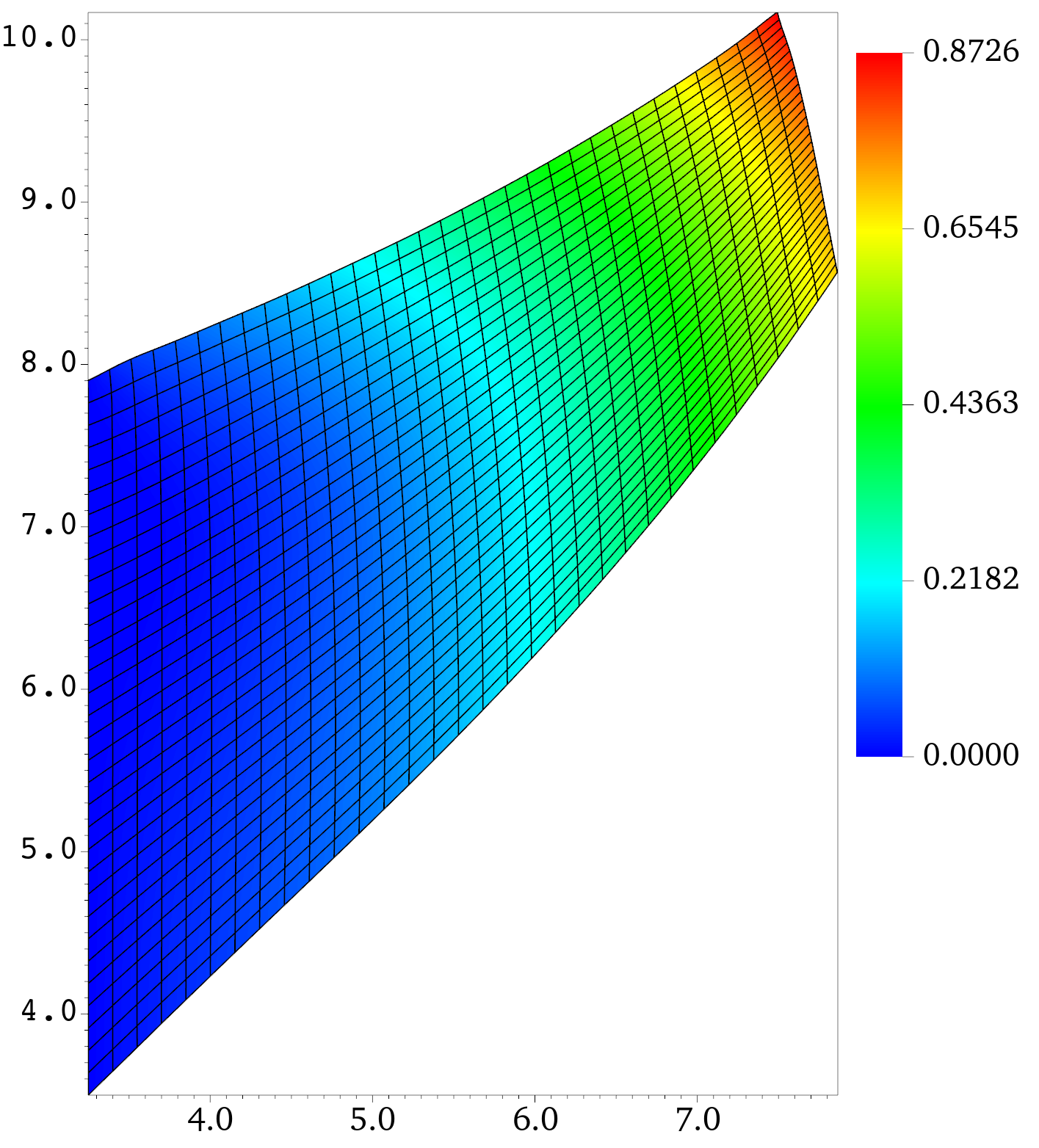}
    \caption{$\text{BS}_3$, $\nu=0.4$, modified invariants}
    \label{fig:cook_membrane_disp_BS3_stab}
  \end{subfigure}  
\end{figure}

\begin{figure}[H]\ContinuedFloat    
  \begin{subfigure}{0.49\textwidth}
    \centering
    \includegraphics[width=\textwidth]{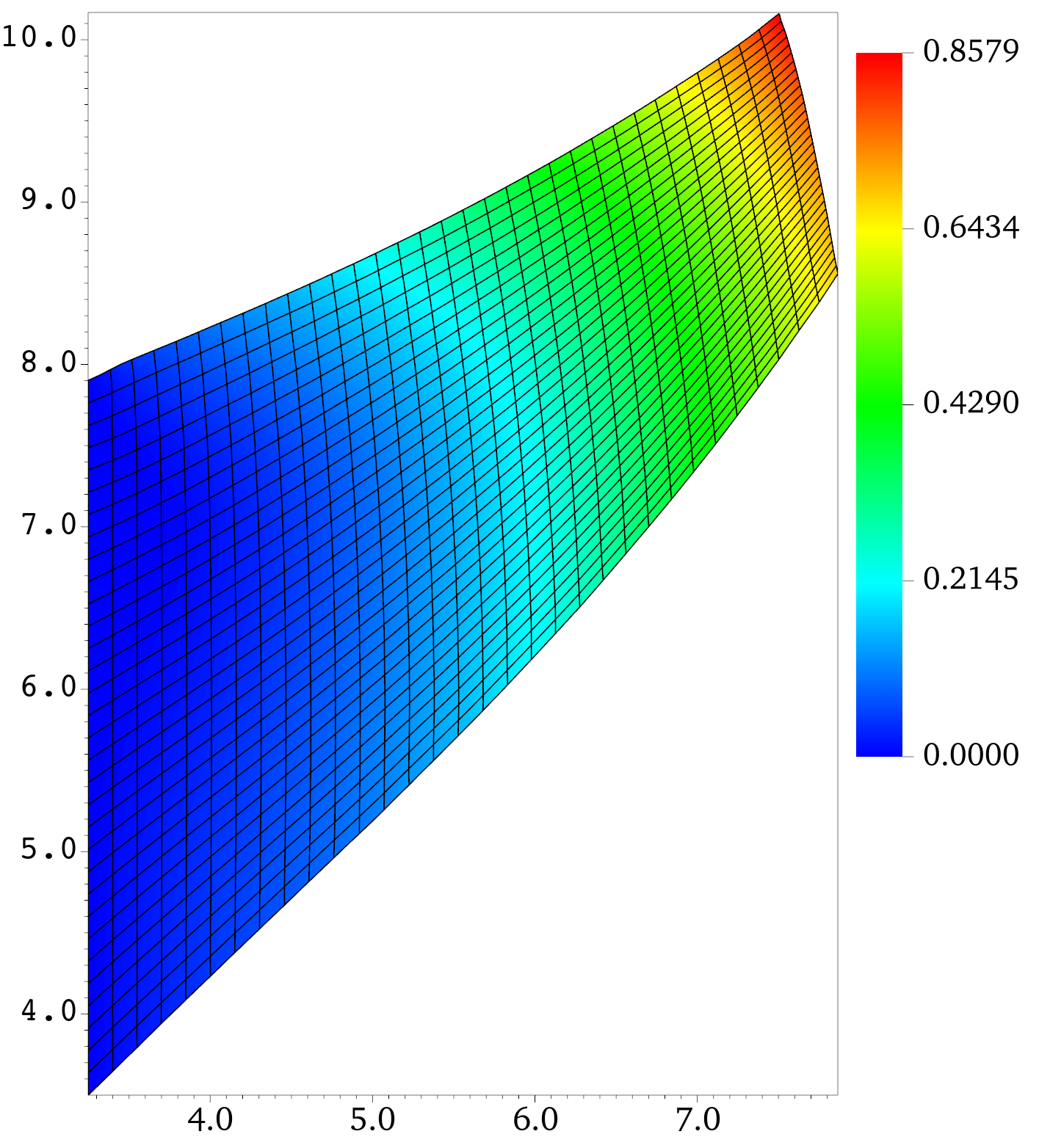}
    \caption{$\text{CBS}_{32}$}
    \label{fig:cook_membrane_disp_CBS32_no_stab}
  \end{subfigure}
  \hfill
  \begin{subfigure}{0.49\textwidth}
    \centering
    \includegraphics[width=\textwidth]{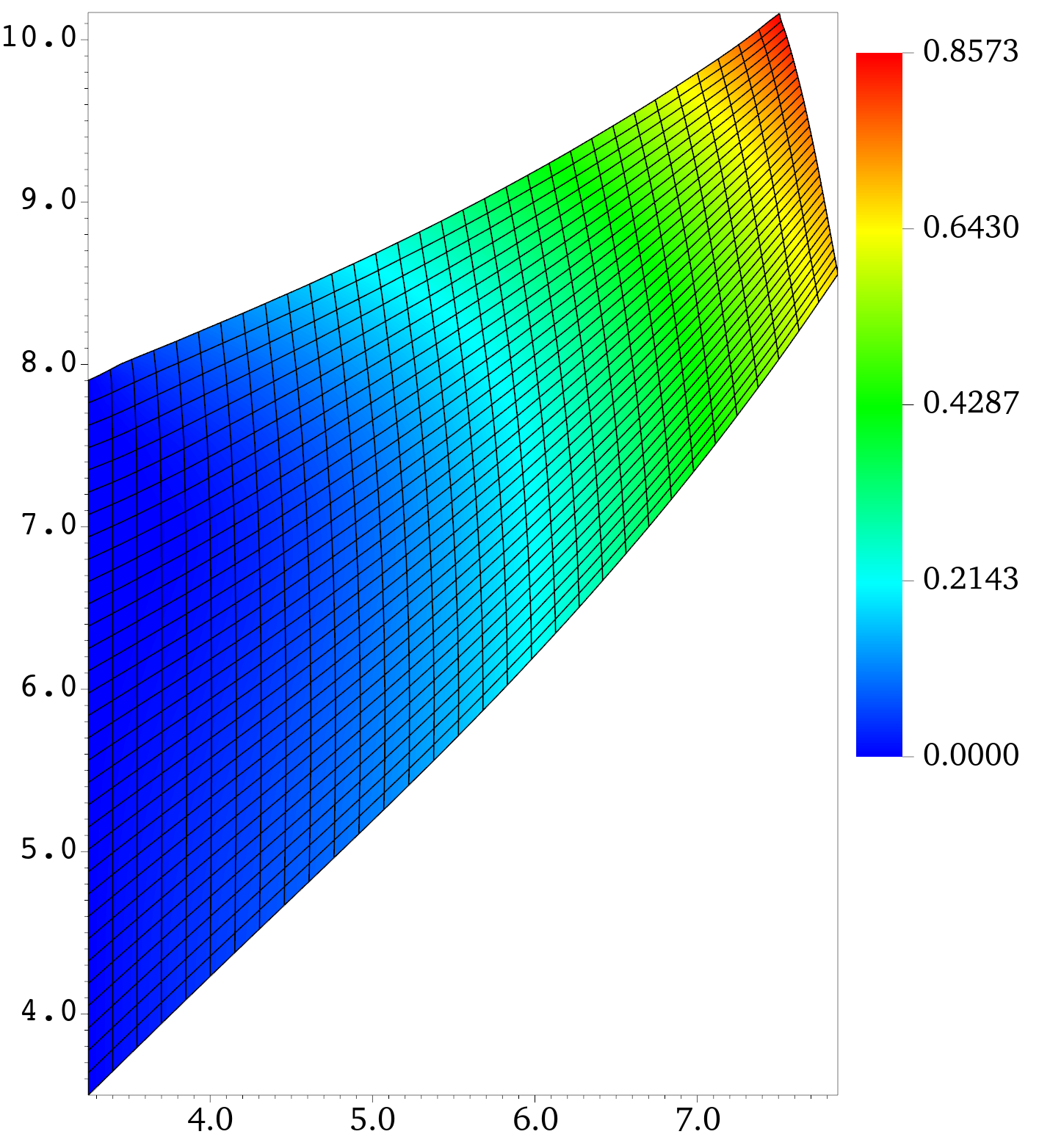}
    \caption{$\text{CBS}_{32}$, $\nu=0.4$, modified invariants}
    \label{fig:cook_membrane_disp_CBS32_stab}
  \end{subfigure}    
\end{figure}

\begin{figure}[H]\ContinuedFloat    
  \begin{subfigure}{0.49\textwidth}
    \centering
    \includegraphics[width=\textwidth]{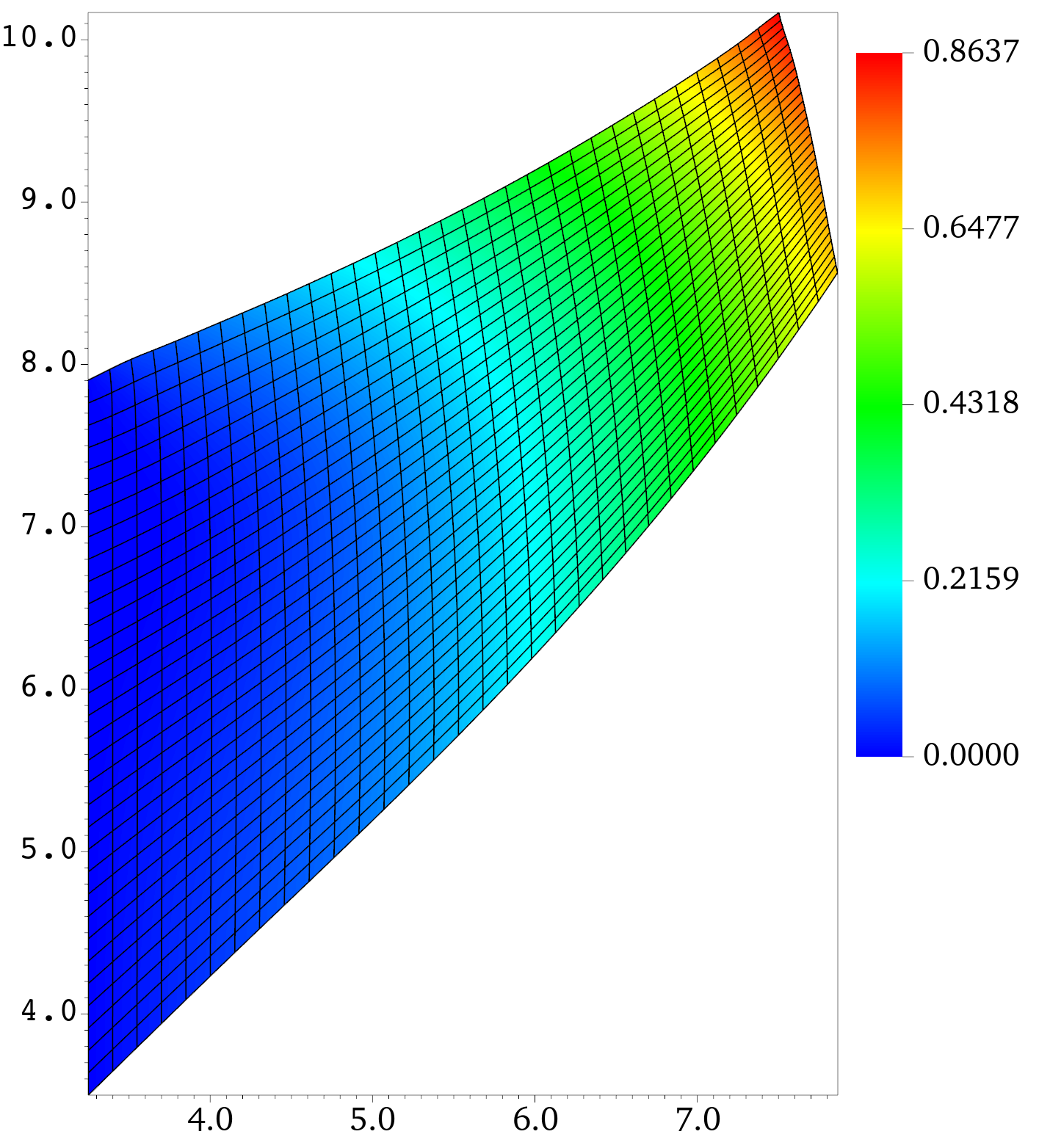}
    \caption{$\text{CBS}_{43}$}
    \label{fig:cook_membrane_disp_CBS43_no_stab}
  \end{subfigure}
  \hfill
  \begin{subfigure}{0.49\textwidth}
    \centering
    \includegraphics[width=\textwidth]{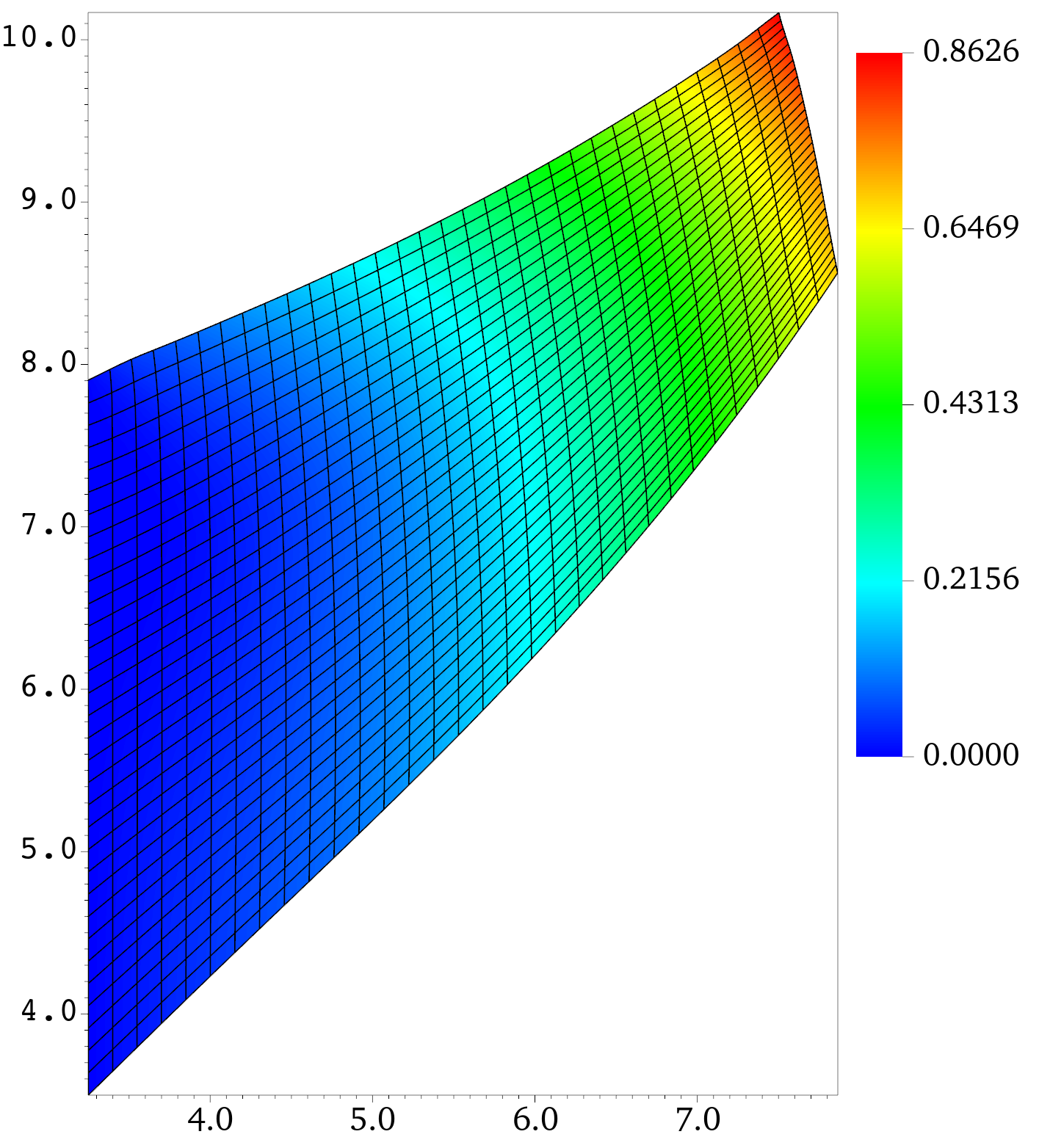}
    \caption{$\text{CBS}_{43}$, $\nu=0.4$, modified invariants}
    \label{fig:cook_membrane_disp_CBS43_stab}
  \end{subfigure}    
  
  \caption{The displacement comparison between the unmodified invariants with no volumetric energy (left column) and with both treatments (right column) for different kernels ($\MFAC$ = 0.5, $t$ = 50 s). The special treatments improve $\text{IB}_3$ and $\text{BS}_3$ a lot for volume conservation but have little effect on $\text{CBS}_{32}$ and $\text{CBS}_{43}$.}
  \label{fig:cook_membrane_disp_field_stab_no_stab}
\end{figure}

\begin{figure}[H]
  \centering
  \begin{subfigure}[b]{0.49\textwidth}
    \centering
    \includegraphics[width=\textwidth]{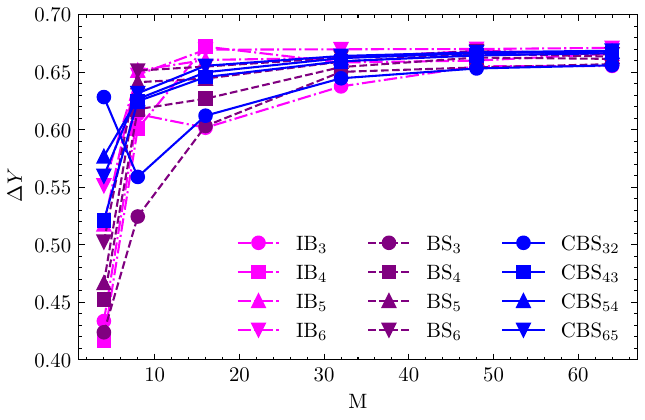}
    \caption{$\MFAC=1.5$}
    \label{fig:cook_membrane_mfac=1.5}
  \end{subfigure}
  \hfill
  \begin{subfigure}[b]{0.49\textwidth}
    \centering
    \includegraphics[width=\textwidth]{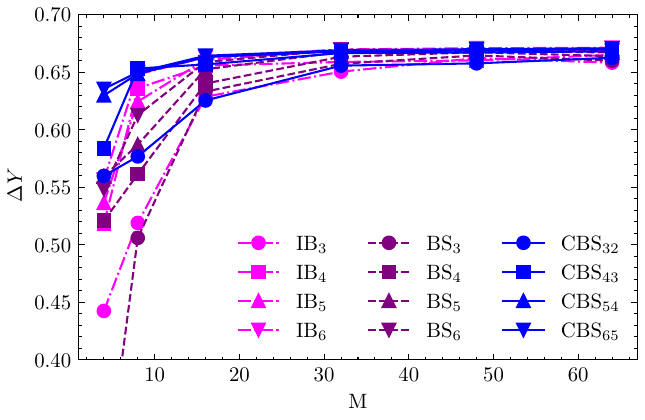}
    \caption{$\MFAC=1.25$}
    \label{fig:cook_membrane_mfac=1.25}
  \end{subfigure}
  
  \begin{subfigure}[b]{0.49\textwidth}
    \centering
    \includegraphics[width=\textwidth]{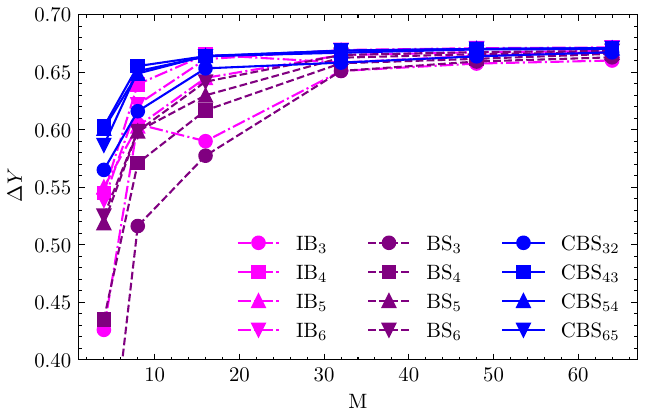}
    \caption{$\MFAC=1.0$}
    \label{fig:cook_membrane_mfac=1.0}
  \end{subfigure}
  \hfill
  \begin{subfigure}[b]{0.49\textwidth}
    \centering
    \includegraphics[width=\textwidth]{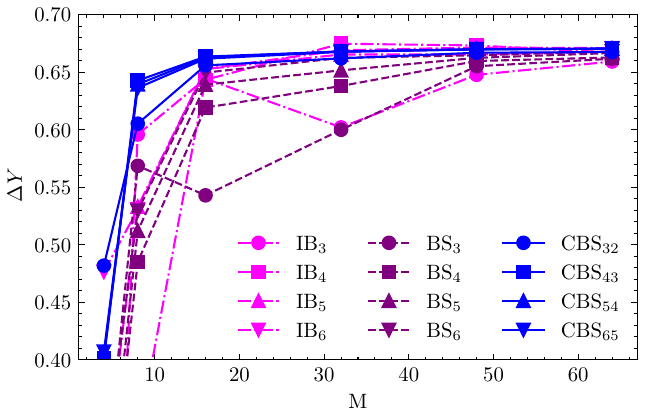}
    \caption{$\MFAC=0.75$}
    \label{fig:cook_membrane_mfac=0.75}
  \end{subfigure}  
  
  \begin{subfigure}[b]{0.49\textwidth}
    \centering
    \includegraphics[width=\textwidth]{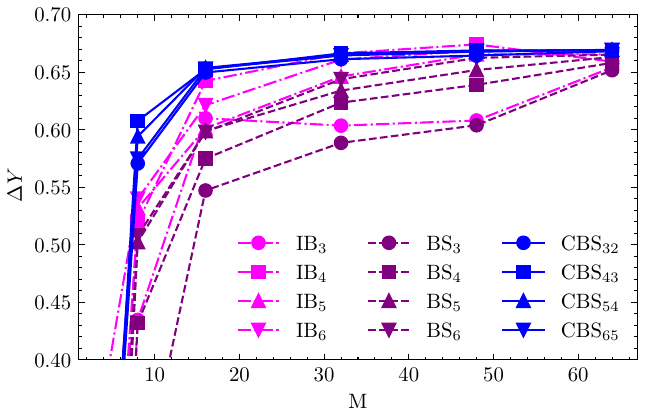}
    \caption{$\MFAC=0.5$}
    \label{fig:cook_membrane_mfac=0.5}
  \end{subfigure}
  \hfill
  
  \caption{All kernels achieve grid convergence as the grid is refined across different MFAC values. CBS kernels demonstrate consistent performance and are less sensitive to MFAC, achieving convergence on coarser grids with better results at smaller MFAC values. In contrast, IB and BS kernels perform better with larger MFAC values and require finer meshes at smaller MFAC values to achieve results comparable to CBS. Among IB and BS kernels, those with wider supports converge more quickly than those with narrower supports.}
  \label{fig:jacobian_stab_no_stab}
\end{figure}

Fig. \ref{fig:jacobian_stab_no_stab} shows the grid convergence with different MFACs in terms of the vertical displacement of the top corner. All kernels achieve grid convergence as the grid is refined across different MFAC values. All CBS kernels show consistent performance, while for IB and BS kernels, those with wider supports converge faster than those with narrower supports. CBS kernels are less sensitive to MFAC and converge at a coarse grid level, with smaller MFAC values providing better convergence. In contrast, IB and BS kernels perform better with larger MFAC values, and for small MFAC values, they require finer meshes to reach results comparable to those of CBS. 

\begin{figure}[H]
\centering
\includegraphics[width=0.6\textwidth]{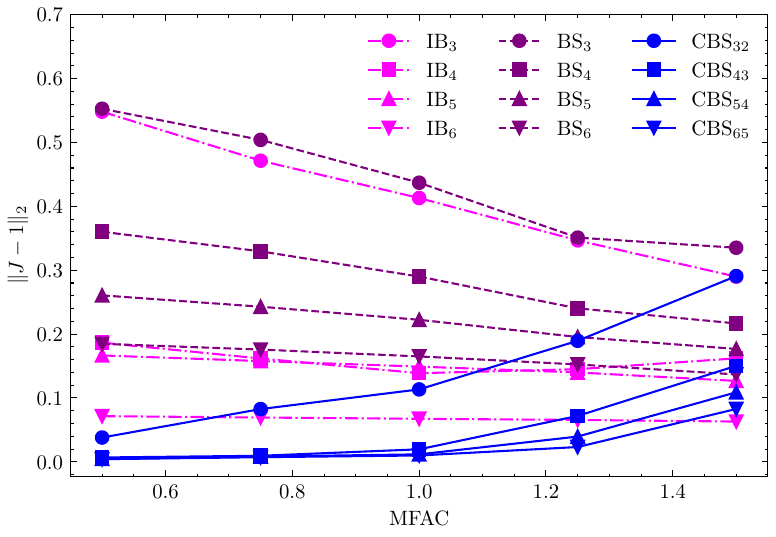}
\caption{Error norms of the Jacobian as a function of MFAC with $N=90$. IB and BS kernels show improved results with larger MFAC values, while CBS kernels display the opposite trend, with smaller MFAC values yielding lower errors. Additionally, different CBS kernels exhibit less discrepancy in error compared to IB and BS kernels.}
\label{fig:cook_membrane_mfac_J}
\end{figure}

Fig. \ref{fig:cook_membrane_mfac_J} shows the error norms of the Jacobian as a function of MFAC with $N=90$. IB and BS kernels show improved results with larger MFAC values, while CBS kernels display the opposite trend, with smaller MFAC values yielding lower errors. Additionally, different CBS kernels exhibit less discrepancy in error compared to IB and BS kernels.

\subsection{Shear-dominated Flows}
\subsubsection{Slanted Channel Flow}
\label{subsec:channel}

\begin{figure}[H]
 \centering
 \includegraphics[scale = 0.5]{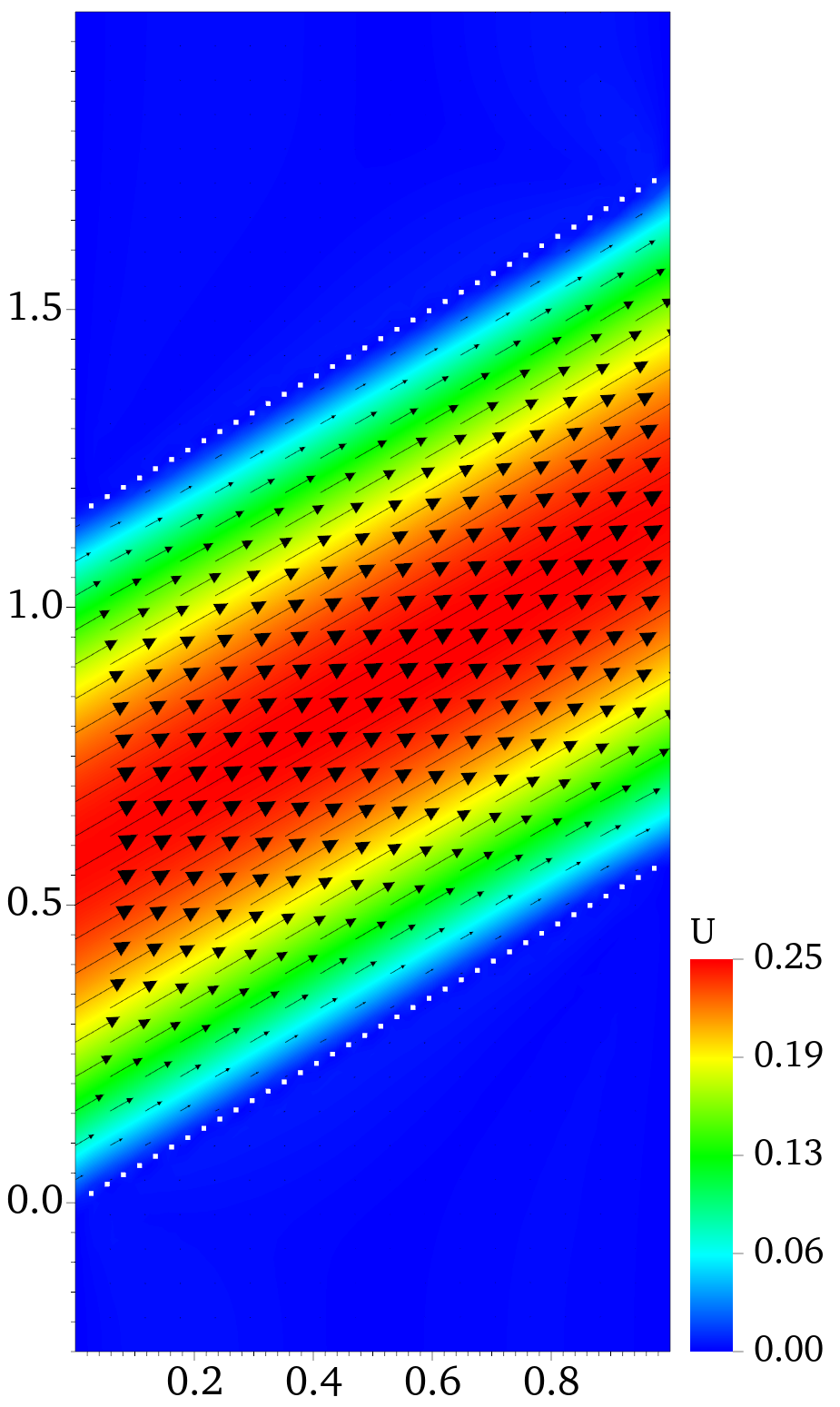}
 \caption{Computational setup of the slanted channel flow problem (inclination angle $\pi/6$). The color map shows the velocity field computed using $\text{CBS}_{32}$ kernel. The white vertical line at $x = 0.5$ indicates the location of velocity profile measurements, and the dots represent Lagrangian markers defining the top and bottom channel plates.}
 \label{fig:slanted_channel} 
 \end{figure}

To evaluate the performance of different kernels in non-grid-aligned configurations, we examine flow through a slanted channel, following the two-dimensional benchmark presented in Gruninger et al. \cite{gruninger2024benchmarking}. The computational domain $\Omega = [0, 1]\times[-0.25, 2]$ contains two parallel plates separated by channel width $D$. Figure~\ref{fig:slanted_channel} illustrates the channel geometry, inclined at angle $\theta=\pi/6$, and shows a representative velocity field computed using $\text{CBS}_{32}$.

 The channel is inclined rather than vertical or horizontal to avoid grid-aligned discretization of the channel walls. The exact steady-state solution for flow through the inclined channel can be derived through the coordinate transformation of the plane Poiseuille equation. As in previous work \cite{gruninger2024benchmarking}, the analytic solution is given by:
\begin{equation}
\label{eq:plane_poiseuille}
\begin{aligned}
   u(x,y) &= -\frac{\Delta p \cos\theta}{2\mu L}\left(y\cos\theta - x\sin\theta\right)\left(y\cos\theta - x\sin\theta - h\right),\\ 
   v(x,y) &= -\frac{\Delta p \sin\theta}{2\mu L}\left(y\cos\theta - x\sin\theta\right)\left(y\cos\theta - x\sin\theta - h\right),
\end{aligned}    
\end{equation}
where $h$ denotes the height of the channel in its horizontal configuration and $\theta$ the counterclockwise rotation angle about the origin.

The simulation uses the following parameters: the horizontal channel width $h = 1$, and the slanted channel width $D=h/\cos(\theta )$, dynamic viscosity $\mu = 0.5$, density $\rho = 1.0$, pressure gradient $\Delta p/L = 1.0$, and maximum velocity $U_\text{max} = 0.25$. The fine-grid Cartesian cell size is set to $\euleriandx = L/N$, where $N = 32$ with a time step size of $\Delta t = 0.15 \euleriandx$.

We impose velocity boundary conditions at the inlet and outlet using the analytical steady-state solution from Eq.~\eqref{eq:plane_poiseuille}. To maintain the rigidity and stationary position of the structure, we employ both the penalty stiffness and the penalty body and damping forces described in Eq.~\eqref{eq:tether_force}. The penalty parameters are empirically determined as the largest values that maintain numerical stability for each combination of kernel, grid spacing, and time step size.

\begin{figure}[H]
 \centering
 \includegraphics[width=0.49\textwidth]{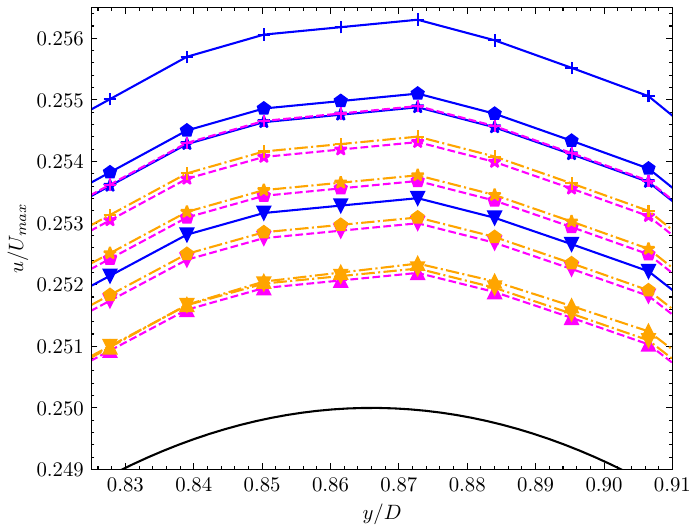}
 \includegraphics[width=0.49\textwidth]{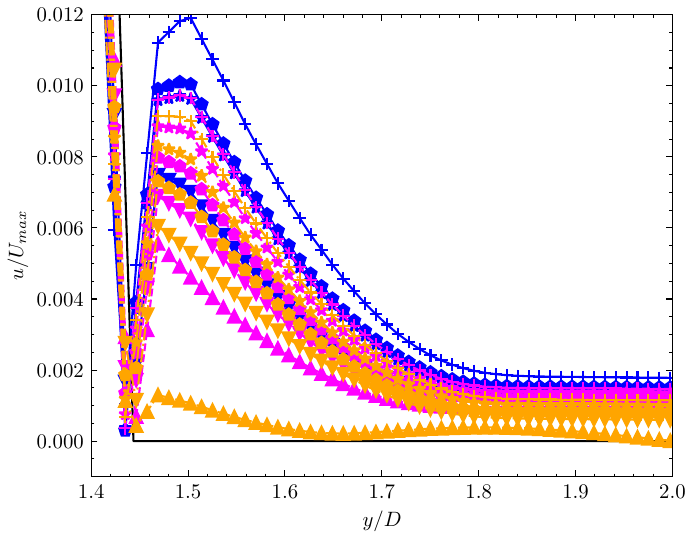}
\caption{Comparison of velocity profiles at $x=0.5$ for different kernel types. Kernels with narrower support regions show better accuracy in capturing peak velocities, with BS and CBS kernels of equal support performing similarly and outperforming their IB counterparts. The width of the numerical boundary layer decreases with decreasing kernel support size, with $\text{CBS}_{21}$ demonstrating the thinnest numerical boundary layer among all tested kernels}
 \label{fig:slanted_channel_velocity_profile} 
 \end{figure}

This benchmark evaluates explicitly how different kernel choices affect the accuracy of flow computations within a confined, stationary geometry.
Fig. \ref{fig:slanted_channel_velocity_profile} compares the velocity profiles at $x=0.5$ for different kernel types. Kernels with narrower support regions show better accuracy in capturing peak velocities, with BS and CBS kernels of equal support performing similarly and outperforming their IB counterparts. The width of the numerical boundary layer decreases with decreasing kernel support size, with $\text{CBS}_{21}$ demonstrating the thinnest numerical boundary layer among all tested kernels.

\begin{figure}[H]
 \centering
 \includegraphics[width=0.5\textwidth]{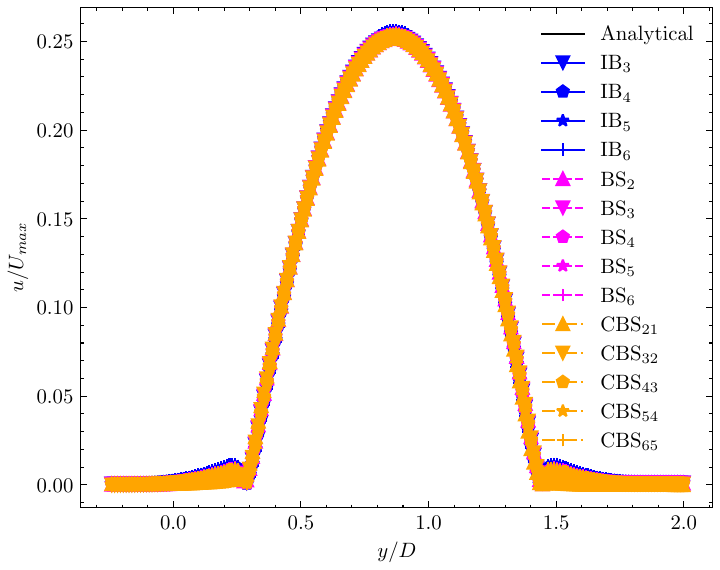}
\end{figure}
\subsubsection{Modified Turek-Hron}

\begin{figure}[H]
\centering
\includegraphics[width=0.8\linewidth]{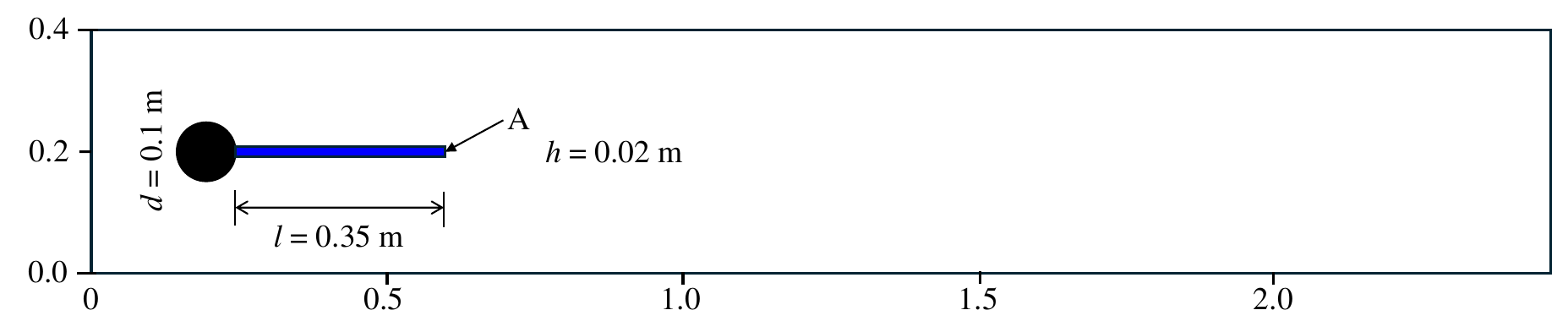}
\caption{Schematic of the Turek-Hron benchmark}
\label{fig:turek_schematic}
\end{figure}

We investigate a modified version of the Turek-Hron fluid-structure interaction (FSI) benchmark \cite{turek2006proposal}, which simulates flow around a flexible elastic beam attached to a fixed circular cylinder \cite{lee2022lagrangian}. While the original benchmark specifies domain dimensions of $L = 2.5$ and $H = 0.41$, we extend the length to $L = 2.46 = 6.0H$ to accommodate square Cartesian grid cells. This modification has a negligible impact on the benchmark results.
The computational setup uses a fine-grid Cartesian cell size of $\Delta x = L/N$ with a time step of $\Delta t = 0.00164 \Delta x$, where $N$ is the grid number along the longest dimension of the fluid domain. The structure consists of a circular cylinder (diameter d = 0.1) centered at (0.2, 0.2); (2) and an elastic beam (length $l = 0.35$, height $h = 0.02$) fixed to the cylinder's rear. A control point $A$ (initial position is (0.6, 0.2)) is used for monitoring displacement. Fig. \ref{fig:turek_schematic} shows the setup schematic.
The boundary conditions are specified as follows: at the inlet ($x = 0$), $u(0, y) = 1.5U y(H - y)/(H/2)^2$, where $U = 2$ is the average velocity; at the outlet ($x = L$), zero normal traction and zero tangential velocity are imposed; and along the top and bottom walls ($y = 0, H$), zero velocity conditions are enforced.

The flow parameters yield $Re = \rho Ud/\mu$ = 200, with $\rho$ = 1000 and $\mu$ = 1.

The structure is modeled using the Saint Venant-Kirchhoff constitutive law:
$\mathbb{S} = \lambda_s \text{tr}(\mathbb{E}) \mathbb{I} + 2\mu_s \mathbb{E}$
where $\mathbb{S}$ is the second Piola-Kirchhoff stress tensor, $\mathbb{E}$ is the Green-Lagrange strain tensor, $\mathbb{I}$ is the second-order identity tensor, and material parameters are $\mu_s = 1\times 10^6$, and $\lambda = 8\times 10^6$. The cylinder is constrained using a spring tether force with penalty parameter $\kappa_s = 5.0\times 10^4 \Delta x / \Delta t^2$.

We examine three kernels: $ \text{IB}_3$, $ \text{BS}_{3}$, and $ \text{CBS}_{32}$. The fluid domain uses $N = 128$ grid points along its longest dimension, providing sufficient resolution to isolate the effects of solid mesh refinement. We investigate MFAC values of 0.5, 0.75, 1.0, 1.25, and 1.5. 

\begin{figure}[H]
\centering
\includegraphics[width=0.8\linewidth]{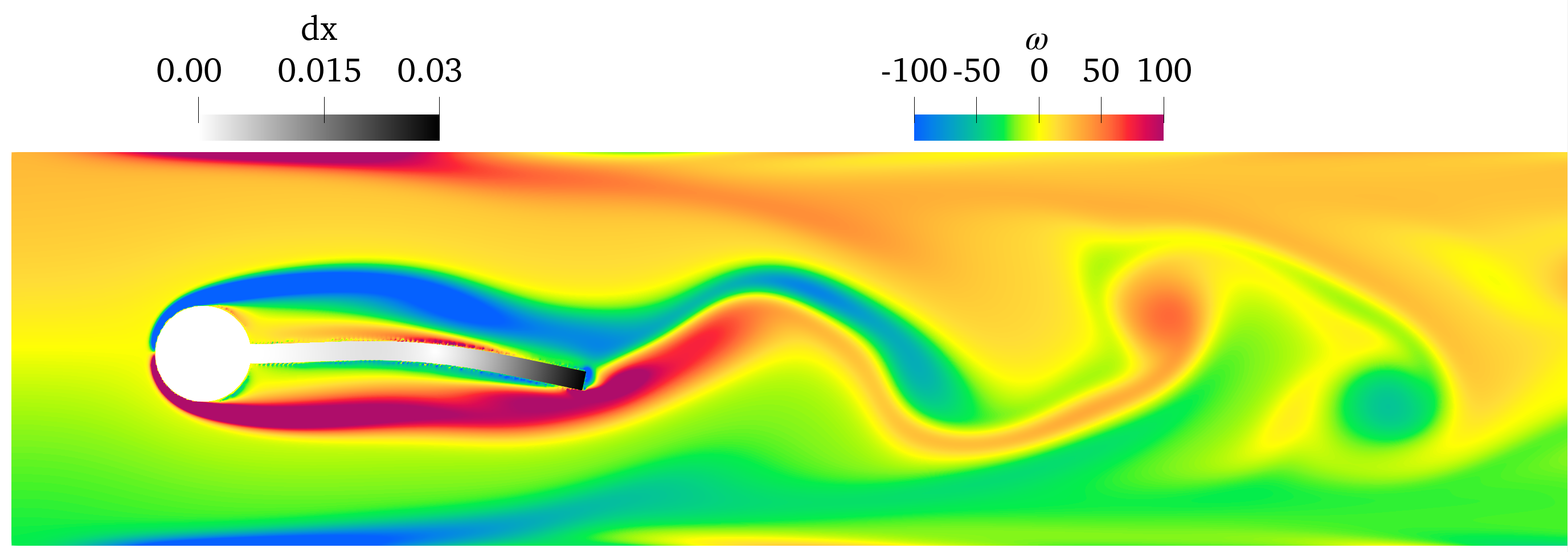}
\caption{Vorticity field of $\text{CBS}_{32}$ with MFAC=0.5. The gray colormap shows the displacement magnitude.}
\label{fig:turek_omega}
\end{figure}

Figure \ref{fig:turek_omega} presents a representative color map of the vorticity field, highlighting the deformed beam simulated with the $\text{CBS}_{32}$ kernel and MFAC = 0.5.

\begin{figure}[H]
\centering
\includegraphics[width=0.49\linewidth]{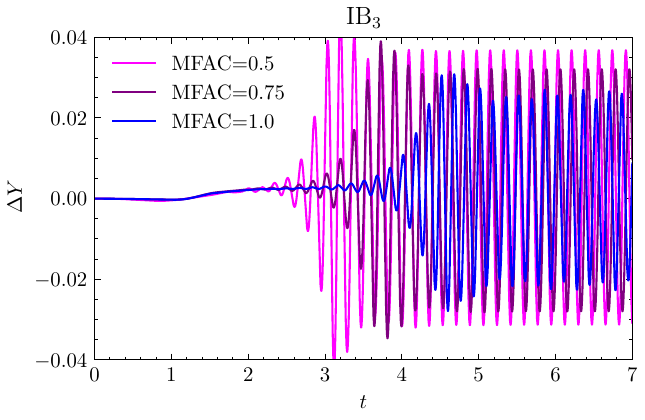}
\includegraphics[width=0.15\linewidth]{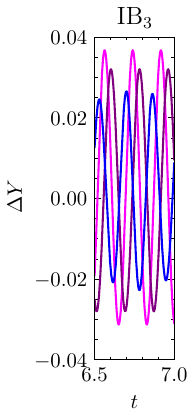}
\includegraphics[width=0.49\linewidth]{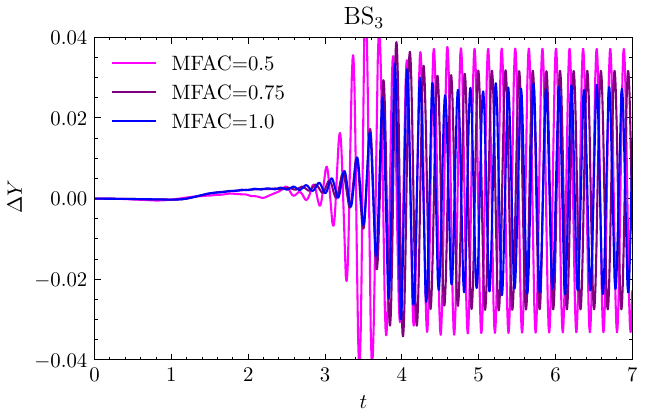}
\includegraphics[width=0.15\linewidth]{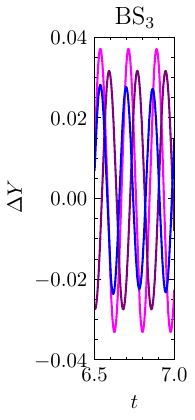}
\includegraphics[width=0.49\linewidth]{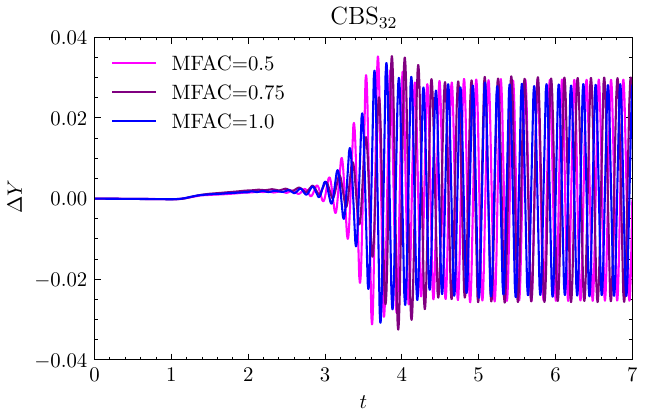}
\includegraphics[width=0.15\linewidth]{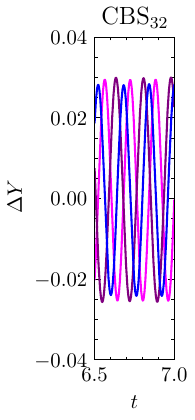}
\caption{Vertical displacement of point $A$ (as shown in Fig. \ref{fig:turek_schematic} under varying MFAC values for different kernels. Right panels show detailed oscillations during $t=6.5$--7.0.}
\label{fig:turek_dy_history}
\end{figure}

Table \ref{tab:turek_dy} summarizes the maximum vertical displacements ($\Delta Y$) at point A (shown in Fig. \ref{fig:turek_schematic}) for each kernel type across different MFAC values. The $\text{IB}_3$ kernel exhibits the most stable behavior concerning MFAC variations, while CBS kernels require smaller MFAC values and fail when MFAC exceeds 1, consistent with observations from other benchmarks. Figure \ref{fig:turek_dy_history} illustrates the oscillation histories of the vertical displacement for three MFAC values. The IB and BS kernels yield similar results, with smaller MFAC values generally predicting larger displacements. In contrast, the CBS kernel is less sensitive to the solid mesh resolution. The observed phase shift across different MFAC values is attributed to time step variations.

\begin{table}[H]
  \centering
    \caption{Maximum vertical displacements for different kernel types across MFAC values. Missing data points indicate timestepping instabilities encountered when using a time step size of $\Delta t = 10^{-6}$ s.}
    \begin{tabular}{lccccc}
    \toprule
    \multicolumn{1}{c}{\multirow{2}[4]{*}{Kernel}} & \multicolumn{5}{c}{MFAC} \\
\cmidrule{2-6}      & \multicolumn{1}{r}{0.5} & \multicolumn{1}{r}{0.75} & \multicolumn{1}{r}{1} & \multicolumn{1}{r}{1.25} & \multicolumn{1}{r}{1.5} \\
    \midrule
    $\text{IB}_3$ & 0.03686 & 0.03215 & 0.02794 & 0.02633 & 0.03087 \\
    $BS_3$ & 0.03719 & 0.03182 & 0.02860 & 0.02723 & \textbackslash{} \\
    $CBS_{32}$ & 0.02964 & 0.03030 & 0.02879 & \textbackslash{} & \textbackslash{} \\
    \bottomrule
    \end{tabular}%
  \label{tab:turek_dy}%
\end{table}%

\subsection{Bioprosthetic Heart Valve Dynamics}
\label{BHV}

Previous studies by Lee et al.\cite{lee2022lagrangian} and Wells et al.\cite{wells2023nodal} have examined bioprosthetic heart valve (BHV) dynamics in a pulse-duplicator using isotropic kernels. We adopt their setup to evaluate CBS kernel performance. This fluid-structure interaction model presents a complex combination of substantial pressure-loads (when the valve is closed) and shear-dominant flows (when the valve is open).

The bovine pericardial BHV leaflets are modeled using a modified version~\cite{Lee2020} of the Holzapfel--Gasser--Ogden model~\cite{Gasser2006}:
\begin{equation}
W_{\text{BHV}} = C_{10}{\exp{\left[C_{01}(\bar{I}_1-3)\right]}-1} + \frac{k_1}{2k_2}{\exp{\left[k_2(\kappa\bar{I}_1 + (1-3\kappa)\bar{I}_{4}^{\star}-1)^2\right]-1}},
\end{equation}
where $\bar{I}_{4}^{\star} = \max(\bar{I}_{4}, 1) = \max(\e_{0}^{\text{T}} \bar{\CC}\e_{0}, 1)$, and $\e_{0}$ represents a unit vector aligned with the mean fiber direction in the reference configuration. The parameter $\kappa \in [0,\frac{1}{3}]$ characterizes collagen fiber angle dispersion. Following Lee et al.~\cite{Lee2020}, we use material parameters $C_{10} = 0.119$~kPa, $C_{01} = 22.59$, $k_1 = 2.38$~MPa, $k_2 = 149.8$, and $\kappa = 0.292$.

The fluid properties are set to $\rho = 1.0$ g/cm$^3$ and $\mu = 1.0$ cP. The peak Reynolds number~\cite{Lee2020} is calculated as $Re_\text{peak}=\frac{\rho Q_\text{peak}D}{\mu A}\approx 14,800$, with $D=25$~mm and $A$ indicating the valve's geometrical diameter and cross-sectional area, respectively.

The computational domain spans 5.05 cm $\times$ 10.1 cm $\times$ 5.05 cm. We use a three-level locally refined grid with a refinement ratio of two between levels and an $N/2 \times N \times N/2$ coarse grid ($N = 64$), yielding a fine-grid Cartesian resolution of 0.4~mm. The time step size is $\Delta t = 2.5 \times 10^{-6} \ \text{s}$ and is systematically reduced as needed to maintain stability throughout the simulation. First-order tetrahedral elements are used for the test section, while second-order tetrahedral elements are used for the valve leaflets, with MFAC being approximately 1 and 0.5 respectively.  The piecewise-linear kernel \cite{lee2022lagrangian} is used for the test section, and $\text{BS}_3$, $\text{CBS}_{32}$, and $\text{CBS}_{43}$ kernels are evaluated for the valve leaflets.  

\begin{figure}[H]
    \centering
    \begin{subfigure} {0.07\textwidth}
        \centering
            \includegraphics[width=\linewidth]{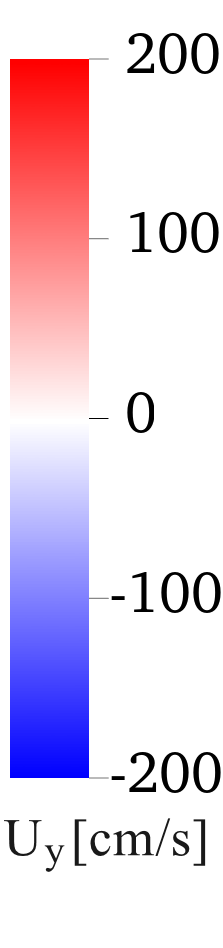}
    \end{subfigure}
    \hfill
    \begin{subfigure}{0.22\textwidth}
    \includegraphics[width=\linewidth]{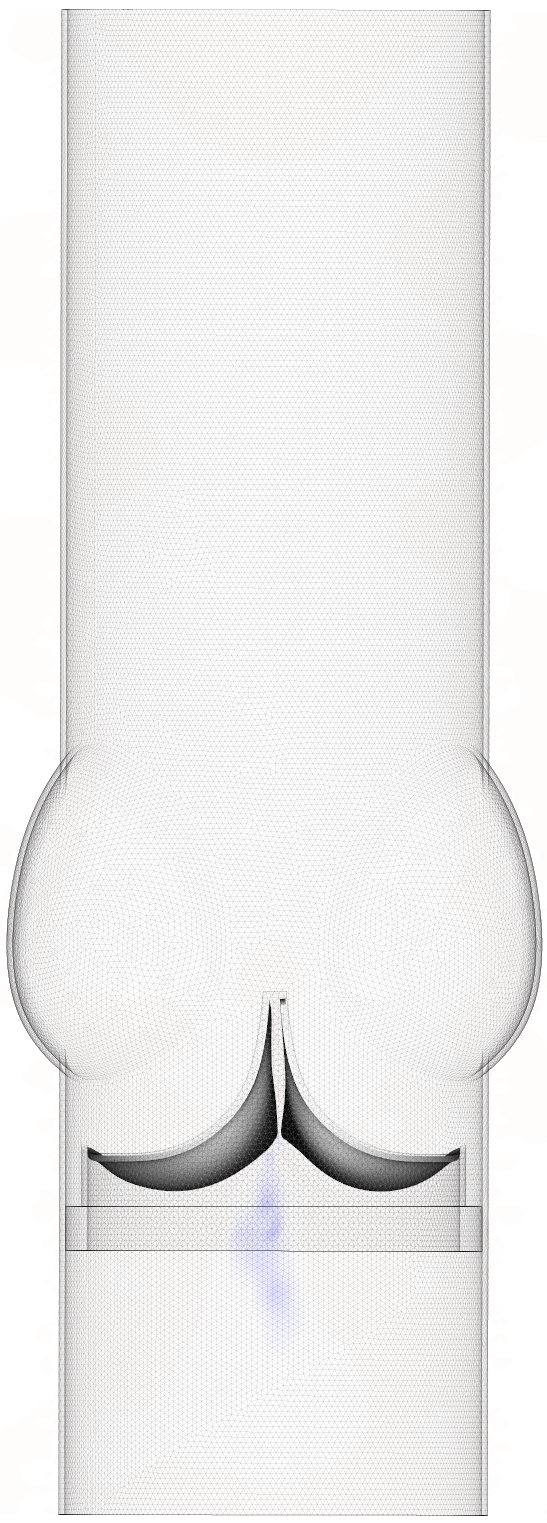}
    \caption{$ t= 0.05$s}    
    \end{subfigure}
    \hfill
    \begin{subfigure}{0.22\textwidth}
            \includegraphics[width=\linewidth]{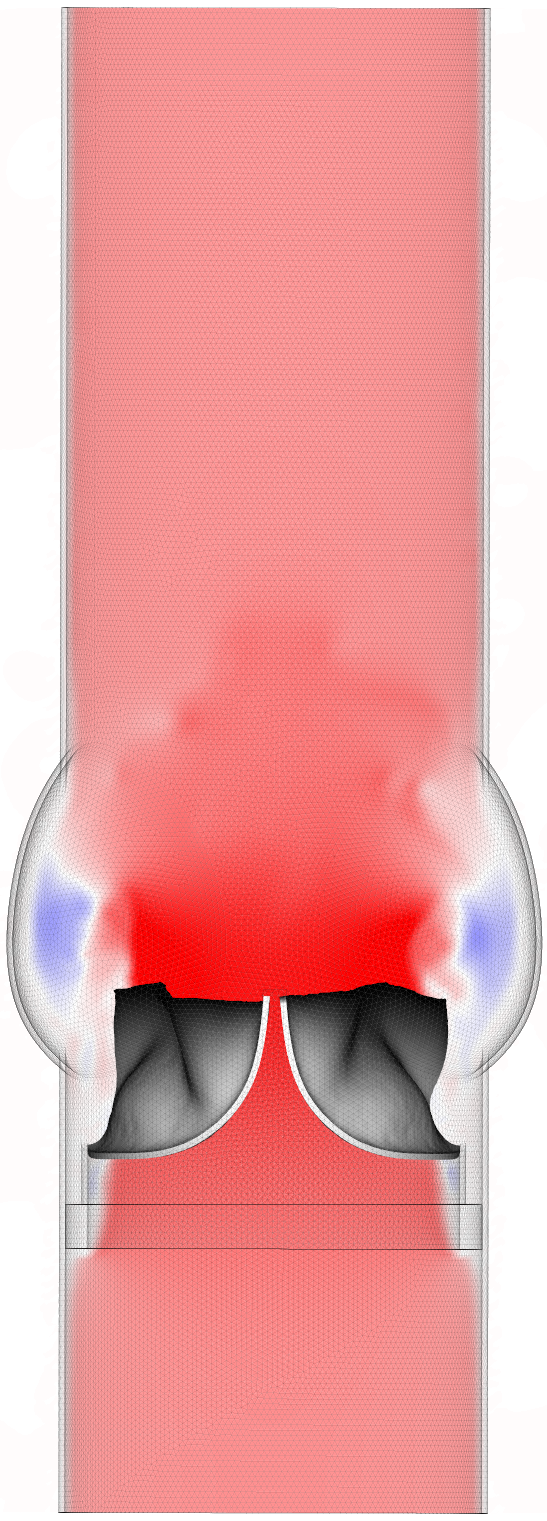}
        \caption{$ t= 0.15$ s}
            
    \end{subfigure}
    \hfill    
    \begin{subfigure}{0.22\textwidth}
        \includegraphics[width=\linewidth]{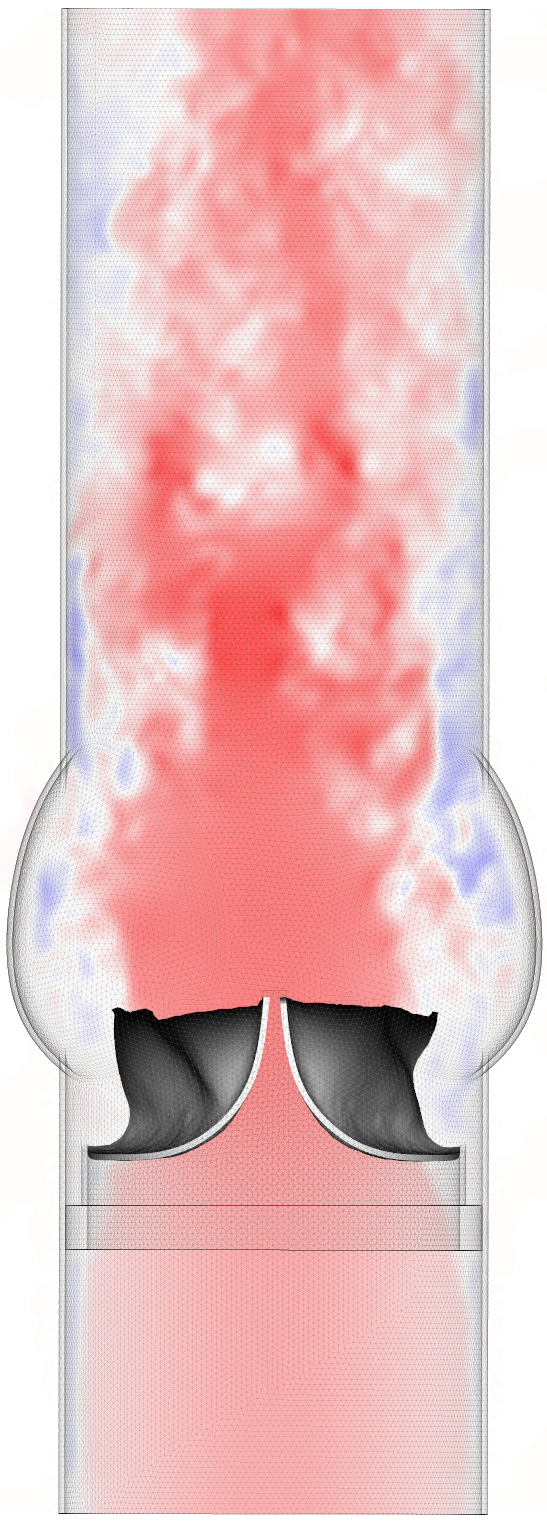}
        \caption{$ t= 0.30$s}    
    \end{subfigure}
    \hfill
    \begin{subfigure}{0.22\textwidth}
        \includegraphics[width=\linewidth]{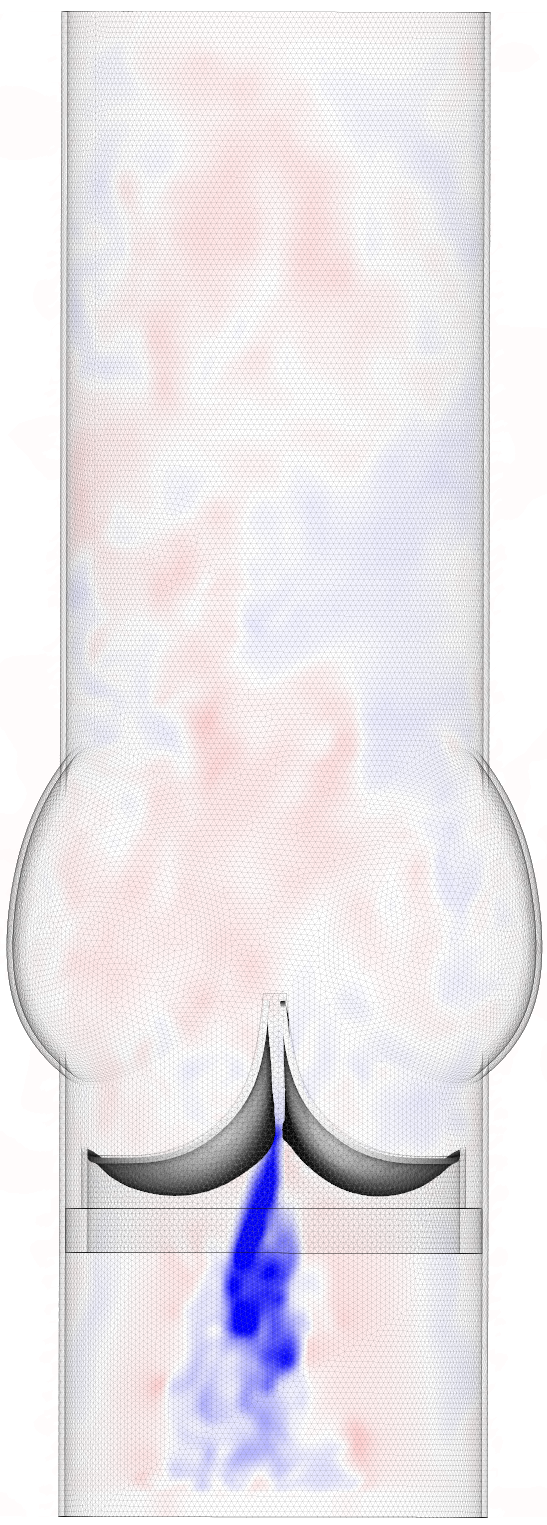}
        \caption{$ t= 0.45$s}    
    \end{subfigure}

    \caption{Representative cross section views of simulated axial velocity for the bioprosthetic heart valve model}
    \label{fig:bhv-flow-slice}
\end{figure}

The boundary conditions include three-element Windkessel models for upstream driving and downstream loading conditions of the aortic test section, as detailed previously \cite{Lee2020}. The inlet and outlet couple these reduced-order models to the detailed flow description through a combination of normal traction and zero tangential velocity boundary conditions. Solid wall boundary conditions are imposed on the remaining domain boundaries.

Fig. \ref{fig:bhv-flow-slice} shows representative cross-sectional views of the flow field during valve opening and closing. During systole, the valve reaches its maximum opening at $t=0.15$ s, corresponding to peak flow velocity. The flow velocity subsequently decreases during diastole ($t=0.30$ s). The flow patterns are similar across all tested kernels, demonstrating that all kernels reasonably capture the physiological dynamics of valve motion.

\begin{figure}[H]
    \centering
    \begin{subfigure}[t]{0.49\textwidth}
        \centering
        \includegraphics[width=\linewidth]{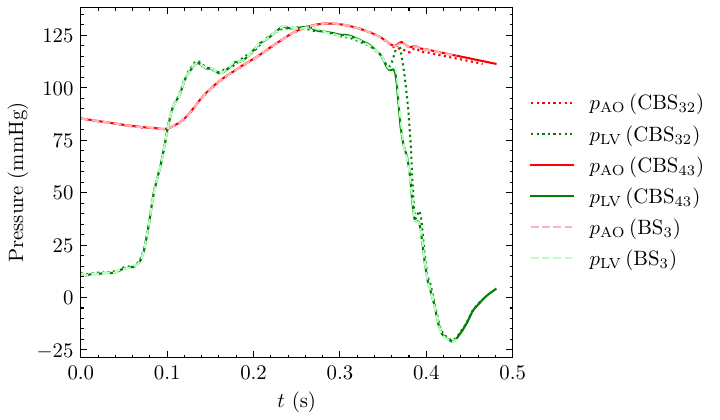}
        \caption{Pressure waveform}
    \end{subfigure}%
    \hfill 
    \begin{subfigure}[t]{0.49\textwidth}
        \centering
        \includegraphics[width=\linewidth]{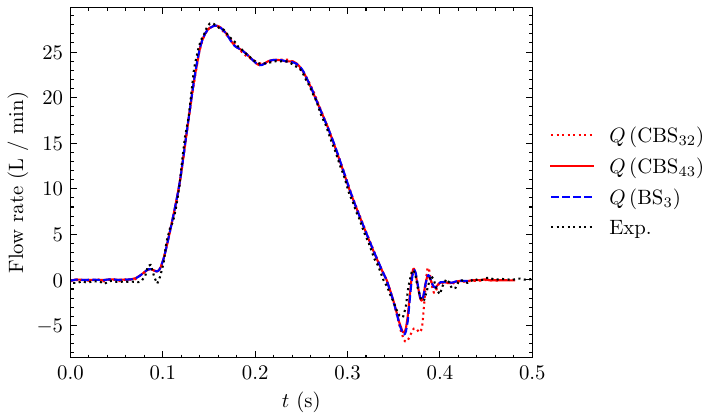}
        \caption{Flow rate waveform}
    \end{subfigure}%
\caption{Comparison of (a) pressure waveforms and (b) flow rates for different kernels in the bioprosthetic heart valve model. All kernels show minimal discrepancy and good agreement with experimental data before valve flutter ($t > 0.35$ s). During early diastole (0.34 s to 0.39 s), CBS$_{32}$ shows notable deviation from other kernels in both pressure and flow rate predictions. The wider CBS$_{43}$ kernel demonstrates superior performance in capturing high-frequency flow characteristics associated with valve motion.}
    \label{fig:bhv-pressure-flowrate-comparison}
\end{figure}

Fig. \ref{fig:bhv-pressure-flowrate-comparison} compares the aortic and left ventricular pressures, along with flow rates obtained using different kernels. Before valve flutter ($t < 0.35$ s), the results show a minimal discrepancy between the kernels, and the predicted flow rates demonstrate good agreement with experimental data. During early diastole when the valve begins to flutter (0.35 s $ < t < 0.39$ s), CBS$_{32}$ shows notable deviation from other kernels and experimental data. This discrepancy can be attributed to CBS$_{32}$'s lower regularity compared to BS$_3$ and CBS$_{43}$, resulting in reduced filtering of noisy components in the interpolated velocity and spread force. The wider CBS$_{43}$ kernel demonstrates superior performance in capturing high-frequency flow characteristics associated with valve motion.

\section{Discussion and Conclusion}  
\label{sec:discussion}  

This study evaluates the performance of CBS kernels in FSI simulations using the IFED method with nodal coupling. Our analysis focuses on volume conservation and accuracy across different mesh resolutions. These kernels are compared with conventional isotropic kernels, including IB and BS kernels, using established FSI benchmark problems. These kernels are compared with conventional isotropic kernels, including IB and BS kernels, using established FSI benchmark problems.

CBS kernels exhibit superior volume conservation performance by maintaining a divergence-free solid velocity field. This leads to reduced volume conservation errors and smaller numerical boundary layer errors compared to isotropic kernels. In the tests considered in this study, CBS kernels achieve satisfactory results without requiring volumetric energy terms or modified invariants in the IFED method, whereas these additional terms are essential for isotropic kernels.

Regarding structural mesh resolution, shear-dominated problems show different optimal configurations for different kernel types. Isotropic kernels perform better with coarser structural meshes relative to the background Cartesian grid. In contrast, CBS kernels perform better with finer structural meshes ($\text{MFAC} \le 1.0$) but show less sensitivity to mesh resolution overall. For problems involving pressure-loads with significant normal forces at the fluid-structure interface, both CBS and isotropic kernels achieve optimal results with structural meshes at least as fine as the background grid.
To address scenarios that combine substantial shear-stresses and pressure-loads, we recommend using structural meshes with resolutions matching or exceeding the background grid, combined with narrower kernels. Based on our findings in this and earlier studies, we recommend using composite kernels over isotropic kernels for FSI simulations due to their superior balance of accuracy and computational efficiency. For most applications, the $\text{CBS}_{32}$ kernel paired with a $\text{MFAC} = 0.5$ value achieves a preferable balance between accuracy and computational cost in the IFED framework.

This study focuses mainly on two-dimensional problems using lower-order elements with MAC grid discretization. Future research should examine the performance of CBS kernels with higher-order elements in three-dimensional problems and explore alternative quadrature schemes beyond nodal coupling scheme investigated here. Additionally, we intend to study more challenging benchmark problems featuring large pressure loads and geometric singularities, which will further test the robustness of these kernels under more extreme conditions.

In conclusion, CBS kernels provide superior performance in volume conservation while offering significant advantages for FSI simulations. By addressing the limitations of this study and exploring new applications, future work can further realize the potential of CBS kernels in complex multi-physics problems.

\section*{Acknowledgments}
L.L. acknowledges funding from  NIH (Awards 5-R01-DK132328-03). C.G. is grateful for support from the Department of Defense (DoD) through the National Defense Science and Engineering Graduate (NDSEG) Fellowship Program. B.E.G. acknowledges funding from the NIH (Awards R01HL157631 and U01HL143336) and NSF (Awards OAC 1450327, OAC 1652541, OAC 1931516, and CBET 1757193).
Simulations were performed using computational facilities provided by the University of North Carolina at Chapel Hill through the Research Computing Division of UNC Information Technology Services.

\bibliographystyle{unsrt}
\bibliography{references}  
\end{document}